%% file: pap.tex
\documentclass[final,leqno,onefignum,onetabnum]{siamltex}

\usepackage{fullpage}
\usepackage{comment}
\usepackage{upgreek}
\usepackage{cite}
\usepackage[mathscr]{eucal}
\usepackage{color}
\usepackage[usenames,dvipsnames,svgnames,table]{xcolor}
\usepackage{amsmath,amssymb,amsopn,mathtools}
\usepackage{graphicx}
\usepackage{hyperref}       
\hypersetup{
	colorlinks=true,   
	citecolor=blue,     
	filecolor=blue,     
	linkcolor=red,    
	urlcolor=blue}      

\usepackage{relsize}
\numberwithin{equation}{section}
\usepackage{enumitem}

\newlist{Assumption}{enumerate}{1}
\setlist[Assumption]{label=A\arabic*}

\definecolor{Blue}{rgb}{0,0,1}
\definecolor{Red}{rgb}{1,0,0}
\definecolor{Green}{rgb}{0,1,0}
\definecolor{Cyan}{rgb}{0,0.72,0.92}
\definecolor{Amethyst}{rgb}{0.6,0.4,0.8}
\definecolor{Bronze}{rgb}{0.8,0.5,0.2}
\definecolor{Violet}{rgb}{0.54,0.17,0.89}

\input{commands}

\usepackage{amsmath,amssymb}
\usepackage{epsfig, algpseudocode, algorithm}
\usepackage{multirow, booktabs}
\usepackage{siunitx}
\usepackage{xr}
\mathtoolsset{showonlyrefs=true}

\newlist{steps}{enumerate}{1}
\setlist[steps, 1]{label = Step \arabic*:}

\title{Local Lagrangian reduced-order modeling for Rayleigh--Taylor instability by solution manifold decomposition}

\author{
  Siu Wun Cheung\thanks{Center for Applied Scientific Computing, Lawrence
  Livermore National Laboratory, Livermore, CA 94550 (cheung26@llnl.gov)}
  \and
  Youngsoo Choi\thanks{Center for Applied Scientific Computing, Lawrence
  Livermore National Laboratory, Livermore, CA 94550 (choi15@llnl.gov)}
  \and
  Dylan Matthew Copeland\thanks{Center for Applied Scientific Computing,
  Lawrence Livermore National Laboratory, Livermore, CA 94550
  (copeland11@llnl.gov)}
  \and
  Kevin Huynh\thanks{Applications, Simulations, and Quality, Lawrence Livermore
  National Laboratory, Livermore, CA 94550 (huynh24@llnl.gov)}
}

\begin{document}
\setlength{\abovedisplayskip}{3pt}
\setlength{\belowdisplayskip}{3pt} 
\setlength{\abovedisplayshortskip}{3pt} 
\setlength{\belowdisplayshortskip}{3pt}

\maketitle

\begin{abstract}
Rayleigh--Taylor instability is a classical hydrodynamic instability of great interest in 
  various disciplines of science and engineering, 
including astrophyics, atmospheric sciences and climate, geophysics, and fusion energy. 
Analytical methods cannot be applied to explain the long-time behavior of Rayleigh--Taylor instability, 
and therefore numerical simulation of the full problem is required. 
However, in order to capture the growth of amplitude of perturbations accurately, 
both the spatial and temporal discretization need to be extremely fine for traditional numerical methods, 
and the long-time simulation may become prohibitively expensive. 
In this paper, we propose efficient reduced order model techniques to accelerate 
the simulation of Rayleigh--Taylor instability in compressible gas dynamics. 
We introduce a general framework for decomposing the solution manifold 
to construct the temporal domain partition and 
temporally-local reduced order model construction with varying Atwood number. 
We propose two practical approaches in this framework, namely  
decomposition by using physical time and by penetration distance respectively. 
Numerical results are presented to examine the performance of the proposed approaches. 
\end{abstract}


\section{Introduction}\label{sec:intro}

Rayleigh--Taylor instability is a classical hydrodynamic instability 
occuring at the interface between two fluids of different densities. 
A denser fluid is initially supported on top by a less dense fluid 
with a counterbalancing pressure gradient under the effects of a gravitational field. 
It was first introduced by Rayleigh \cite{rayleigh1882investigation} in 1882 and 
further investigated by Taylor \cite{taylor1950instability}.  
In recent years, it has been applied to various disciplines of science and engineering, 
including astrophyics \cite{livio1980rayleigh,yamada1991rayleigh,blondin20011rayleigh,ribeyre2004compressible}, 
atmospheric sciences and climate \cite{keskinen1981nonlinear,huang1993nonlinear}, 
geophysics \cite{wilcock1991rayleigh,plag1995rayleigh,conrad1997growth}, 
and fusion energy \cite{betti1998growth,remington2019rayleigh}. 
The readers are referred to \cite{chandrasekhar1961hydrodynamic,sharp1984overview, casner2019icf} 
for a comprehensive overview of Rayleigh--Taylor instability.

In recent years, numerous research efforts have been devoted to modeling the growth of 
Rayleigh--Taylor instability mathematically and numerically, so as to provide better understanding of the nonlinear dynamics. 
Linear stability theory can be used to analytically derive the exponential growth in the amplitude 
in the subwavelength regime in the short-time dynamics \cite{taylor1950instability}.  
However, as the nonlinearity dominates and the dynamics become chaotic in the later stages, 
numerical simulation plays an important role in the studies of the long-time dynamics 
\cite{tryggvason1990computations,youngs1991three,dimonte2004comparative,lee2013numerical,zanella2020two}.  
Depending on the nature of the fluids, 
the physical process is governed by a set of conservation laws 
and mathematically modeled by a system of partial differential equations. 
In this work, we consider the scenario of a single-mode two-dimensional Rayleigh--Taylor instability 
between two compressible and isotropic ideal gases with different densities, in a parametric setting. 
The Euler equation is used to model the compressible gas dynamics 
in a complex multimaterial setting. 
Moreover, we study the effects of Atwood number as a problem parameter on the physical quantities. 
The governing equations are numerically solved in a moving Lagrangian frame, 
where the computational mesh is moved along with the fluid velocity. 
Our work is built on \cite{dobrev2012high}, where a general framework of high-order curvilinear finite
elements and adaptive time stepping of explicit time integrators is proposed for
numerical discretization of the Lagrangian hydrodynamics problem over general
unstructured two-dimensional and three-dimensional computational domains. 
The method shows great capability and lots of advantages, including 
high-order accuracy, total energy conservation, and the Rankine-Hugoniot jump conditions 
fulfilled at shock boundaries. However, forward simulations of Lagrangian hydrodynamics can be computationally very expensive, 
especially in long-time simulation of highly nonlinear problems like Rayleigh--Taylor instability. 
It is therefore desirable to develop efficient techniques for accelerating the computationally expensive simulations. 
In this work, we propose new reduced order model (ROM) techniques for 
Rayleigh--Taylor instability in compressible Euler equations. 

A reduced order model has nowadays become a popular and actively researched 
computational technique to reduce the computational
cost of simulations while minimizing the error introduced in the reduction process. 
Many physics-based models are complex and nonlinear, 
and are formulated on large spatio-temporal domains, 
in which the computational cost can be prohibitively high. It may take a long
time to run one forward simulation even with high performance computing.  In
decision-making applications where multiple forward simulations are needed, such
as parameter study, design optimization \cite{wang2007large, de2020three,
de2018adaptive, white2020dual}, optimal control \cite{choi2015practical,
choi2012simultaneous}, uncertainty quantification \cite{smith2013uncertainty,
biegler2011large}, and inverse problems \cite{galbally2010non,
biegler2011large}, the computationally expensive simulations are not desirable.
To this end, a reduced order model can be useful to obtain
sufficiently accurate approximate solutions with considerable speed-up
compared to a corresponding full order model (FOM). 

Many model reduction schemes have been developed to reduce the computational
cost of simulations while minimizing the error introduced in the reduction
process. Most of these approaches seek to extract an intrinsic solution subspace
for condensed solution representation by a linear combination of reduced basis
vectors. The reduced basis vectors are extracted from performing proper
orthogonal decomposition (POD) on the snapshot data of the FOM simulations.  The
number of degrees of freedom is then reduced by substituting the ROM solution
representation into the (semi-)discretized governing equation.  These approaches
take advantage of both the \textit{known governing equation} and the \textit{solution data}
generated from the corresponding FOM simulations to form linear subspace reduced
order models (LS-ROM).  Example applications include, but are not limited to, the nonlinear
diffusion equations \cite{hoang2021domain, fritzen2018algorithmic}, the Burgers
equation and the Euler equations in small-scale \cite{choi2019space,
choi2020sns, carlberg2018conservative}, the convection--diffusion equations
\cite{mojgani2017lagrangian, kim2021efficientII}, the Navier--Stokes equations
\cite{xiao2014non, burkardt2006pod}, rocket nozzle shape design
\cite{amsallem2015design}, flutter avoidance wing shape optimization
\cite{choi2020gradient}, topology optimization of wind turbine blades
\cite{choi2019accelerating}, lattice structure design
\cite{mcbane2021component}, porous media flow/reservoir simulations
\cite{ghasemi2015localized, jiang2019implementation, yang2016fast,
wang2020generalized}, computational electro-cardiology \cite{yang2017efficient},
inverse problems \cite{fu2018pod}, shallow water equations \cite{zhao2014pod,
cstefuanescu2013pod}, Boltzmann transport problems \cite{choi2021space},
computing electromyography \cite{mordhorst2017pod}, spatio-temporal dynamics of
a predator--prey system \cite{dimitriu2013application}, acoustic wave-driven
microfluidic biochips \cite{antil2012reduced}, and Schr{\"o}dinger equation
\cite{cheng2016reduced}.  Survey papers for the projection-based LS-ROM
techniques can be found in \cite{gugercin2004survey, benner2015survey}. 

In spite of successes of the classical LS-ROM in many applications, these
approaches are limited to the assumption that the intrinsic solution space falls
into a subspace with a small dimension, i.e., the solution space with a
Kolmogorov $n$-width decaying fast.  This assumption is violated in advection-dominated
problems, due to features such as sharp gradients, moving shock fronts, and turbulence, which
hinder these model reduction schemes from being practical.  Our goal in this
paper is to develop an efficient reduced order model for hydrodynamics
simulation with advection-dominated solutions.  Some reduced order model techniques for hydrodynamics or turbulence
models in the literature include \cite{mou2020data, parish2017non,
gadalla2020comparison, bergmann2009enablers, osth2014need, baiges2015reduced,
san2018extreme}, which are mostly built on the Eulerian formulation, i.e., the
computational mesh is stationary with respect to the fluid motion.  In contrast,
numerical methods in the Lagrangian formulation, which are characterized by a
computational mesh that moves along with the fluid velocity, are developed for
better capturing the shocks and preserving the conserved quantities in
advection-dominated problems.  It therefore becomes natural to develop
Lagrangian-based reduced order models to overcome the challenges posed by
advection-dominated problems.  Some existing work in this research direction
include \cite{mojgani2017lagrangian,lu2020lagrangian}, where a Lagrangian POD
and dynamic mode decomposition (DMD) reduced order model are introduced
respectively for the one-dimensional nonlinear advection-diffusion equation.  We
remark that there are some similarities and differences between our work and
\cite{mojgani2017lagrangian} in using POD for developing Lagrangian-based
reduced order models.  It is important to note that our work is based on the
more complicated and challenging two-dimensional compressible Euler equations. 
The reduced bases are built independently for each state variable in our work, while a single basis is built
for the whole state in \cite{mojgani2017lagrangian}.  Furthermore, we introduce
the solution manifold decomposition concept so as to ensure adequate ROM speed-up by 
classifying solution snapshot samples for constructing ROM spaces with low dimension and 
assigning the appropriate ROM to be used in the time marching.

Recently, there have been many attempts to develop efficient ROMs for the
advection-dominated or sharp gradient problems.  The attempts can be divided
mainly into two categories.  The first category enhances the solution
representability of the linear subspace by introducing some special treatments
and adaptive schemes.  A dictionary-based model reduction method for the
approximation of nonlinear hyperbolic equations is developed in
\cite{abgrall2016robust}, where the reduced approximation is obtained from the
minimization of the residual in the $L_1$ norm for the reduced linear subspace.
In \cite{carlberg2015adaptive}, a fail-safe $h$-adaptive algorithm is developed.
The algorithm enables ROMs to be incrementally refined to capture the shock
phenomena which are unobserved in the original reduced basis through
a-posteriori online enrichment of the reduced-basis space by decomposing a given
basis vector into several vectors with disjoint support.  The windowed
least-squares Petrov--Galerkin model reduction for dynamical systems with
implicit time integrators is introduced in \cite{parish2019windowed,
shimizu2020windowed}, which can overcome the challenges arising from the
advection-dominated problems by representing only a small time window with a
local ROM.  Our previous work \cite{copeland2022reduced} adopts a similar approach 
for Lagrangian hydrodynamics. 
Another active research direction is to exploit the sharp gradients
and represent spatially local features in ROM, such as the online adaptivity
bases and adaptive sampling approach \cite{peherstorfer2018model} and the shock
reconstruction surrogate approach \cite{constantine2012reduced}.  In
\cite{taddei2020space}, an adaptive space-time registration-based model
reduction is used to align local features of parameterized hyperbolic PDEs in a
fixed one-dimensional reference domain.  Some new approaches have been developed
for aligning the sharp gradients by using a superposition of snapshots with
shifts or transforms.  In \cite{reiss2018shifted}, the shifted proper orthogonal
decomposition (sPOD) introduces time-dependent shifts of the snapshot matrix in
POD in an attempt to separate different transport velocities in advection-dominated
problems.  The practicality of this approach relies heavily on accurate
determination of shifted velocities.  In \cite{rim2018transport}, an iterative
transport reversal algorithm is proposed to decompose the snapshot matrix into
multiple shifting profiles.  In \cite{welper2020transformed},  inspired by the
template fitting \cite{kirby1992reconstructing}, a high resolution transformed
snapshot interpolation with an appropriate behavior near singularities is
considered. 

The second category replaces the linear subspace solution representation with
the nonlinear manifold, which is a very active research direction.  Recently, a
neural network-based reduced order model is developed in \cite{lee2020model} and
extended to preserve the conserved quantities in the physical conservation laws
\cite{lee2019deep}.  In these approaches, the weights and biases in the neural
network are determined in the training phase, and existing numerical methods,
such as finite difference and finite element methods, are utilized.  However,
since the nonlinear terms need to be updated every time step or Newton step, and
the computation of the nonlinear terms still scale with the FOM size, these
approaches do not achieve any speed-up with respect to the corresponding FOM.
Recently, Kim, et al., have achieved a considerable speed-up with the nonlinear
manifold reduced order model \cite{kim2021fast, kim2020efficient}, but it was
only applied to small problems.  Manifold approximations via transported
subspaces in \cite{rim2019manifold} introduced a low-rank approximation to the
transport dynamics by approximating the solution manifold with a transported
subspace generated by low-rank transport modes.  However, their approach is
limited to a one-dimensional problem setting.  In \cite{rim2020depth}, a depth
separation approach for reduced deep networks in nonlinear model reduction is
presented, in which the reduced order model is composed with hidden layers with
low-rank representation.

The method presented in this paper belongs to the first category. 
The underlying concept is to build small and accurate projection-based reduced-order models 
by decomposing the solution manifold into submanifolds.
These reduced-order models are local in the sense that 
each of them will be valid only over a certain sub-interval of the temporal domain. 
The appropriate local reduced order model is chosen based on the current state of the system,  
and all the local reduced order models will cover the whole time marching in the online phase. 
The concept of a local reduced order model was introduced in \cite{washabaugh2012nonlinear, 
amsallem2012nonlinear}, where unsupervised clustering is used for the solution manifold decomposition. 
In our work, the solution manifold decomposition is based on a suitably defined physics-based indicator,
which generalizes the time-windowing ROM approach for Lagrangian hydrodynamics simulations
in our previous work \cite{copeland2022reduced}. 
The idea of time-windowing ROM is to construct temporally-local ROM spaces which are small but accurate within a
short period in advection-dominated problems. 
The method provided tremendous speed-up and accurate approximated solutions in various benchmark experiments 
including Sedov blast, Gresho vortices, Taylor-Green vortices, and triple-point problems. 
However, in order to approximate the solutions of Rayleigh--Taylor instability over a range of Atwood number, 
the time-windowing approach is insufficient to provide good approximations, 
since the penetration speed of the fluid interface varies with the Atwood number. 
In this work, we modify our ROM approach by introducing new techniques 
to capture the fast-moving shock boundaries, which provide significant 
improvements in Rayleigh--Taylor instability.

Our contribution in this paper is to propose an efficient model order reduction scheme 
for accelerating the simulation of Rayleigh--Taylor instability in compressible gas dynamics with varying Atwood number. 
Similar to \cite{copeland2022reduced}, 
the idea is to construct temporally-local ROM spaces which are small but accurate within a
short period in advection-dominated problems to achieve a good speed-up and solution accuracy. 
Proper orthogonal decomposition (POD) is used to extract the dominant modes in solution
representability, the solution nonlinear subspace (SNS) method is used to establish the subspace relations and 
construct the nonlinear term bases as in \cite{choi2020sns}, 
and an oversampling discrete empirical
interpolation method (DEIM) serves as a hyper-reduction technique to
reduce the complexity due to the nonlinear terms in the governing equations.  
The difference between this work and \cite{copeland2022reduced} is the introduction of 
a new temporal domain partition scheme that allows a smaller dimension of reduced basis in each temporal subdomain and 
use of temporally-local ROMs at different Atwood numbers for better solution accuracy. 
In this work, we propose a general framework to decompose the solution manifold by 
a suitably defined indicator. 
The indicator is used to classify solution snapshot samples for constructing 
ROM spaces with low dimension and to assign the appropriate ROM to be used in the time marching. 
We consider two such indicators, namely the physical time and the penetration distance, 
and compare the performance in speed-up and solution accuracy.
As we will see in our numerical experiments, the penetration distance is a good alternative indicator 
to resolve the deficiency of degenerating solution accuracy with respect to Atwood number observed in the time-windowing approach. 
We remark that, by taking the advantage of Lagrangian formulation of the Euler equations,
the penetration distance is easily accessible and plays an important role in the newly proposed distance-windowing approach. 

The rest of the paper is organized as follows. 
In Section~\ref{sec:FOM}, we introduce the governing equations and the 
numerical discretization which will be used as the full order model. 
Next, a projection-based ROM is described in Section~\ref{sec:ROM}. 
In Section~\ref{sec:decompose}, a general framework of solution manifold decomposition is introduced,  
and two practical examples will be discussed.
Numerical results are presented in Section~\ref{sec:numerical}. 
Finally, conclusions are given in Section~\ref{sec:conclusion}. 

\section{Numerical modeling}\label{sec:FOM}
In this section, we describe a direct numerical simulation methodology
for Rayleigh--Taylor instability in compressible gas dynamics, using a finite element method, 
which will serve as the full order model in this work. 
We first present the governing equations,
initial and boundary conditions which give rise to Rayleigh--Taylor instability.
We also present the key ingredients of the numerical solver.

\subsection{Governing equations}\label{sec:prob}
We consider the two-dimensional mathematical model of Rayleigh--Taylor
hydrodynamic instability caused by a gravitational field acting on 
compressible isotropic ideal gas with stratified densities. 
Suppose $\finalTime > 0$ is the final time of the hydrodynamic process. 
For $t \in [0,\finalTime]$, let $\solDomainSymbol(\timeSymbol) \subset \mathbb{R}^2$ 
denote a continuous medium which  
deforms in time, with the initial configuration at $\timeSymbol = 0$
being a rectangular domain $\initialDomain = \solDomainSymbol(0) = [0, 1/2] \times [-1, 1]$ with aspect ratio 1:4. 
The compressible gas dynamics is modeled by the Euler equations in a Lagrangian
reference frame \cite{harlow1971fluid}, driven by a constant gravitational acceleration $\gravitySymbol$: 
\begin{equation}\label{eq:euler}
  \begin{aligned}
    \text{momentum conservation}:& & \densitySymbol\frac{d\velocitySymbol}{d\timeSymbol} &=
    \gradientSymbol \cdot \stressSymbol + \gravitySymbol \\
    \text{mass conservation}:& & \dfrac{1}{\densitySymbol}\frac{d\densitySymbol}{d\timeSymbol} &=
    -\gradientSymbol \cdot \velocitySymbol \\
   \text{energy conservation}:& & \densitySymbol\frac{d\energySymbol}{d\timeSymbol} &=
    \stressSymbol :  \gradientSymbol \velocitySymbol \\
    \text{equation of motion}:& & \frac{d\positionSymbol}{d\timeSymbol} &= \velocitySymbol.
  \end{aligned}
\end{equation}
Here, $\frac{d}{d\timeSymbol} = \frac{\partial}{\partial \timeSymbol} + 
\velocitySymbol \cdot \gradientSymbol$ is the material derivative, 
$\densitySymbol$ denotes the density of the fluid, 
$\stressSymbol$ denotes the deformation stress tensor, 
$\positionSymbol$ denotes the Lagrangian position, 
$\velocitySymbol$ denotes the velocity, 
and $\energySymbol$ denotes the internal energy per unit mass.  
Each of these physical quantities depends on the time $\timeSymbol \in [0,\finalTime]$ and 
the Eulerian coordinates of the particle $\initialPosition \in \initialDomain$. 
In gas dynamics, the stress tensor is isotropic, and we write $\stressSymbol = -\pressureSymbol
\identitySymbol + \artificialStressSymbol$, where $\pressureSymbol$ denotes the
thermodynamic pressure, and $\artificialStressSymbol$ denotes the artificial
viscosity stress. Due to the ideal gas assumption, the thermodynamic pressure is 
related to the density and the internal energy by the equation of state
\begin{equation}\label{eq:EOS}
  \pressureSymbol = (\adiabaticIndexSymbol - 1) \densitySymbol \energySymbol, 
\end{equation}
where the adiabatic index is $\adiabaticIndexSymbol = 5/3$. 
For $\timeSymbol > 0$, the system is prescribed with a boundary condition
$\velocitySymbol \cdot \normalSymbol = 0$, where $\normalSymbol$ is
the outward normal unit vector on the domain boundary. 
In order to close the system, an initial condition needs to be imposed. 
We assume the gas is initially stratified at the interface $\initialPosition_2 = 0$, 
supported by a counterbalancing pressure gradient and perturbed by a velocity. 
In mathematical terms, at the time $\timeSymbol = 0$, 
for $\initialPosition = (\initialPosition_1, \initialPosition_2) \in \initialDomain$, we take 
\begin{equation}
\begin{split}
\velocitySymbol(0,\initialPosition) & = 
\left(0, 0.02 \cos(2 \pi \initialPosition_1) \exp(-2 \pi \initialPosition_2^2)\right) \\
\densitySymbol(0,\initialPosition) & = 
\begin{cases}
R& \text{ if } \initialPosition_2 \geq 0 \\
1 & \text{ if } \initialPosition_2 < 0
\end{cases} \\
\energySymbol(0,\initialPosition) & = 
\dfrac{4 + R- \densitySymbol \initialPosition_2}{(\adiabaticIndexSymbol - 1)\densitySymbol}, 
\end{split}
\label{eq:ic}
\end{equation}
where $R> 1$ is the density ratio of the fluids. 
In the context of hydrodynamic instability, we define the dimensionless Atwood number by 
\begin{equation}
\param = \dfrac{R - 1}{R + 1}, 
\end{equation}
which increases mildly with the density ratio $R$. 
The density ratio and the Atwood number play an important role in the perturbation of the fluid interface. 
In our work, the mathematical model is considered in the parametric setting 
by regarding the Atwood number $\param$ as 
a problem parameter which lies in a parametric domain $\paramDomain \subset \mathbb{R}$. 
In this way, the physical quantities are treated as functions 
not only of the time $\timeSymbol \in [0,\finalTime]$ and 
the Eulerian coordinates of the particle $\initialPosition \in \initialDomain$, 
but also the problem parameter $\param \in \paramDomain$. 

\subsection{Finite element modeling}\label{sec:FEM}
In \cite{dobrev2012high}, a general framework of high-order curvilinear finite
elements and adaptive time stepping of explicit time integrators is proposed for
numerical discretization of the Lagrangian hydrodynamics problem over general
unstructured two-dimensional and three-dimensional computational domains.  Using
general high-order polynomial basis functions for approximating the state
variables, and curvilinear meshes for capturing the geometry of the flow and
maintaining robustness with respect to mesh motion, the method achieves
high-order accuracy.  On the other hand, a modification is made to the
second-order Runge-Kutta method to compensate for the lack of total energy
conservation in standard high-order time integration techniques.  The
introduction of an artificial viscosity tensor further generates the appropriate
entropy and ensures the Rankine-Hugoniot jump conditions at shock boundaries.

Following \cite{dobrev2012high}, we adopt a spatial discretization for
\eqref{eq:euler} using a kinematic space $\kinematicFE \subset
[H^1(\initialDomain)]^\dimensionSymbol$ for approximating the position and the
velocity, and a thermodynamic space $\thermodynamicFE \subset
L_2(\initialDomain)$ for approximating the energy. The density can be
eliminated, and the equation of mass conservation can be decoupled from
\eqref{eq:euler}. We assume high-order finite element (FE) discretization in
space, and the finite dimensions $\sizeKinematicFE$ and
$\sizeThermodynamicFE$ are the global numbers of FE degrees of freedom in the
corresponding discrete FE spaces.  For more details, see \cite{dobrev2012high}.
The FE coefficient vector functions for velocity and position are denoted as
$\velocity, \position:[0,\finalTime]\times \paramDomain \mapto
\RR{\sizeKinematicFE}$,
and the coefficient vector function for energy is denoted as
$\energy:[0,\finalTime] \times \paramDomain \mapto \RR{\sizeThermodynamicFE}$.
The semidiscrete Lagrangian conservation laws can be expressed as a nonlinear
system of differential equations in the coefficients with respect to the bases
for the kinematic and thermodynamic spaces:
\begin{equation}\label{eq:fom}
  \begin{aligned}
    \text{momentum conservation}:& & \kinematicMassMat\frac{d\velocity}{d\timeSymbol} &=
    -\forceMat(\velocity, \energy, \position; \param) \cdot\oneVec +  \kinematicMassMat \mathbf{g} \\
    \text{energy conservation}:& & \thermodynamicMassMat\frac{d\energy}{d\timeSymbol} &=
    \forceMat(\velocity, \energy, \position; \param)^T\cdot\velocity \\
    \text{equation of motion}:& & \frac{d\position}{d\timeSymbol} &= \velocity, 
  \end{aligned}
\end{equation}
where $\kinematicMassMat \mathbf{g}$ denotes the effects of the gravitation force in the discrete system.  

Let $\state \equiv (\velocity; \energy; \position)^T\in\RR{\sizeWholeFE}$,
$\sizeWholeFE = 2\sizeKinematicFE + \sizeThermodynamicFE$, be the hydrodynamic
state vector. Then the semidiscrete conservation equation of \eqref{eq:fom} can
be written in a compact form as
\begin{equation}\label{eq:combinedFOM}
  \frac{d\state}{d\timeSymbol} = \forceSys(\state;\param),  
\end{equation}
where the nonlinear force term, $\forceSys:\RR{\sizeWholeFE} \times
\paramDomain \mapto \RR{\sizeWholeFE} $, is defined as 
\begin{align}\label{eq:combinedForce}
  \forceSys(\state;\param) &\equiv
  \pmat{\forceSysVelocity(\velocity, \energy, \position; \param) \\
        \forceSysEnergy(\velocity, \energy, \position; \param) \\
        \forceSysPosition(\velocity, \energy, \position; \param)} 
        \equiv
  \pmat{-\kinematicMassMat^{-1}\forceOne(\state;\param) + \mathbf{g} \\
        \thermodynamicMassMat^{-1}\forceTv(\state;\param) \\
        \velocity},
\end{align}
where $\forceOne: \RR{\sizeWholeFE} \times 
\paramDomain \mapto \RR{\sizeKinematicFE}$ and $\forceTv: \RR{\sizeWholeFE}
\times \paramDomain \mapto \RR{\sizeThermodynamicFE}$ are
nonlinear vector functions that are defined respectively as 
\begin{align}\label{eq:forceOneforceTv}
  \forceOne(\state;\param) &\equiv \forceMat(\velocity, \energy, \position; \param)\cdot\oneVec, & 
  \forceTv(\state;\param) &\equiv \forceMat(\velocity, \energy, \position; \param)^T\cdot\velocity.
\end{align}
From now on, we drop the dependence on the parameter $\param$ to simplify the notations 
where there is no ambiguity. 

In order to obtain a fully discretized system of equations, one needs to apply a
time integrator. We consider an explicit Runge-Kutta scheme called the RK2-average scheme,
which is proved to conserve the discrete total energy exactly (see Proposition 7.1 of
\cite{dobrev2012high}). 
The temporal domain is discretized as $\{
  \timek{\timeIndex} \}_{\timeIndex=0}^{\ntimestep}$, where $\timek{n}$ denotes
a discrete moment in time with $\timek{0} = 0$, $\timek{\ntimestep} =
\finalTime$, and $\timek{n-1} < \timek{n}$ for
$\timeIndex\in\nat{\ntimestep}$, where $\nat{N}\equiv \{1,\ldots,N\}$.  The
computational domain at time $\timek{\timeIndex}$ is denoted as
$\Omega^{\timeIndex} \equiv \Omega(\timek{\timeIndex})$. We denote the
quantities of interest defined on $\Omega^{\timeIndex}$ with a subscript
$\timeIndex$. Starting with $\velocityt{0} = \velocity(0)$, $\energyt{0} = \energy(0)$, 
and $\positiont{0} = \position(0)$, the discrete state is updated by 
\begin{align}\label{eq:RK2-avg}
  \velocityt{\timeIndex+\frac{1}{2}} &= \velocityt{\timeIndex} + (\timestepk{\timeIndex}/2)
    (-\kinematicMassMat^{-1} \forceOnek{\timeIndex} + \mathbf{g}), &
    \velocityt{\timeIndex+1} &= \velocityt{\timeIndex} + \timestepk{\timeIndex}
    (-\kinematicMassMat^{-1} \forceOnek{\timeIndex+\frac{1}{2}} + \mathbf{g}), \\
  \energyt{\timeIndex+\frac{1}{2}} &= \energyt{\timeIndex} + (\timestepk{\timeIndex}/2)
    \thermodynamicMassMat^{-1} \forceTvk{\timeIndex}, &
    \energyt{\timeIndex+1} &= \energyt{\timeIndex} + \timestepk{\timeIndex}
    \thermodynamicMassMat^{-1} \avgforceTvk{\timeIndex+\frac{1}{2}}, \\
  \positiont{\timeIndex+\frac{1}{2}} &= \positiont{\timeIndex} + (\timestepk{\timeIndex}/2)
    \velocityt{\timeIndex+\frac{1}{2}}, & \positiont{\timeIndex+1} &=
    \positiont{\timeIndex} + \timestepk{\timeIndex} \avgvelocityt{\timeIndex+\frac{1}{2}},
\end{align}
where the state $\statet{\timeIndex} 
= (\velocityt{\timeIndex}; \energyt{\timeIndex}; 
\positiont{\timeIndex})^T \in\RR{\sizeWholeFE}$ 
is used to compute the updates 
\begin{align}
\forceOnek{\timeIndex} & = \left(\forceMat (\statet{\timeIndex}) \right ) \cdot\oneVec, &
\forceTvk{\timeIndex} & = \left(\forceMat (\statet{\timeIndex}) \right )^T \cdot 
\velocityt{\timeIndex+\frac{1}{2}},
\end{align}
in the first stage. Similarly, $\statet{\timeIndex+\frac{1}{2}} = 
(\velocityt{\timeIndex+\frac{1}{2}}; \energyt{\timeIndex+\frac{1}{2}}; 
\positiont{\timeIndex+\frac{1}{2}})^T \in\RR{\sizeWholeFE}$ 
is used to compute the updates 
\begin{align}
\forceOnek{\timeIndex+\frac{1}{2}} 
& = \left(\forceMat (\statet{\timeIndex+\frac{1}{2}}) \right ) \cdot \oneVec, &
\avgforceTvk{\timeIndex+\frac{1}{2}} & = \left
(\forceMat (\statet{\timeIndex + \frac{1}{2}}) \right )^T \cdot
\avgvelocityt{\timeIndex+\frac{1}{2}},
\end{align}
with $\avgvelocityt{\timeIndex+\frac{1}{2}} = (\velocityt{\timeIndex} +
\velocityt{\timeIndex+1})/2$ in the second stage.  Note that the RK2-average
scheme is different from the midpoint RK2 scheme in the updates for energy and
position.  The RK2-average scheme uses the midpoint velocity
$\velocityt{\timeIndex+\frac{1}{2}}$ and the average velocity
$\avgvelocityt{\timeIndex+\frac{1}{2}}$ to update energy and position in the
first stage and the second stage respectively, while the midpoint RK2 uses the
initial velocity $\velocityt{\timeIndex}$ and the midpoint velocity
$\velocityt{\timeIndex+\frac{1}{2}}$. 
Since an explicit Runge-Kutta method is used, we need to control the time step
size in order to maintain the stability of the fully discrete scheme.  We
follow the automatic time step control algorithm described in Section 7.3 of
\cite{dobrev2012high}, where the time step size is controlled by estimates at
all quadrature points in the mesh used in the evaluation of the force matrix $\forceMat$.

\section{Reduced order model}\label{sec:ROM}

In this section, we present the details of the projection-based reduced oder
model for the semi-discrete Lagrangian conservation laws \eqref{eq:fom}. 
The reduced order model is constructed in the offline phase 
and deployed in the online phase. 
In what follows, we first discuss all the essential ingredients of the reduced order model, 
and move on to discuss the construction. 

\subsection{ROM simulation}
In the reduced order model, 
we restrict our solution space to a subspace spanned by a reduced basis for each
field. That is, the subspaces for velocity, energy, and position fields are
defined as
\begin{align}\label{eq:subspces}
  \velocitySubspace &\equiv
  \Span{\velocityBasisVeck{\basisIndex}}_{\basisIndex=1}^{\sizeROMvelocity} \subseteq
  \RR{\sizeKinematicFE},&
  \energySubspace &\equiv
  \Span{\energyBasisVeck{\basisIndex}}_{\basisIndex=1}^{\sizeROMenergy} \subseteq
  \RR{\sizeThermodynamicFE},&
  \positionSubspace &\equiv
  \Span{\positionBasisVeck{\basisIndex}}_{\basisIndex=1}^{\sizeROMposition} \subseteq
  \RR{\sizeKinematicFE}, 
\end{align}
with $\dim(\velocitySubspace)=\sizeROMvelocity\ll\sizeKinematicFE$,
$\dim(\energySubspace)=\sizeROMenergy\ll\sizeThermodynamicFE$,  and
$\dim(\positionSubspace)=\sizeROMposition\ll\sizeKinematicFE$. Using these
subspaces, each discrete field is approximated in trial subspaces, $\velocity
\approx \velocityApprox \in \velocityOS + \velocitySubspace$, $\energy \approx
\energyApprox \in \energyOS + \energySubspace$, and $\position \approx
\positionApprox \in \positionOS + \positionSubspace$. 
For a generic problem parameter $\param \in \paramDomain$, we write 
\begin{align}\label{eq:solrepresentation}
  \velocityApprox(\timeSymbol; \param) &= \velocityOS(\param) + \velocityBasis\ROMvelocity(\timeSymbol; \param), \\
  \energyApprox(\timeSymbol; \param) &= \energyOS(\param) + \energyBasis\ROMenergy(\timeSymbol; \param), \\
  \positionApprox(\timeSymbol; \param) &= \positionOS(\param) + \positionBasis\ROMposition(\timeSymbol; \param), 
\end{align}
where 
$\velocityOS(\param)\in\RR{\sizeKinematicFE}$,
$\energyOS(\param)\in\RR{\sizeThermodynamicFE}$, and
$\positionOS(\param)\in\RR{\sizeKinematicFE}$ denote the prescribed offset
vectors for velocity, energy, and position fields respectively; the orthonormal basis
matrices
$\velocityBasis\in\RR{\sizeKinematicFE\times\sizeROMvelocity}$,
$\energyBasis\in\RR{\sizeThermodynamicFE\times\sizeROMenergy}$, and
$\positionBasis\in\RR{\sizeKinematicFE\times\sizeROMposition}$
are defined as
\begin{align}\label{eq:basisMats}
  \velocityBasis &\equiv \bmat{\velocityBasisVeck{1} & \cdots & \velocityBasisVeck{\sizeROMvelocity}},&
  \energyBasis &\equiv \bmat{\energyBasisVeck{1} & \cdots & \energyBasisVeck{\sizeROMenergy}},&
  \positionBasis &\equiv \bmat{\positionBasisVeck{1} & \cdots & \positionBasisVeck{\sizeROMposition}};
\end{align}
and $\ROMvelocity:[0,\finalTime]\times \paramDomain \mapto
\RR{\sizeROMvelocity}$, $\ROMenergy:[0,\finalTime]\times \paramDomain \mapto
\RR{\sizeROMenergy}$, and $\ROMposition:[0,\finalTime]\times \paramDomain \mapto
\RR{\sizeROMposition}$ denote the time-dependent generalized coordinates for
velocity, energy, and position fields, respectively.  One natural choice of the
offset vectors is to use the initial values, i.e.  $\velocityOS(\param) = \velocity(0;
\param)$, $\energyOS(\param) = \energy(0; \param)$, and $\positionOS(\param) =
\position(0; \param)$. 

The nonlinear matrix function, $\forceMat$, changes every time the state
variables evolve.  Additionally, it needs to be multiplied by the
basis matrices whenever updates in the nonlinear term occur, which scales
with the FOM size. Therefore, we cannot expect any speed-up
without special treatment of the nonlinear terms.  To overcome this issue,
a hyper-reduction technique needs to be applied (cf. \cite{chaturantabut2010nonlinear}),
where $\forceOne$ and $\forceTv$ are approximated as 
 \begin{align}\label{eq:DEIM_approx}
   \forceOne &\approx \forceOneBasis \ROMforceOne,& 
   \forceTv  &\approx \forceTvBasis \ROMforceTv.
 \end{align}
That is, $\forceOne$ and $\forceTv$ are projected onto subspaces
$\forceOneSubspace \equiv
\Span{\forceOneBasisVeck{\basisIndex}}_{\basisIndex=1}^{\sizeROMforceOne}$ and
$\forceTvSubspace \equiv
\Span{\forceTvBasisVeck{\basisIndex}}_{\basisIndex=1}^{\sizeROMforceTv}$, where
$\forceOneBasis \equiv \bmat{\forceOneBasisVeck{1} & \ldots &
\forceOneBasisVeck{\sizeROMforceOne}} \in
\RR{\sizeKinematicFE\times\sizeROMforceOne}$, $\sizeROMforceOne \ll
\sizeKinematicFE$ and $\forceTvBasis \equiv \bmat{\forceTvBasisVeck{1} & \ldots
& \forceTvBasisVeck{\sizeROMforceTv}} \in
\RR{\sizeThermodynamicFE\times\sizeROMforceTv}$, $\sizeROMforceTv \ll
\sizeThermodynamicFE$, denote the nonlinear term basis matrices, and $\ROMforceOne
\in \RR{\sizeROMforceOne}$ and $\ROMforceTv \in \RR{\sizeROMforceTv}$ denote the
generalized coordinates of the nonlinear terms. 
Now we show how the generalized coordinates, $\ROMforceOne$, can be
determined by the following interpolation:
 \begin{align}\label{eq:DEIM_interpolation}
   \forceOneSamplingMat^T \forceOne = \forceOneSamplingMat^T
   \forceOneBasis\ROMforceOne,
 \end{align}
where $\forceOneSamplingMat \equiv \bmat{\unitveck{p_1},
\ldots,\unitveck{p_{\sizeROMforceOneSample}}} \in \RR{\sizeKinematicFE \times
\sizeROMforceOneSample}$, $\sizeROMforceOne \leq \sizeROMforceOneSample \ll
\sizeKinematicFE$, is the sampling matrix and $\unitveck{p_i}$ is the $p_i$-th
column vector of the identity matrix
$\identityMat{\sizeKinematicFE}\in\RR{\sizeKinematicFE\times\sizeKinematicFE}$.
Note that Eq.~\eqref{eq:DEIM_interpolation} is an over-determined system.  Thus,
we solve the least-squares problem, i.e.,
 \begin{align}\label{eq:ls-hyper}
   \ROMforceOne &= \argmin_{\dummyVec\in\RR{\sizeROMforceOne}} \|
   \forceOneSamplingMat^T ( \forceOne - \forceOneBasis\dummyVec ) \|_2^2.    
 \end{align}
The solution to the least-squares problem~\eqref{eq:ls-hyper} is 
 \begin{align}\label{eq:sol-ls-hyper}
   \ROMforceOne = (\forceOneSamplingMat^T\forceOneBasis)^{\dagger}\forceOneSamplingMat^T  \forceOne,
 \end{align}
where the Moore--Penrose inverse of a matrix $\dummyMat \in \RR{I \times J}$,
$I\geq J$, with full column rank is defined as $\dummyMat^{\dagger} :=
(\dummyMat^T\dummyMat)^{-1}\dummyMat^T$. 
Instead of constructing the sampling matrix
$\forceOneSamplingMat$, for efficiency we simply store the sampling indices
$\{p_1,\ldots,p_{\sizeROMforceOneSample}\} \subset \nat{\sizeKinematicFE}$.
More precisely, the reduced matrix $\left(\forceOneSamplingMat^T\forceOneBasis\right)^{\dagger}
  \in \RR{\sizeROMforceOne\times\sizeROMforceOneSample}$
 can be precomputed and stored in the offline phase,
 and is multiplied to the sampled entries
 $\forceOneSamplingMat^T \forceOne \in \RR{\sizeROMforceOneSample}$
 to obtain $\ROMforceOne$ by \eqref{eq:sol-ls-hyper} in the online phase. 
Similarly, following the same procedure 
to compute the generalized coordinates of the nonlinear term of the energy
conservation equation, we have
 \begin{align}\label{eq:sol-ls-hyper-2}
   \ROMforceTv = (\forceTvSamplingMat^T\forceTvBasis)^{\dagger}\forceTvSamplingMat^T  \forceTv,
 \end{align}
where we denote the sampling matrix by $\forceTvSamplingMat \in
\RR{\sizeThermodynamicFE \times \sizeROMforceTvSample}$. 

Furthermore, we apply the solution nonlinear subspace (SNS)
method in \cite{choi2020sns} to establish the subspace relations $\forceOneBasis =
\kinematicMassMat \velocityBasis$ and $\forceTvBasis = \thermodynamicMassMat
\energyBasis$. 
Using the SNS relation together with the above reduced approximations, 
the hyper-reduced system of \eqref{eq:fom} can be written as :
\begin{align}\label{eq:rom-hr-sns}
   \frac{d\ROMvelocity}{d\timeSymbol} &=
    -\ROMforceOne(\velocityOS + \velocityBasis\ROMvelocity, \energyOS + \energyBasis\ROMenergy, \positionOS + \positionBasis\ROMposition, \timeSymbol; \param) + \velocityBasis^T \mathbf{g} \\
    \frac{d\ROMenergy}{d\timeSymbol} &=
    \ROMforceTv(\velocityOS + \velocityBasis\ROMvelocity, \energyOS + \energyBasis\ROMenergy, \positionOS + \positionBasis\ROMposition, \timeSymbol; \param) \\
    \frac{d\ROMposition}{d\timeSymbol} &= \positionBasis^T \velocityOS +
    \positionBasis^T \velocityBasis \ROMvelocity.
\end{align}

Let $\ROMstate \equiv (\ROMvelocity; \ROMenergy;
\ROMposition)^T\in\RR{\sizeWholeROM}$, $\sizeWholeROM = \sizeROMvelocity +
\sizeROMenergy + \sizeROMposition$, be the reduced order hydrodynamic state
vector. Then the semidiscrete hyper-reduced system \eqref{eq:rom-hr-sns} can
be written in a compact form as
\begin{equation}\label{eq:combinedROM}
  \frac{d\ROMstate}{d\timeSymbol} = \ROMforceSys(\ROMstate,\timeSymbol;\param),  
\end{equation}
where the nonlinear force term, $\ROMforceSys:\RR{\sizeWholeROM} \times
\paramDomain \mapto \RR{\sizeWholeROM} $, is defined as 
\begin{align}\label{eq:combinedROMForce}
  \ROMforceSys(\ROMstate;\param) &\equiv
  \pmat{\ROMforceSysVelocity(\ROMvelocity, \ROMenergy, \ROMposition) \\
        \ROMforceSysEnergy(\ROMvelocity, \ROMenergy, \ROMposition) \\
        \ROMforceSysPosition(\ROMvelocity, \ROMenergy, \ROMposition)} 
        \equiv
  \pmat{-(\forceOneSamplingMat^T\forceOneBasis)^{\dagger}\forceOneSamplingMat^T \forceOne + \velocityBasis^T \mathbf{g}\\
        (\forceTvSamplingMat^T\forceTvBasis)^{\dagger}\forceTvSamplingMat^T \forceTv \\
        \positionBasis^T \velocityApprox}.
\end{align}
Applying the RK2-average scheme to the
hyper-reduced system \eqref{eq:rom-hr-sns}, the RK2-average fully discrete
hyper-reduced system reads:
\begin{align}\label{eq:RK2-avg-rom-hr}
  \ROMvelocityt{\timeIndex+\frac{1}{2}} &= \ROMvelocityt{\timeIndex} + (\timestepk{\timeIndex}/2)
    (-(\forceOneSamplingMat^T\forceOneBasis)^{\dagger}\forceOneSamplingMat^T \forceOneApproxk{\timeIndex} + \velocityBasis^T \mathbf{g}), &
    \ROMvelocityt{\timeIndex+1} &= \ROMvelocityt{\timeIndex} + \timestepk{\timeIndex}
    (-(\forceOneSamplingMat^T\forceOneBasis)^{\dagger}\forceOneSamplingMat^T \forceOneApproxk{\timeIndex+\frac{1}{2}} + \velocityBasis^T \mathbf{g}), \\
  \ROMenergyt{\timeIndex+\frac{1}{2}} &= \ROMenergyt{\timeIndex} + (\timestepk{\timeIndex}/2)
    (\forceTvSamplingMat^T\forceTvBasis)^{\dagger}\forceTvSamplingMat^T \forceTvApproxk{\timeIndex}, &
    \ROMenergyt{\timeIndex+1} &= \ROMenergyt{\timeIndex} + \timestepk{\timeIndex}
    (\forceTvSamplingMat^T\forceTvBasis)^{\dagger}\forceTvSamplingMat^T \avgforceTvApproxk{\timeIndex+\frac{1}{2}}, \\
  \ROMpositiont{\timeIndex+\frac{1}{2}} &= \ROMpositiont{\timeIndex} + (\timestepk{\timeIndex}/2)
    \positionBasis^T \velocityApproxt{\timeIndex+\frac{1}{2}}, & \ROMpositiont{\timeIndex+1} &=
    \ROMpositiont{\timeIndex} + \timestepk{\timeIndex} \positionBasis^T \avgvelocityApproxt{\timeIndex+\frac{1}{2}},
\end{align}
where the lifted ROM approximation $\stateApproxt{\timeIndex} 
= (\velocityApproxt{\timeIndex}; \energyApproxt{\timeIndex}; 
\positionApproxt{\timeIndex})^T \in\RR{\sizeWholeFE}$ given by 
\begin{align}\label{eq:discretesolrepresentation}
  \velocityApproxt{\timeIndex} &= \velocityOS + \velocityBasis\ROMvelocityt{\timeIndex},&
  \energyApproxt{\timeIndex} &= \energyOS + \energyBasis\ROMenergyt{\timeIndex},&
  \positionApproxt{\timeIndex} &= \positionOS + \positionBasis\ROMpositiont{\timeIndex}, 
\end{align}
is used to compute the updates 
\begin{align}
\forceOneApproxk{\timeIndex} & = \left
(\forceMat (\stateApproxt{\timeIndex}) \right ) \cdot \oneVec, &  
\forceTvApproxk{\timeIndex} & = \left
(\forceMat (\stateApproxt{\timeIndex}) \right )^T \cdot
\velocityApproxt{\timeIndex+\frac{1}{2}}, 
\end{align}
in the first stage. Similarly, $\stateApproxt{\timeIndex+\frac{1}{2}}
= (\velocityApproxt{\timeIndex+\frac{1}{2}}; \energyApproxt{\timeIndex+\frac{1}{2}}; 
\positionApproxt{\timeIndex+\frac{1}{2}})^T \in\RR{\sizeWholeFE}$ 
is used to computed the updates  
\begin{align}
\forceOneApproxk{\timeIndex+\frac{1}{2}} & = \left
(\forceMat (\stateApproxt{\timeIndex+\frac{1}{2}}) \right ) \cdot \oneVec, &
\avgforceTvApproxk{\timeIndex+\frac{1}{2}} & = \left
(\forceMat (\stateApproxt{\timeIndex + \frac{1}{2}}) \right )^T \cdot
\avgvelocityApproxt{\timeIndex+\frac{1}{2}}, 
\end{align}
with $\avgvelocityApproxt{\timeIndex+\frac{1}{2}} =
(\velocityApproxt{\timeIndex} + \velocityApproxt{\timeIndex+1})/2$ in the second
stage. The lifting is computed only for the sampled degrees of freedom, avoiding
full order computation. Again, the time step size $\timestepk{\timeIndex}$ is determined
adaptively using the automatic time step control algorithm,
with the state $\statet{\timeIndex}$ replaced
by the lifted ROM approximation $\stateApproxt{\timeIndex}$, 
and the time step size controlled by estimates at all quadrature points 
used in the evaluation of the force matrices in the hyper-reduction sample mesh.

\subsection{Solution bases construction}\label{sec:RBconstruction}
Now we describe how to obtain the reduced basis matrices 
$\velocityBasis$,  $\energyBasis$, and $\positionBasis$
for the solution variables. 
It suffices to describe how to construct
the reduced basis for the energy field only, i.e., $\energyBasis$, because other
bases will be constructed in the same way.
Proper orthogonal decomposition (POD) is commonly used to
construct a reduced basis.  POD
\cite{berkooz1993proper} obtains $\energyBasis$ from a truncated singular
value decomposition (SVD) approximation to a FOM solution snapshot matrix. It is
related to principal component analysis in statistical analysis
\cite{hotelling1933analysis} and Karhunen--Lo\`{e}ve expansion \cite{loeve1955}
in stochastic analysis.  In order to collect solution data for performing POD,
we run FOM simulations on a set of problem parameters, namely
$\{\param_{\paramIndex}\}_{\paramIndex=1}^{\nparam}$. 
For $\paramIndex\in\nat{\nparam}$, let $\finalTime(\param_\paramIndex)$ and $\ntimestep(\param_{\paramIndex})$ 
be the final time and the number of time steps in the training FOM simulation 
with the problem parameter $\param_\paramIndex$. 
By choosing $\energyOS(\param_\paramIndex) =
\energy(0;\param_\paramIndex)$, a solution snapshot matrix is formed by
assembling all the FOM solution data including the intermediate
Runge-Kutta stages, i.e.
 \begin{equation}
 \label{eq:snapshot_serial}
 \snapshots\equiv\bmat{\energyt{1}(\param_1)-\energyOS(\param_1) & \cdots &
 \energyt{\ntimestep(\param_{\nparam})}(\param_{\nparam})-\energyOS(\param_{\nparam})} \in
 \RR{\sizeThermodynamicFE\times\sizeSnapshot},
 \end{equation}
where $\energyt{\timeIndex}(\param_{\paramIndex})$ is the energy state at
$\timeIndex$th time step with problem parameter $\param_\paramIndex$ for
$\timeIndex\in\nat{\ntimestep(\param_{\paramIndex})}$ computed from the FOM
simulation, e.g. the fully discrete RK2-average scheme \eqref{eq:RK2-avg}, and
$\sizeSnapshot = \rkStage \sum_{\paramIndex=1}^{\nparam}
\ntimestep(\param_{\paramIndex})$ with $\rkStage$ being the number of 
Runge-Kutta stages in a time step. 
Then, POD computes its thin SVD: 
 \begin{align}\label{eq:SVD} 
   \snapshots &= \leftSingularMat\singularValueMat\rightSingularMat^T,
 \end{align} 
where $\leftSingularMat\in\RR{\sizeThermodynamicFE\times\sizeSnapshot}$ and
$\rightSingularMat\in\RR{\sizeSnapshot\times\sizeSnapshot}$ are orthogonal matrices,
and $\singularValueMat\in\RR{\sizeSnapshot\times\sizeSnapshot}$ is the diagonal
singular value matrix.  Then POD chooses the leading
$\sizeROMenergy$ columns of $\leftSingularMat$ to set
$\energyBasis = \bmat{\leftSingularVeck{1} & \ldots &
\leftSingularVeck{\sizeROMenergy}}$, where $\leftSingularVeck{\basisIndex}$ is
$\basisIndex$-th column vector of $\leftSingularMat$. The basis size,
$\sizeROMenergy$, is determined by the energy criteria, i.e., we find the
minimum $\sizeROMenergy\in\nat{\sizeSnapshot}$ such that the following
condition is satisfied:
 \begin{align}\label{eq:energy_criteria}
   \frac{\sum_{\basisIndex=1}^{\sizeROMenergy}
   \singularValue_{\basisIndex}}{\sum_{\basisIndex=1}^{\sizeSnapshot}
   \singularValue_{\basisIndex}} &\geq 1 - \singularValueThreshold,
 \end{align}
where $\singularValue_{\basisIndex}$ is a $\basisIndex$-th largest singular
value in the singular matrix, $\singularValueMat$, and $\singularValueThreshold
\in [0,1]$ denotes a threshold.
\footnote{We use the default value $\singularValueThreshold =
10^{-4}$ unless stated otherwise. }
The POD basis minimizes $\|\snapshots - \energyBasis\energyBasis^T\snapshots
\|_F^2$ over all $\energyBasis\in\RR{\sizeThermodynamicFE \times
\sizeROMenergy}$ with orthonormal columns, where $\|\dummyMat\|_F$ denotes the
Frobenius norm of a matrix $\dummyMat\in\RR{I\times J}$, defined as
$\|\dummyMat\|_F = \sqrt{\sum_{i=1}^{I}\sum_{j=1}^{J} \dummyElemSymbol_{ij}^2}$
with $\dummyElemSymbol_{ij}$ being the $(i,j)$ element of $\dummyMat$.  Since
the objective function does not change if $\energyBasis$ is post-multiplied by
an arbitrary $\sizeROMenergy\times\sizeROMenergy$ orthogonal matrix, the POD
procedure seeks the optimal $\sizeROMenergy$-dimensional subspace that
captures the snapshots in the least-squares sense.  For more details on POD, we
refer to \cite{hinze2005proper,kunisch2002galerkin}.
The same procedure can be used to construct the other solution bases $\velocityBasis$ and $\positionBasis$.
Finally, the reduced bases for the nonlinear terms are 
constructed by the subspace relations $\forceOneBasis =
\kinematicMassMat \velocityBasis$ and $\forceTvBasis = \thermodynamicMassMat \energyBasis$. 

\subsection{Sampling indices selection}\label{sec:deim}
It remains to describe how to obtain the sampling matrices, i.e.
$\forceOneSamplingMat$ and $\forceTvSamplingMat$.  The discrete empirical
interpolation method (DEIM) is a popular choice for nonlinear model reduction.
It suffices to describe how to construct the sampling matrix for the
momentum nonlinear term only, i.e., $\forceOneSamplingMat$, as the other
matrix will be constructed in the same way.  The sampling matrix
$\forceOneSamplingMat$ is characterized by the sampling indices
$\{p_1,\ldots,p_{\sizeROMforceOneSample}\}$, which can be found either by a row
pivoted LU decomposition \cite{chaturantabut2010nonlinear} or the strong column
pivoted rank-revealing QR (sRRQR) decomposition \cite{drmac2016new,
drmac2018discrete}.  Algorithm 1 of \cite{chaturantabut2010nonlinear} uses the
greedy algorithm to sequentially seek additional interpolating indices
corresponding to the entry with the largest magnitude of the residual of
projecting an active POD basis vector onto the preceding basis vectors at the
preceding interpolating indices.  The number of interpolating indices
returned is the same as the number of basis vectors, i.e. $\sizeROMforceOneSample =
\sizeROMforceOne$.  Algorithm 3 of \cite{carlberg2013gnat} and Algorithm 5 of
\cite{carlberg2011efficient} use the greedy procedure to minimize the error in
the gappy reconstruction of the POD basis vectors $\forceOneBasis$.  These
algorithms allow over-sampling, i.e. $\sizeROMforceOneSample \geq
\sizeROMforceOne$, and can be regarded as extensions of Algorithm 1 of
\cite{chaturantabut2010nonlinear}.  Instead of using the greedy algorithm,
Q-DEIM is introduced in \cite{drmac2016new} as a new framework for constructing
the DEIM projection operator via the QR factorization with column pivoting.
Depending on the algorithm for selecting the sampling indices, the DEIM
projection error bound is determined. For example, the row pivoted LU
decomposition in \cite{chaturantabut2010nonlinear} results in the following
error bound:
 \begin{align}\label{eq:errorbound_DEIM}
   \|\forceOne - \forceOneObliqProjMat \forceOne \|_2 \leq \conditionnumber
   \|(\identityMat{\sizeKinematicFE}-\forceOneBasis\forceOneBasis^T)\forceOne \|_2,
 \end{align}
where $\|\cdot\|_2$ denotes $\ell_2$ norm of a vector, and $\conditionnumber$ is
the condition number of $(\forceOneSamplingMat^T\forceTvBasis)^{-1}$, bounded by
 \begin{align}\label{eq:crudeBound}
   \conditionnumber \leq
   (1+\sqrt{2\sizeKinematicFE})^{\sizeROMforceOne-1}\|\forceOneBasisVeck{1}\|_\infty^{-1}.
 \end{align}
On the other hand, the sRRQR factorization in \cite{drmac2018discrete} reveals a
tighter bound than \eqref{eq:crudeBound}:
 \begin{align}\label{eq:tighterBound}
   \conditionnumber \leq
   \sqrt{1+\tuningparam^2\sizeROMforceOne(\sizeKinematicFE-\sizeROMforceOne)}
 \end{align}
where $\tuningparam$ is a tuning parameter in the sRRQR factorization (i.e.,
$f$ in Algorithm 4 of \cite{gu1996efficient}). 
If SNS is used, the nonlinear term basis $\forceOneBasis$ is not necessarily orthogonal, 
but we can still obtain similar error bounds. For example, when the sRRQR factorization is used,
\eqref{eq:errorbound_DEIM} becomes
 \begin{align}\label{eq:errorbound_SNS_DEIM}
   \|\forceOne - \forceOneObliqProjMat \forceOne \|_2 \leq \conditionnumber
   \|(\identityMat{\sizeKinematicFE}-\forceOneBasis\forceOneBasis^\dagger)\forceOne \|_2,
 \end{align}
 with $\conditionnumber$ satisfying \eqref{eq:tighterBound} (c.f. Theorem 5.2 of \cite{choi2020sns}). 
Once the sampling indices are determined, the reduced matrix 
$\left(\forceOneSamplingMat^T\forceOneBasis\right)^{\dagger}
  \in \RR{\sizeROMforceOne\times\sizeROMforceOneSample}$ 
 can be precomputed and stored. 

\section{Decomposition of solution manifold}\label{sec:decompose}

Although the reduced order model reported in Section~\ref{sec:ROM} 
can serve as a powerful technique in short-time simulation of Rayleigh--Taylor instability, 
it is not applicable to long-time simulation. 
The nonlinear evolution leads to the formation of mushroom-shaped vortices, 
which is the main interest of numerical studies in the hydrodynamic instability. 
The main difficulty is that, in order to capture the formation of the vortex, 
the time step is adaptively reduced. In some numerical examples presented, 
the full order model simulation ends up with $10^5$ to $10^7$ time steps. 
This gives rise to a huge number of snapshot samples, which imposes a heavy 
burden in storage and computational cost for the SVD computations. 
On the other hand, the penetration of the fluid interface grows with time and 
characterizes the solution dynamics, i.e., advection-dominated. As a result, the linear dependence among the 
snapshots is weak, and therefore there is no intrinsic low-dimensional subspace 
which can universally approximate the solution manifold, comprised of 
all the solutions over the parameter domain and temporal domain. 
Therefore, it is impossible for the reduced order model approach 
to achieve any meaningful speed-up and good accuracy. 

To this end, we propose a framework to overcome these difficulties by 
employing multiple reduced order models. 
In the offline phase, we construct each of these reduced order models 
from a small subset of the snapshot samples, to ensure low dimension. 
In the online phase, each of these reduced order models will be used 
in a subdomain of time and parameter where it is supposed to provide good approximation. 
This framework involves a decomposition of the solution manifold 
and relies on an indicator which will be used to classify the snapshot samples 
and assign the reduced-order models. 
The idea is to decompose the solution manifold into submanifolds 
where the Kolmogorov $n$-width decaying fast with respect to the subspace dimension, within which we can 
collect snapshots with strong linear dependence. This enables us to build accurate multiple low-dimensional subspaces.

\subsection{General framework}
We describe the general framework of the decomposition of the solution manifold, 
from which we will derive two practical examples later in this section. 
Let $\windowIndicator: \RR{\sizeWholeFE} \times \mathbb{R}^+ \times 
\paramDomain \to \mathbb{R}$ be an indicator 
which maps the triplet $(\state, \timeSymbol, \param)$ to a real value. 
The range of the indicator $[\windowIndicator_\text{min}, \windowIndicator_\text{max})$ 
will be partitioned into $\nwindow^\text{off}$ subintervals, i.e. 
\begin{equation}
\windowIndicator_\text{min} = \windowIndicator_{0} < \windowIndicator_{1} < \cdots <
\windowIndicator_{\nwindow^\text{off}-1} < \windowIndicator_{\nwindow^\text{off}} = \windowIndicator_\text{max}. 
\end{equation}
This partition can be either prescribed or 
determined by the snapshot data collected in the offline phase. 
For simplicity, we assume that 
$\windowIndicator(\state(0), 0, \param) = \Psi_\text{0}$, and 
$\windowIndicator(\state(\timeSymbol), \timeSymbol, \param)$ 
is increasing with time $\timeSymbol$ for any initial state $\state(0)$ and any 
parameter $\param \in \paramDomain$. 
These assumptions are valid in the concrete examples we present in this paper. 
In general, one can relax these assumptions.

With the partition of the indicator range, instead of directly assembling all the snapshot samples into a single huge 
snapshot matrix as in \eqref{eq:snapshot_serial}, 
the FOM states will be first classified into groups. 
Let $\windowIndex \in \nat{\nwindow^\text{off}}$ 
be a group index. 
We denote the subset of paired indices of time and offline parameter 
whose corresponding snapshot belongs to the $\windowIndex$th group as 
\begin{equation}\label{eq:snapshot_cluster}
\snapshotGroup_{\windowIndex} 
= \left\{ (\timeIndex, \paramIndex) \in \mathbb{Z} \times \nat{\nparam}: 
0 \leq \timeIndex \leq \ntimestep(\param_{\paramIndex}) \text{ and }
\windowIndicator\left(\statet{\timeIndex}(\param_\paramIndex), \timeSymbol_{\timeIndex}(\param_\paramIndex), 
\param_\paramIndex\right) \in [\windowIndicator_{\windowIndex-1}, \windowIndicator_{\windowIndex})
\right\}. 
\end{equation} 
Then the snapshot matrix $\snapshots_\windowIndex$ 
of the energy variable in the $\windowIndex$th group will be formed by assembling the corresponding snapshots, i.e. 
 \begin{equation}
 \label{eq:snapshot_group}
 \snapshots_\windowIndex \equiv\bmat{\energyt{\timeIndex}(\param_\paramIndex)-\energyOSWindow{\windowIndex}(\param_\paramIndex)
 }_{(\timeIndex, \paramIndex) \in \snapshotGroup_{\windowIndex}}, 
 \end{equation}
and POD is used to construct the energy solution basis $\energyBasisWindow{\windowIndex}$ as described in Section~\ref{sec:RBconstruction}.
The same procedure can be used to construct the other solution bases $\velocityBasisWindow{\windowIndex}$ and 
$\positionBasisWindow{\windowIndex}$.

For a generic problem parameter $\param \in \paramDomain$, 
let $\finalTimeROM(\param)$ be the final time in the ROM simulation 
with the problem parameter $\param$, 
where the temporal domain is discretized as $\{
  \timekROM{\timeIndex} \}_{\timeIndex=0}^{\ntimestepROM(\param)}$.
We remark that $\finalTimeROM(\param)$ can be different from the 
final time $\{\finalTime(\param_\paramIndex)\}_{\paramIndex \in \nat{\nparam}}$ used in the snapshot sampling, 
and that even with the same problem setting and same final time, it is very likely that the temporal
discretization used in the hyper-reduced system is different from the full
order model, since the temporal discretization is adaptively controlled by the states.  
The computation in the online phase will be performed 
using different reduced bases 
in $\nwindow(\param)$ subintervals 
of the temporal domain $[0,\finalTimeROM(\param)]$, i.e. 
\begin{equation}
0 = \windowk{0}(\param) < \windowk{1}(\param) < \cdots 
\windowk{\nwindow(\param)-1}(\param) < \windowk{\nwindow(\param)}(\param) = \finalTimeROM(\param), 
\end{equation}
with $\{  \windowk{\windowIndex}(\param)  \}_{\windowIndex=1}^{\nwindow(\param)-1}$ being the 
partition points of the temporal domain.
We remark that $\nwindow(\param)$ is the number of subintervals that 
the temporal domain $[0,\finalTimeROM(\param)]$ is decomposed into for the parameter $\param$, 
and does not exceed the number of groups $\nwindow^\text{off}$. 
In general, the temporal domain partition is parameter-dependent, i.e. 
the number of subintervals $\nwindow(\param)$ and the end-point of each subinterval  
$\{  \windowk{\windowIndex}(\param)  \}_{\windowIndex=1}^{\nwindow(\param)}$ 
depend on the problem parameter $\param$, 
and they are assigned by the indicator $\windowIndicator$ 
iteratively in the time marching of the ROM simulation. 
More precisely, starting with 
$\windowk{0}(\param) = 0$, for $\windowIndex \geq 1$ and 
$\timeSymbol \in \timeWindow_{\windowIndex}(\param) 
\equiv [\windowk{\windowIndex-1}(\param), \windowk{\windowIndex}(\param))$, 
we employ the reduced order model in the form of \eqref{eq:rom-hr-sns} 
by projecting onto the reduced subspaces spanned by the corresponding bases 
$\velocityBasisWindow{\windowIndex}$, $\energyBasisWindow{\windowIndex}$, and 
$\positionBasisWindow{\windowIndex}$, i.e. 
\begin{align}\label{eq:rom-hr-tw}
   \frac{d\ROMvelocityWindow{\windowIndex}}{d\timeSymbol} &=
    -((\forceOneSamplingMatWindow{\windowIndex} )^T\forceOneBasisWindow{\windowIndex} )^{\dagger}(\forceOneSamplingMatWindow{\windowIndex})^T \forceOne(\velocityOSWindow{\windowIndex} + \velocityBasisWindow{\windowIndex}\ROMvelocityWindow{\windowIndex}, \energyOSWindow{\windowIndex} + \energyBasisWindow{\windowIndex}\ROMenergyWindow{\windowIndex}, \positionOSWindow{\windowIndex} + \positionBasisWindow{\windowIndex}\ROMpositionWindow{\windowIndex}, \timeSymbol; \param) 
    + (\velocityBasisWindow{\windowIndex})^T \mathbf{g} \\
   \frac{d\ROMenergyWindow{\windowIndex}}{d\timeSymbol} &=
    \left((\forceTvSamplingMatWindow{\windowIndex})^T
      \forceTvBasisWindow{\windowIndex} \right)^{\dagger} (\forceTvSamplingMatWindow{\windowIndex})^T \forceTv(\velocityOSWindow{\windowIndex} + \velocityBasisWindow{\windowIndex}\ROMvelocityWindow{\windowIndex}, \energyOSWindow{\windowIndex} + \energyBasisWindow{\windowIndex}\ROMenergyWindow{\windowIndex}, \positionOSWindow{\windowIndex} + \positionBasisWindow{\windowIndex}\ROMpositionWindow{\windowIndex}, \timeSymbol; \param) \\
    \frac{d\ROMpositionWindow{\windowIndex}}{d\timeSymbol} &= (\positionBasisWindow{\windowIndex})^T \velocityOSWindow{\windowIndex} +
    (\positionBasisWindow{\windowIndex})^T \velocityBasisWindow{\windowIndex} \ROMvelocityWindow{\windowIndex}, 
\end{align}
where the offset vectors $\velocityOSWindow{\windowIndex}(\param)$, 
$\energyOSWindow{\windowIndex}(\param)$, and $\positionOSWindow{\windowIndex}(\param)$ are prescribed. 
The initial condition is given by
lifting the ROM solution to the FOM spaces, using the ROM bases in the time
window $\windowIndex-1$, and then projecting onto the ROM spaces in the time window $\windowIndex$, i.e.
\begin{align}\label{eq:ic-tw}
\ROMvelocityWindow{\windowIndex}(\windowk{\windowIndex-1}(\param), \param) & = 
(\velocityBasisWindow{\windowIndex})^T ( 
\velocityOSWindow{\windowIndex-1}(\param) + 
  \velocityBasisWindow{\windowIndex-1}\ROMvelocityWindow{\windowIndex-1}(\windowk{\windowIndex-1}(\param), \param) - \velocityOSWindow{\windowIndex}(\param) ), \\
\ROMenergyWindow{\windowIndex}(\windowk{\windowIndex-1}(\param), \param) & = 
(\energyBasisWindow{\windowIndex})^T ( 
\energyOSWindow{\windowIndex-1} (\param) + 
  \energyBasisWindow{\windowIndex-1}\ROMenergyWindow{\windowIndex-1}(\windowk{\windowIndex-1}(\param), \param) - \energyOSWindow{\windowIndex}(\param) ), \\
\ROMpositionWindow{\windowIndex}(\windowk{\windowIndex-1}(\param), \param) & = 
(\positionBasisWindow{\windowIndex})^T ( 
\positionOSWindow{\windowIndex-1} (\param) + 
  \positionBasisWindow{\windowIndex-1}\ROMpositionWindow{\windowIndex-1}(\windowk{\windowIndex-1}(\param), \param) - \positionOSWindow{\windowIndex}(\param) ).
\end{align}
It should be remarked that the end-point $\windowk{\windowIndex}(\param)$ 
is determined adaptively during the ROM simulation. 
As in \eqref{eq:RK2-avg-rom-hr}, 
we can derive the RK2-average fully discrete
hyper-reduced system by applying the RK2-average scheme to the
hyper-reduced system \eqref{eq:rom-hr-tw}. 
We denote by $\ntimestepROM(\param)$ the number of time steps in the ROM simulation 
with the problem parameter $\param$, 
where the temporal domain is discretized as $\{
  \timekROM{\timeIndex}(\param) \}_{\timeIndex=0}^{\ntimestepROM(\param)}$.
We remark that even with the same problem setting and same final time, it is very likely that the temporal
discretization used in the hyper-reduced system is different from the full
order model, due to the adaptive time step size.  
The end-point $\windowk{\windowIndex}(\param)$ of the temporal subinterval $\timeWindow_{\windowIndex}(\param)$ 
is defined as the time instance $\timekROM{\timeIndex}(\param)$
when $\windowIndicator(\stateApproxt{\timeIndex}(\param), \timekROM{\timeIndex}(\param), \param)$ 
first exceeds $\windowIndicator_{\windowIndex}$, 
at which we increment to the next time subinterval $\timeWindow_{\windowIndex+1}(\param)$.

\begin{figure}[ht!]
\centering
\includegraphics[height=1.4in]{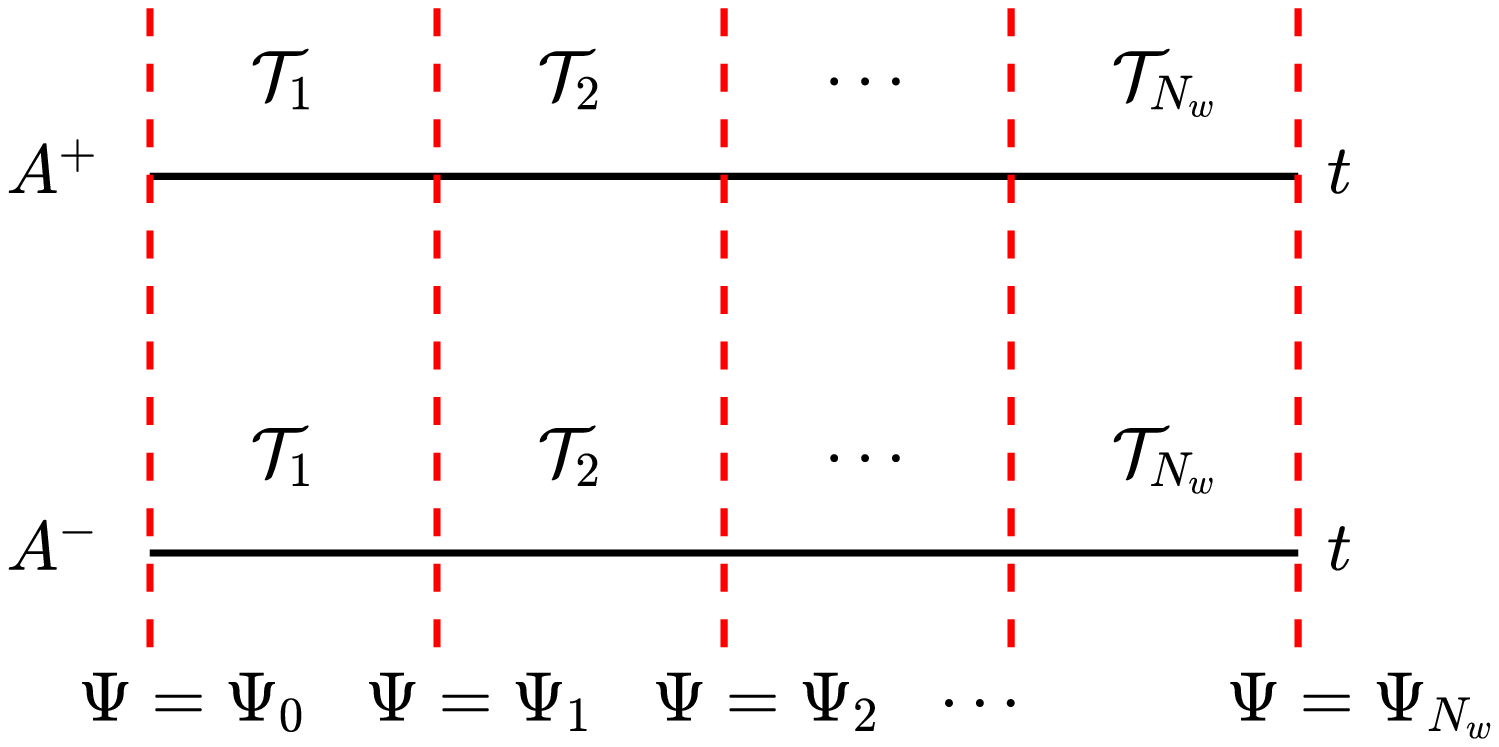}
\hspace{0.05\linewidth}
\includegraphics[height=1.4in]{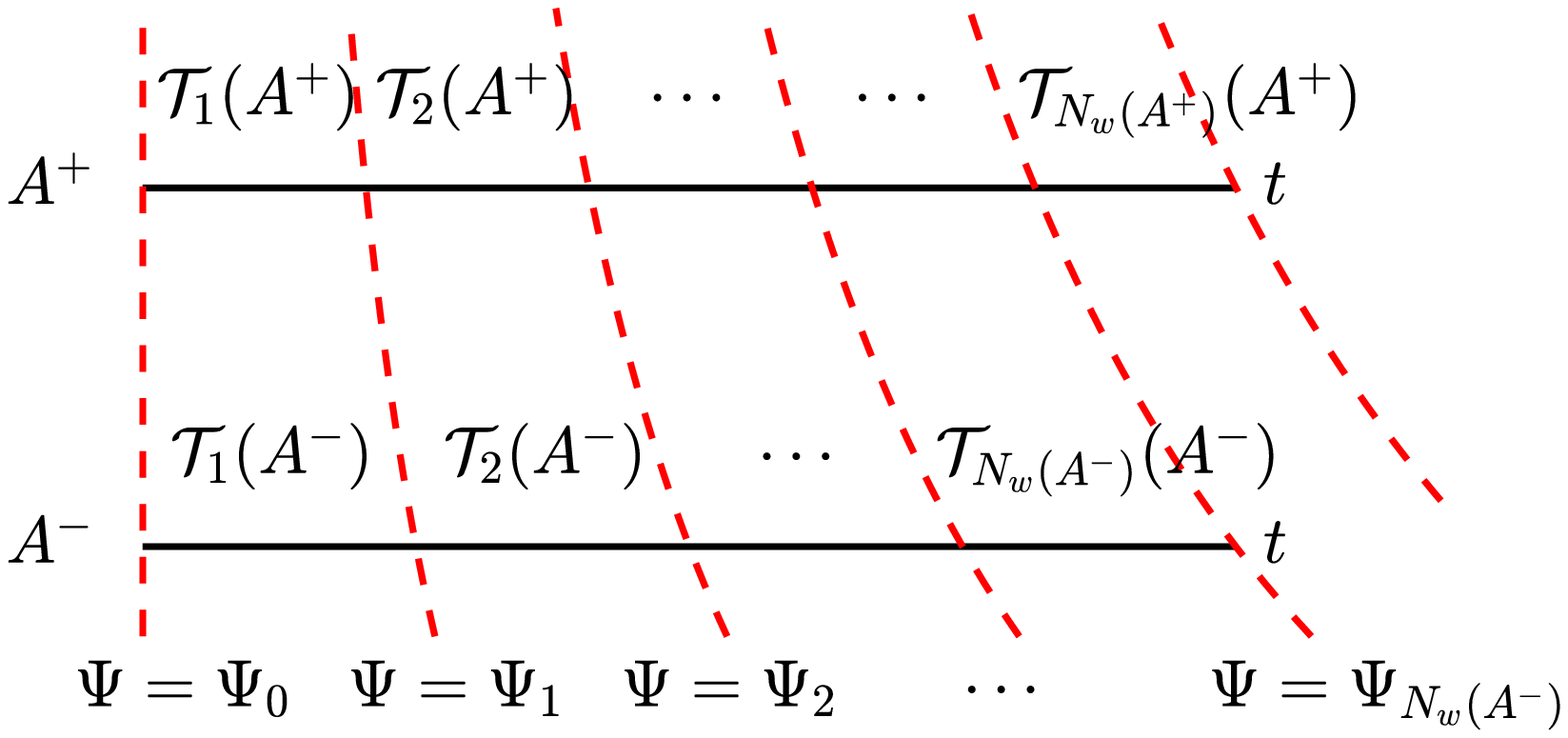}
\caption{Illustration of the parameter-dependence of temporal domain partition
  using physical time (left) and penetration distance (right) for manifold decomposition 
  for two generic testing parameters $\param^-, \param^+ \in \paramDomain$, where $\param^- < \param^+$. 
  Note that we have $\nwindow(\param^-) \leq \nwindow(\param^+)$. 
  }
\label{fig:idea}
\end{figure}

Figure~\ref{fig:idea} illustrates a possible scenario of the temporal domain partition 
for two generic testing parameters $\param^-, \param^+ \in \paramDomain$, where $\param^- < \param^+$,
with two decomposition mechanisms which we are going to 
introduce in the rest of this section. It should be highlighted that  
the partition of the temporal domain is parameter-dependent if penetration distance is used 
as the indicator $\windowIndicator$, while it is parameter-independent if physical time is used. 
This will be further explained in later subsections. 

\subsection{Time-windowing: decomposition by physical time}\label{sec:tw}
A very natural choice of the indicator is the physical time, i.e. 
$\windowIndicator(\state, \timeSymbol, \param) = \timeSymbol$. 
In that case, for any problem parameter $\param \in \paramDomain$, 
we have $\nwindow(\param) = \nwindow^\text{off}$ and 
$\windowk{\windowIndex}(\param) = \windowIndicator_{\windowIndex}$ 
for all $1 \leq \windowIndex \leq \nwindow^\text{off}$. 
In other words, the temporal partition is parameter-independent. 
As in the figure on the left in Figure~\ref{fig:idea}, we can safely remove the parameter dependence
and the superscript for distinguishing offline and online indices in the temporal domain partition.  
We remark that this choice reduces to the time-windowing (TW) reduced order model 
approach in Section 4 of \cite{copeland2022reduced}. 

One may decompose the temporal domain into windows using a prescribed partition,
that is, the end points $\{ \windowIndicator_{\windowIndex}
\}_{\windowIndex=1}^{\nwindow}$ are user-defined. 
Heuristically, the number of snapshots for a variable in a window is proportional to
the window size if a uniform time step is used in the FOM simulation.  In
addition, if a uniform partition is used to decompose the domain into windows,
then the basis sizes will be heuristically balanced among the windows.  However,
in the setting of adaptive time stepping, it becomes unclear how to
prescribe window sizes to balance the ROM basis size among time windows. 
As an alternative, one can adaptively divide the temporal domain into windows using a
prescribed number of snapshots per window, which can be heuristically related to the basis size. 
This adaptive decomposition approach is described in Algorithm~\ref{alg:decompose} in the context of a general indicator $\windowIndicator$ and 
is used in our numerical experiments in Section~\ref{sec:numerical}.

\begin{algorithm} \caption{Adaptive partition of indicator range} \label{alg:decompose}
\begin{algorithmic}[1]
\State \textbf{Input:} 
\begin{itemize}
\item $\left\{ (\statet{\timeIndex}(\param_\paramIndex), \timeSymbol_{\timeIndex}(\param_\paramIndex), 
\param_\paramIndex) : 0 \leq \timeIndex \leq \ntimestep(\param_{\paramIndex}) 
\text{ and } \paramIndex \in \nat{\nparam} \right\}$, the collection of all snapshot triplets, 
\item $\ntimestepWindow$, maximum number of intermediate samples in a subinterval for each offline parameter, 
\item $\windowIndicator: \RR{\sizeWholeFE} \times \mathbb{R}^+ \times 
\paramDomain \to \mathbb{R}$, the partition indicator, 
\item $[\windowIndicator_\text{min}, \windowIndicator_\text{max})$, the range of the collection of all snapshot triplets under the indicator
\end{itemize}
\State Initilize $\windowIndex = 0$, $\windowIndicator_0 = \windowIndicator_{\text{min}}$, $\mathcal{K} = \nat{\nparam}$, and 
$\timeIndexWindow{0}(\param_{\paramIndex}) = 0$ for all $\paramIndex \in \nat{\nparam}$ 
\Repeat
\State Increment $\windowIndex \leftarrow \windowIndex+1$
\For {$\paramIndex \in \mathcal{K}$}
\State Set $\timeIndexWindow{\windowIndex, {\paramIndex}}' = \min\{{\timeIndexWindow{\windowIndex-1}(\param_{\paramIndex})+\ntimestepWindow+1}, \ntimestep(\param_{\paramIndex})\}$
\State Set $\windowIndicator_{\windowIndex, \paramIndex}  = \windowIndicator(\statet{\timeIndexWindow{\windowIndex, {\paramIndex}}'}(\param_{\paramIndex}), 
\timek{\timeIndexWindow{\windowIndex, {\paramIndex}}' }, \param_{\paramIndex})$
\EndFor
\State Set the endpoint of the $\windowIndex$-th subinterval as 
$\windowIndicator_{\windowIndex} = \min_{\paramIndex \in \mathcal{K} } \windowIndicator_{\windowIndex, \paramIndex}$
\State Initialize the $\windowIndex$-th snapshot index group by $\snapshotGroup_\windowIndex = \emptyset$
\For {$\paramIndex \in \mathcal{K}$}
\State Set
$\timeIndexWindow{\windowIndex}(\param_{\paramIndex})$ as the latest time index in the subinterval $\timeWindow_{\windowIndex}$ for the offline parameter $\param_{\paramIndex}$ , i.e.
\begin{equation}
\windowIndicator(\statet{\timeIndexWindow{\windowIndex}(\param_{\paramIndex})}(\param_{\paramIndex}), \timek{\timeIndexWindow{\windowIndex}(\param_{\paramIndex})}, \param_{\paramIndex}) 
\leq \windowIndicator_{\windowIndex} < \windowIndicator(\statet{\timeIndexWindow{\windowIndex}(\param_{\paramIndex})+1}(\param_{\paramIndex}), \timek{\timeIndexWindow{\windowIndex}(\param_{\paramIndex})+1}, \param_{\paramIndex})
\end{equation}
\State Enrich the $\windowIndex$-th snapshot index group by 
\begin{equation}
\snapshotGroup_\windowIndex \leftarrow \snapshotGroup_\windowIndex \cup 
\{ (\timeIndex, \paramIndex) \in \mathbb{Z} \times \nat{\nparam}: 
\timeIndexWindow{\windowIndex-1}(\param_{\paramIndex}) \leq \timeIndex \leq \timeIndexWindow{\windowIndex}(\param_{\paramIndex}) \}
\end{equation}
\If{$\timeIndexWindow{\windowIndex}(\param_{\paramIndex}) = \ntimestep(\param_\paramIndex)$}
\State Remove $\paramIndex$ from the index set of active training parameters by $\mathcal{K} \leftarrow \mathcal{K} \setminus \{ \paramIndex \}$
\EndIf
\EndFor
\Until {$\mathcal{K} = \emptyset$}
\State Set $\nwindow^\text{off} = j$ and verify $\windowIndicator_{\nwindow^\text{off}} < \windowIndicator_\text{max}$
\State \textbf{Output:} 
\begin{itemize}
\item $\nwindow^\text{off}$, number of subintervals, 
\item $\{ \windowIndicator_{\windowIndex} \}_{\windowIndex = 0}^{\nwindow^\text{off}}$, partition of indicator range, and 
\item $\{ \snapshotGroup_{\windowIndex} \}_{\windowIndex = 1}^{\nwindow^\text{off}}$, sequence of snapshot index groups
\end{itemize}
\end{algorithmic}
\end{algorithm}

\begin{figure}[ht!]
\centering
\includegraphics[width=0.98\linewidth]{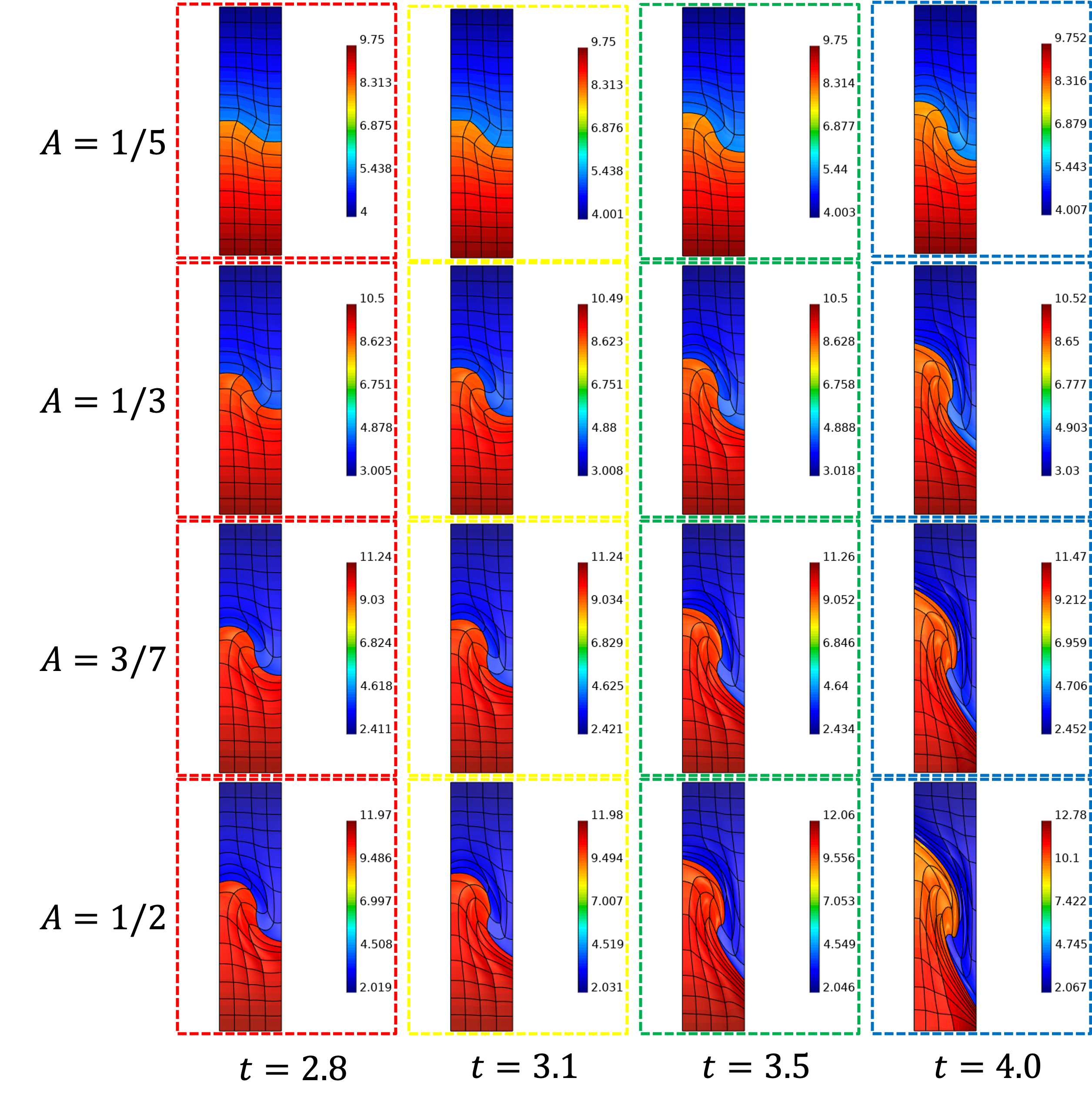}
\caption{An illustrative example to explain the mechanism of 
decomposition of solution manifold using physical time as indicator. 
The samples collected at the same time instance belong to the same group. 
The samples in each group is surrounded by a box with the same color, namely 
red, yellow, green, and blue. 
Each group has weak linear dependence since the samples 
penetrate to different extents.}
\label{fig:tw_example}
\end{figure}

Figure~\ref{fig:tw_example} illustrates a possible scenario of 
the classification of snapshot samples with 
physical time being the indicator. 
The sixteen snapshot samples of the energy variable are collected from 
four different time instances and four different density ratio (and therefore problem parameter). 
With a sufficiently fine partition of the indicator range, i.e. the temporal domain, 
each one of the sixteen snapshot samples will belong to one of the four different groups, 
where the samples in each group are surrounded by a box with the same color, namely 
red, yellow, green, and blue. 
It can be seen that the fluid interface penetrates to different extents  
for the snapshot samples in the same group. This suggests that 
the linear dependence among the snapshots is weak, and 
it is difficult to achieve a linear subspace reduced basis with a small dimension, limiting speed-up.

\subsection{Distance windowing: decomposition by penetration distance}\label{sec:dw}
To reduce the size of the solution basis, we propose to use the downward
penetration distance, which characterizes the Rayleigh-Taylor instability
problem well. The penetration distance is defined as the largest downward
displacement of the interface, i.e.,  
\begin{equation}
\windowIndicator(\state, \timeSymbol, \param) = \max_{\initialPosition_2 = 0} 
\positionSymbol_2(\initialPosition), 
\end{equation}
where $\positionSymbol_2\in\RR{}$ is the vertical component of the functional representation $\positionSymbol:\initialDomain\rightarrow\RR{2}$ 
of the position vector $\position$ of a point on the interface, i.e., $\initialPosition_2=0$,
where $\initialPosition$ is the Eulerian coordinates, or equivalently, the initial Lagrangian coordinates. 
We remark that, by taking the advantage of Lagrangian formulation of the Euler equations,
the penetration distance is easily accessible. 
Since the penetration speed depends on the density ratio and hence the Atwood number, 
the temporal domain partition is parameter-dependent. 
The figure on the right in Figure~\ref{fig:idea} shows the temporal domain partition 
of two generic testing parameters $\param^-, \param^+ \in \paramDomain$, where $\param^- < \param^+$. 
With a higher Atwood number $\param^+$, the fluid interface penetrates faster, and 
it takes shorter time for the penetration distance to reach the endpoints $\windowIndicator_\windowIndex$. 
This generates more but shorter subintervals in the temporal partition of $\param^+$ than those of $\param^-$, i.e. 
$\nwindow(\param^-) \leq \nwindow(\param^+)$.  
We name this approach distance-windowing (DW) reduced order model. 

Again, one may decompose the indicator range of penetration distance 
into subintervals using a prescribed partition of the end points $\{ \windowk{\windowIndex}
\}_{\windowIndex=1}^{\nwindow}$.
Besides the difficulties discussed in TW-ROM, 
the indicator range of penetration distance is not known prior to obtaining all the snapshots, 
and it is even more difficult to determine a good partition. 
Similar to the case of physical time, one can adaptively divide the indicator range of penetration distance into subintervals using a
prescribed number of snapshots per window, which can be heuristically related to the basis size. 
Again, we refer the readers to Algorithm~\ref{alg:decompose} for the details. 
The adaptive decomposition approach is used in our numerical experiments in Section~\ref{sec:numerical}.

\begin{figure}[ht!]
\centering
\includegraphics[width=0.98\linewidth]{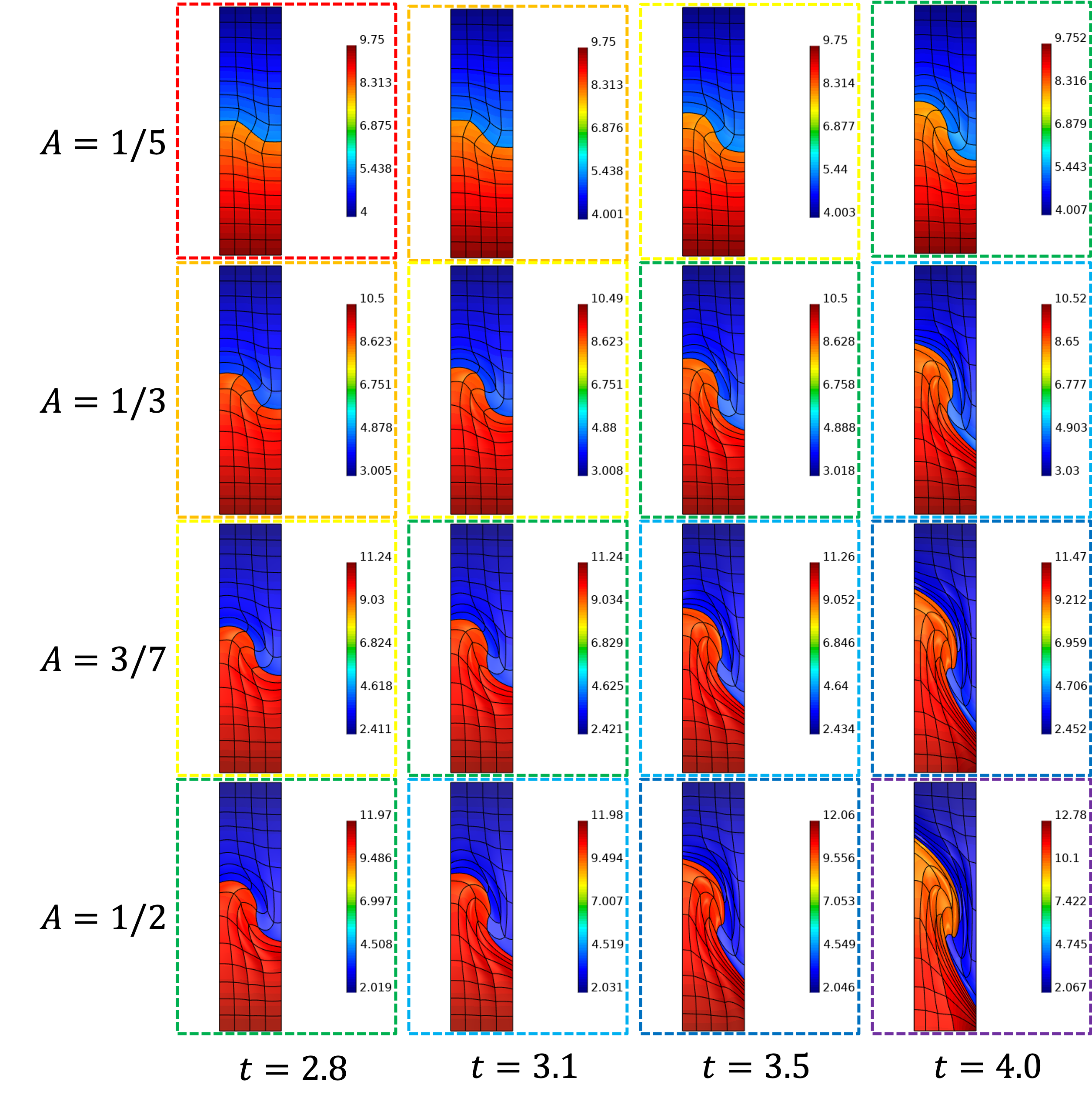}
\caption{An illustrative example to explain the mechanism of 
decomposition of solution manifold using penetration distance as indicator.
The samples penetrating to similar extent belong to the same group.
The samples in each group is surrounded by a box with the same color, namely 
red, orange, yellow, green, cyan, blue, and purple. 
Each group has strong linear dependence.}
\label{fig:dw_example}
\end{figure}

Figure~\ref{fig:dw_example} illustrates a possible scenario of 
the classification of snapshot samples with 
penetration distance being the indicator. 
Each one of the sixteen snapshot samples of the energy variable will belong to one of the seven different groups, 
where the samples in the same group are surrounded by a box with the same color, namely 
red, orange, yellow, green, cyan, blue, and purple. 
A sample collected at the earlier time instance with higher density ratio is 
associated with one collected at the later time instance with lower density ratio, 
since the penetration speed is slower with low density ratio.  
As we can observe, this choice of indicator provides much 
better linear dependence among the snapshots within the same group, as 
the fluid interface penetrates to a similar extent. 

\section{Numerical studies}\label{sec:numerical}
In this section, we present numerical results to test the performance of 
our proposed method applied to Rayleigh--Taylor instability in Lagrangian hydrodynamics
simulated with LaghosROM\footnote{GitHub page, {\it
https://github.com/CEED/Laghos/tree/rom/rom}.} and libROM\footnote{GitHub page, {\it
https://github.com/LLNL/libROM}; libROM webpage, {\it https://www.librom.net}}. 
LaghosROM is a C++ miniapp that
accelerates time-dependent Euler equations of compressible gas dynamics in a
moving Lagrangian frame using unstructured high-order finite element spatial
discretization, explicit high-order time-stepping, and projection-based reduced
order models. libROM is a C++ library for data-driven physical simulation methods
\cite{choi2019librom}.
All the simulations in this section use the machine Quartz in the Livermore Computing Center\footnote{High
performance computing at LLNL, {\it https://hpc.llnl.gov/hardware/platforms/Quartz}},
on Intel Xeon CPUs with 128 GB memory, peak TFLOPS of 3251.4, and peak single
CPU memory bandwidth of 77 GB/s. 

We compare the ROM performance with different FOM discretization levels, 
ROM parameters, and decomposition mechanisms at different Atwood numbers. 
Algorithm~\ref{alg:decompose} reports the adaptive decomposition mechanism in our implementation, 
where the indicator $\windowIndicator$  is chosen as the physical time as in Section~\ref{sec:tw}
or the penetration distance as in Section~\ref{sec:dw}, with $\ntimestepWindow$ being a user-defined constant 
to control the maximum number of intermediate samples to be added to a group for each training parameter. 
Moreover, the over-sampling factors $\{\factorROMforceOneSample, \factorROMforceTvSample\}$ are user-defined
constants over all time windows, which control the number of sampling indices 
$\sizeROMforceOneSampleWindow{\windowIndex}$ and $\sizeROMforceTvSampleWindow{\windowIndex}$ 
in the sampling matrices  $\forceOneSamplingMatWindow{\windowIndex}$ and 
$\forceTvSamplingMatWindow{\windowIndex}$ of the window $\windowIndex$ by 
\begin{equation}
\sizeROMforceOneSampleWindow{\windowIndex} = \min\left\{ \sizeKinematicFE,  \, 
\factorROMforceOneSample\sizeROMforceOneWindow{\windowIndex} \right\}, 
\quad 
\sizeROMforceTvSampleWindow{\windowIndex} = \min\left\{ \sizeThermodynamicFE, \, 
\factorROMforceTvSample\sizeROMforceTvWindow{\windowIndex} \right\}, 
\end{equation}
where $\sizeROMforceOneWindow{\windowIndex}$ and $\sizeROMforceTvWindow{\windowIndex}$ 
are the numbers of columns of the nonlinear term bases $\forceOneBasisWindow{\windowIndex}$ and 
$\forceTvBasisWindow{\windowIndex}$, respectively.
To evaluate the ROM performance, 
the relative error for each ROM field is measured
against the corresponding FOM solution at the final time $\finalTime$, which is
defined as:
\begin{align}\label{eq:relerrors}
  \relerrorVelocityt{\finalTime} &= \frac{\| \velocityt{\ntimestep} -
  \velocityApproxt{\ntimestepROM} \|_2 }{\| \velocityt{\ntimestep} \|_2},& 
  \relerrorEnergyt{\finalTime} &= \frac{\| \energyt{\ntimestep} - 
  \energyApproxt{\ntimestepROM} \|_2 }{\| \energyt{\ntimestep} \|_2},& 
  \relerrorPositiont{\finalTime} &= \frac{\| \positiont{\ntimestep} -
  \positionApproxt{\ntimestepROM} \|_2 }{\| \positiont{\ntimestep} \|_2}.
\end{align}
The speed-up of each ROM simulation is measured by dividing the wall-clock time
for the FOM time loop by the wall-clock time for the corresponding ROM time
loop. 

\subsection{Effects of FOM discretization}\label{sec:exp-fom}
As a first experiment, we investigate the performance in terms of solution accuracy and speed-up of 
the ROM in the same problem set-up with different FOM discretization settings. 
The FOM discretization is set on different mesh refinement levels  
and the same finite element polynomial degree of 2 and 1 for the kinematic and thermodynamic space respectively, 
and hence the dimension of the finite element spaces 
asymptotically scales quadratically with the refinement level. 
The snapshots are taken at $\nparam = 1$ parameter, namely $\param_1 = 1/3$,
with the final time $\finalTime(\param_1) = 1.5$. 
Table~\ref{tab:exp1-fom} summarizes the degress of freedom, number of time steps, and wall clock time 
in Rayleigh--Taylor instability problem with various mesh refinement level of FOM discretization.
\begin{table}[ht!]
\centering
\begin{tabular}{|c||c|c|c|c|}
\hline
Refinement level & 2 & 3 & 4 & 5 \\ 
\hline
$\sizeKinematicFE$ & 594 & 2210 & 8514 & 33410 \\ 
$\sizeThermodynamicFE$ & 256 & 1024 & 4096 & 16384 \\ 
$\ntimestep$ & 435 & 926 & 1870 & 3780 \\ 
FOM wall clock (seconds) & 1.96 & 15.1 & 128 & 1010 \\ \hline
\end{tabular}
\caption{Degress of freedom, number of time steps, and wall clock time in 
  Rayleigh--Taylor instability problem with varying FOM discretization.}
\label{tab:exp1-fom}
\end{table}

In the offline phase, the fully discrete FOM scheme is first used to compute solution snapshots, 
which are classified into groups using the physical time as the indicator $\windowIndicator$ 
(see Section~\ref{sec:tw}). 
In practice, we start with $\windowIndicator_0 = 0$, and 
sequentially determine the indicator range partition in the offline phase using Algorithm~\ref{alg:decompose}. 
We remark that, when $\nparam = 1$, according to \eqref{eq:snapshot_cluster}, we have 
\begin{equation}
\snapshotGroup_{\windowIndex} 
= \{ (n,1) \in \mathbb{Z} \times \{1\}: 
\timeIndexWindow{\windowIndex-1}+1 \leq \timeIndex \leq 
\min\{\timeIndexWindow{\windowIndex-1}+\ntimestepWindow+1, \ntimestep\} \}, 
\end{equation}
where we drop the notation of the dependence on $\param_1$ for brevity. 
In this experiment, we take $\ntimestepWindow = 20$.
Then for $1 \leq \windowIndex \leq \nwindow^\text{off}$, 
the reduced order model is constructed by performing POD on the subset of snapshots 
whose indices belong to the group $\snapshotGroup_{\windowIndex}$. 
The oversampling ratio is taken as $\factorROMforceOneSample=\factorROMforceTvSample=2$ in hyper-reduction.
The offset vectors are taken as the initial state, i.e. 
$\velocityOSWindow{\windowIndex} = \velocity(0)$, 
$\energyOSWindow{\windowIndex} = \energy(0)$,  and 
$\positionOSWindow{\windowIndex} = \position(0)$. 
The resultant reduced solution bases for velocity and energy are left-multiplied by
the corresponding mass matrix to obtain the nonlinear term bases using SNS, which in turn
determine the sampling indices by DEIM (see Section~\ref{sec:deim}).

To avoid complicating the results, we fix all ROM parameters
and consider only the reproductive case, i.e. 
the ROM simulation is performed at the same parameter $\param = 1/3$
and the final time $\finalTimeROM(\param) = 1.5$ in the online phase. 
To simplify the notation, we suppress the parameter dependence 
of all the notations in the remainder of this subsection, without introducing any ambiguity. 
In the online phase, the corresponding reduced order model is used, 
where the end-point $\windowk{\windowIndex}(\param)$ of the temporal subinterval $\timeWindow_{\windowIndex}(\param)$ 
is defined as the time instance $\timekROM{\timeIndex}(\param)$
when $\windowIndicator(\stateApproxt{\timeIndex}(\param), \timekROM{\timeIndex}(\param), \param)$ 
first exceeds $\windowIndicator_{\windowIndex}$, 
at which we increment to the next time subinterval $\timeWindow_{\windowIndex+1}(\param)$.
The reduced basis and the sampling matrices are then used to formulate the
hyper-reduced system \eqref{eq:rom-hr-tw} in each temporal subtinterval. The fully
discrete hyper-reduced system follows from applying the RK2-average scheme.  

Table~\ref{tab:exp1-tw} and Table~\ref{tab:exp1-dw} 
summarize the results in solution accuracy and speed-up of TW-ROM in Section~\ref{sec:tw} 
and DW-ROM in Section~\ref{sec:dw}, respectively. 
In the short-time simulation, both TW-ROM and DW-ROM perform well, and 
they have similar level of solution accuaracy and speed-up.
As the refinement level increases, the solution accuracy in all the variables remains outstanding. 
The relative error ranges from $O(10^{-7})$ to $O(10^{-3})$ and is not sensitive to the FOM discretization. 
Among the three solution components, the velocity has the largest relative error at the final time. 
In fact, the absolute error of the velocity and the specific internal energy has the same order of magnitude, 
but the order of the norm of the velocity is lower than that of the energy by 2.
Meanwhile, the speed-up increases with the refinement level. 
While the FOM computational cost scales as a polynomial with the refinement level, 
the ROM computational expense mostly depends on the distribution of snapshot data and the ROM parameters, 
and is much less sensitive to the FOM discretization. 
The speed-up is very remarkable with higher refinement levels, for which the FOM simulations are more expensive.
With the refinement level of 5, the speed-up is almost 50 times with both TW-ROM and DW-ROM.

\begin{table}[ht!]
\centering
\begin{tabular}{|c||c|c|c|c|}
\hline
Refinement level & 2 & 3 & 4 & 5 \\ 
\hline
$\nwindow$ & 44 & 93 & 187 & 379 \\
$\ntimestepROM$ & 403 & 848 & 1841 & 3502 \\ 
$\relerrorVelocityt{\finalTime}$ & 4.2217e-03 & 4.5689e-03 & 4.1547e-02 & 2.8951e-03 \\
$\relerrorEnergyt{\finalTime}$ & 3.8508e-06 & 1.0837e-05 & 1.4926e-04 & 1.0660e-05 \\ 
$\relerrorPositiont{\finalTime}$ & 9.4068e-06 & 9.1305e-07 & 2.7141e-04 & 2.7827e-05 \\ 
speed-up & 1.1058 & 3.4887 & 12.5000 & 48.6021 \\  \hline
\end{tabular}
\caption{TW-ROM performance comparison for short-time simulation in 
  Rayleigh--Taylor instability problem with varying FOM discretization. 
  }
\label{tab:exp1-tw}
\end{table}

\begin{table}[ht!]
\centering
\begin{tabular}{|c||c|c|c|c|}
\hline
Refinement level & 2 & 3 & 4 & 5 \\ 
\hline 
$\nwindow$ & 44 & 93 & 187 & 379 \\
$\ntimestepROM$ & 403 & 848 & 1752 & 3584 \\ 
$\relerrorVelocityt{\finalTime}$ & 4.2217e-03 & 4.6212e-03 & 3.3420e-03 & 2.4819e-03 \\
$\relerrorEnergyt{\finalTime}$ & 3.8508e-06 & 1.0990e-05 & 8.6944e-06 & 5.7578e-06 \\
$\relerrorPositiont{\finalTime}$ & 9.4068e-06 & 1.3922e-06 & 1.7685e-05 & 1.7119e-05 \\
speed-up & 1.0927 & 3.5467 & 12.0573 & 46.0935 \\  \hline
\end{tabular}
\caption{DW-ROM performance comparison for short-time simulation in 
  Rayleigh--Taylor instability problem with varying FOM discretization. 
  }
\label{tab:exp1-dw}
\end{table}

\subsection{Comparison of partition indicators for extrapolation}\label{sec:exp-atwood}
In this section, we investigate the performance in terms of solution accuracy and speed-up of 
the ROM with the different partition indicators at various Atwood numbers in the parametric problem setting. 
It should be emphasized that the parametric problem setting is challenging as 
the initial state depends on the density ratio $R$ and hence the Atwood number $\param$, and has a significant 
influence on the highly nonlinear dynamics.
The FOM discretization is set on the same mesh refinement level of 4
and the same finite element polynomial degree of 2 and 1 for the kinematic and thermodynamic space respectively. 
We follow the same procedure as in Section~\ref{sec:exp-fom} for sampling snapshots, decomposing the indicator range, 
constructing bases in the offline phase, and ROM simulation in the online phase.

\begin{figure}[ht!]
\centering
\includegraphics[width=0.45\linewidth]{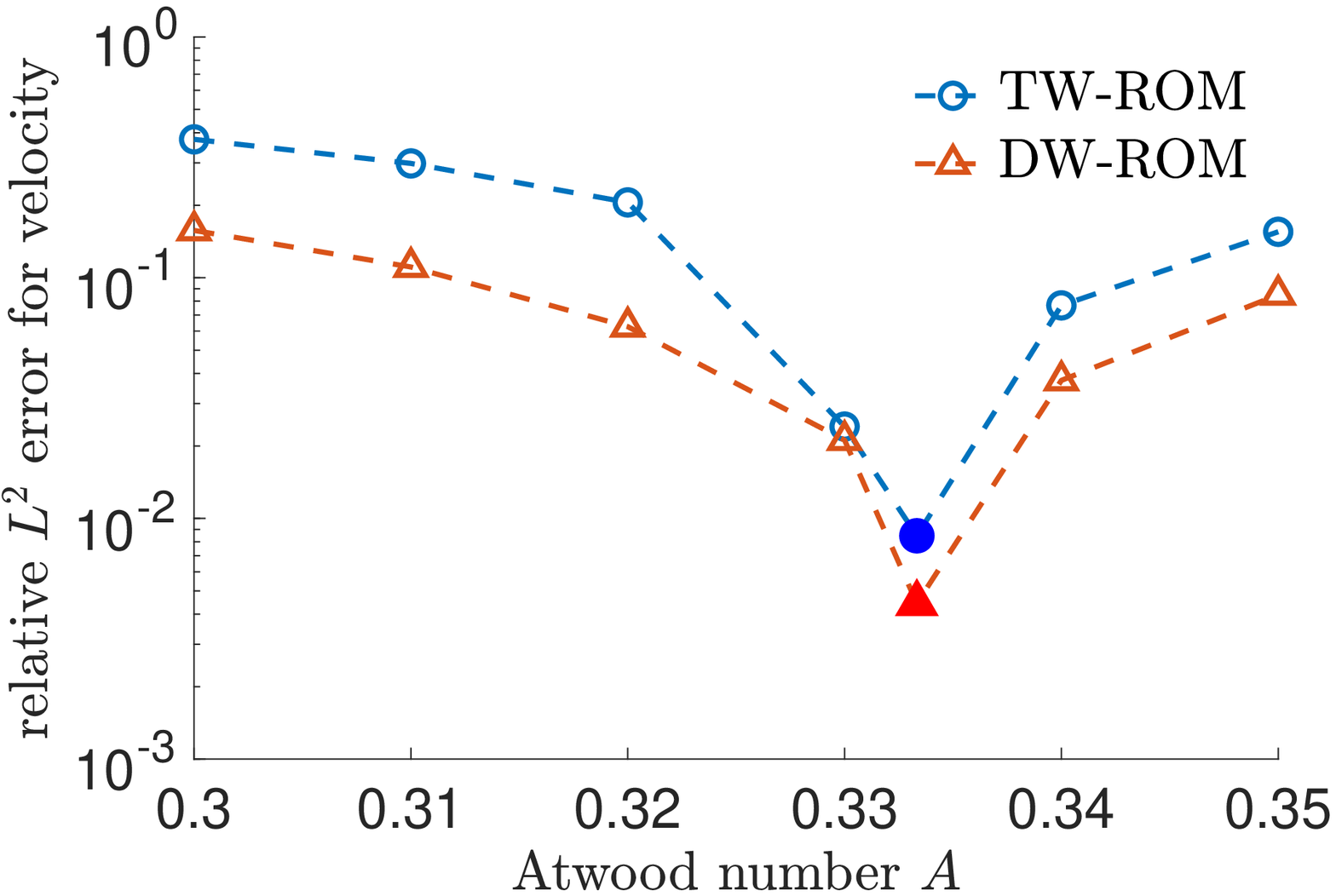}
\hspace{0.05\linewidth}
\includegraphics[width=0.45\linewidth]{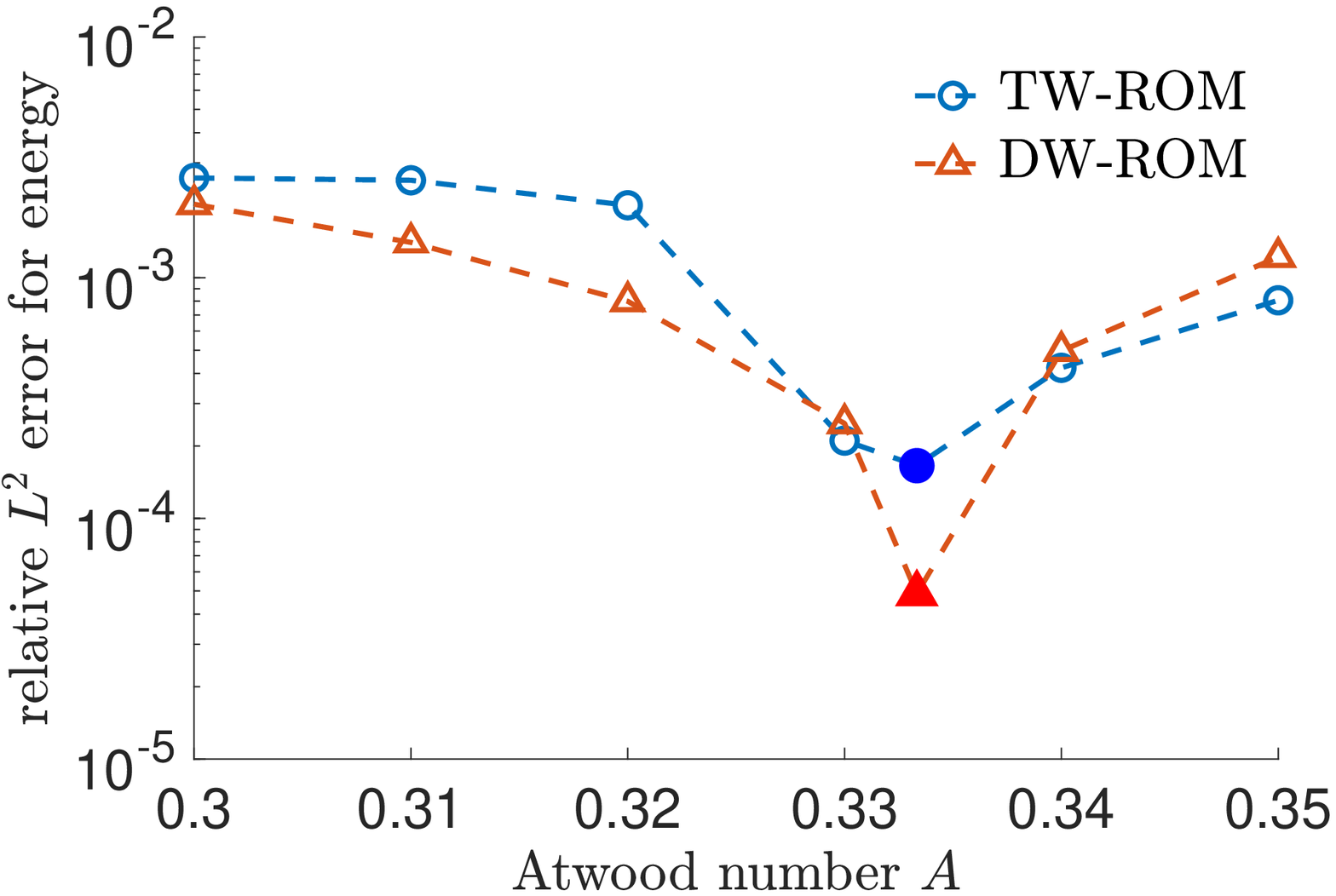}
\vspace{0.2in}
\includegraphics[width=0.45\linewidth]{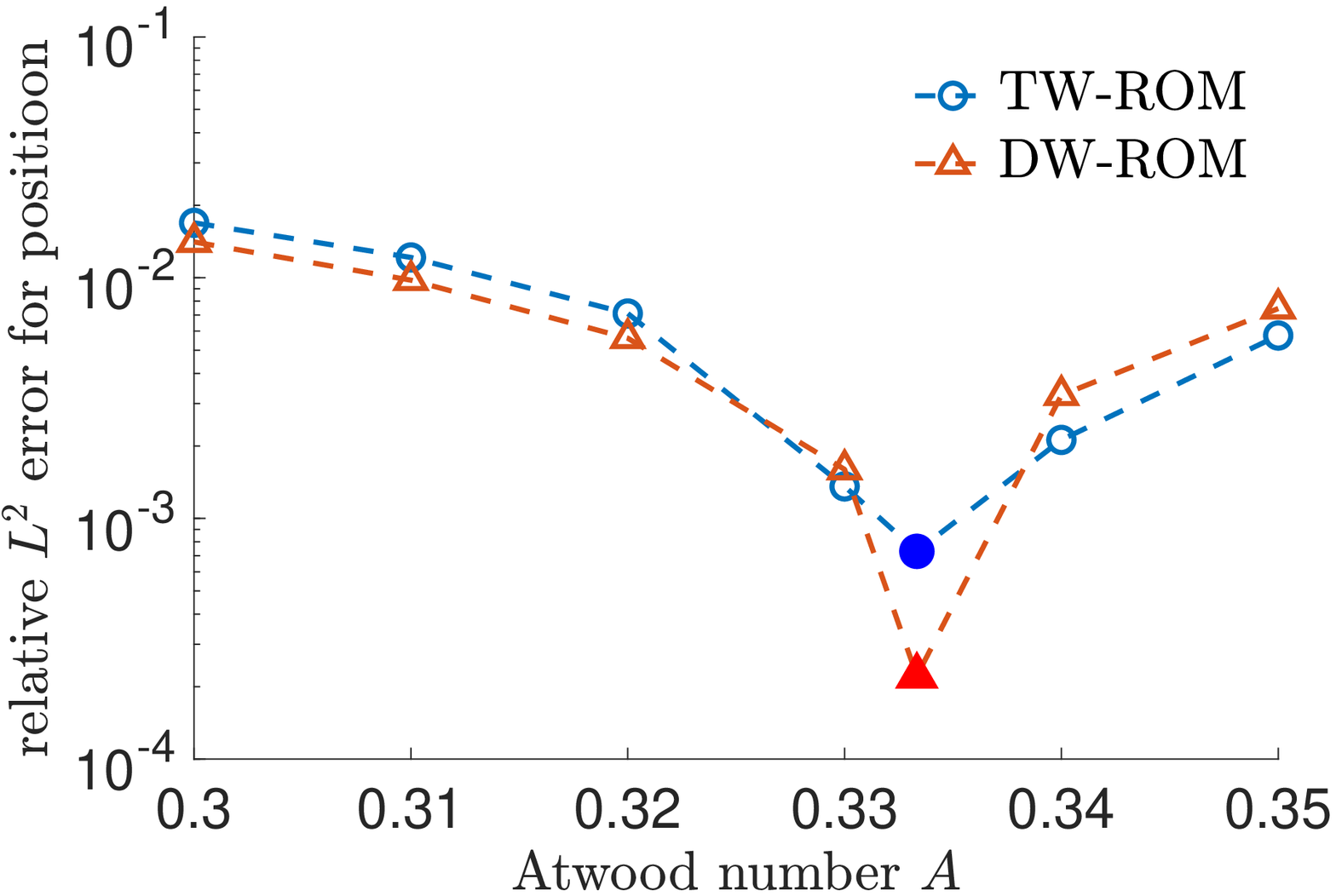}
\hspace{0.05\linewidth}
\includegraphics[width=0.45\linewidth]{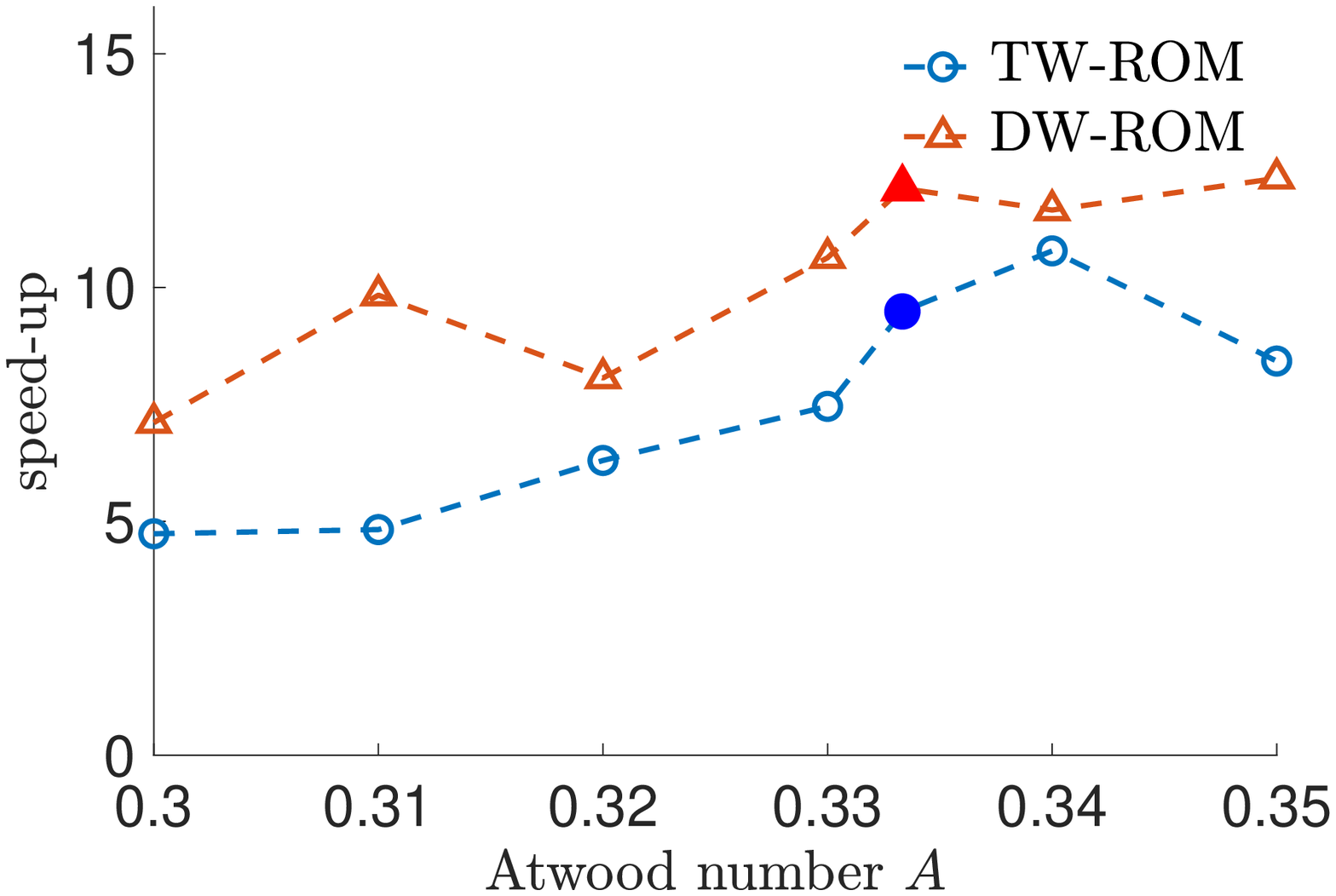}
  \caption{{\bf Extrapolation test}: ROM performance comparison for long-time simulation in 
Rayleigh--Taylor instability problem at varying Atwood number with $\nparam = 1$ and $\singularValueThreshold = 10^{-4}$:
relative $L^2$ error for velocity (top-left), 
relative $L^2$ error for energy (top-right), 
relative $L^2$ error for position (bottom-left), and 
speed-up (bottom-right). 
Using the penetration distance improves the solution accuracy at the extrapolatory cases.
}
\label{fig:mid_par1_ef4}
\end{figure}

\begin{figure}[ht!]
\centering
\includegraphics[width=0.45\linewidth]{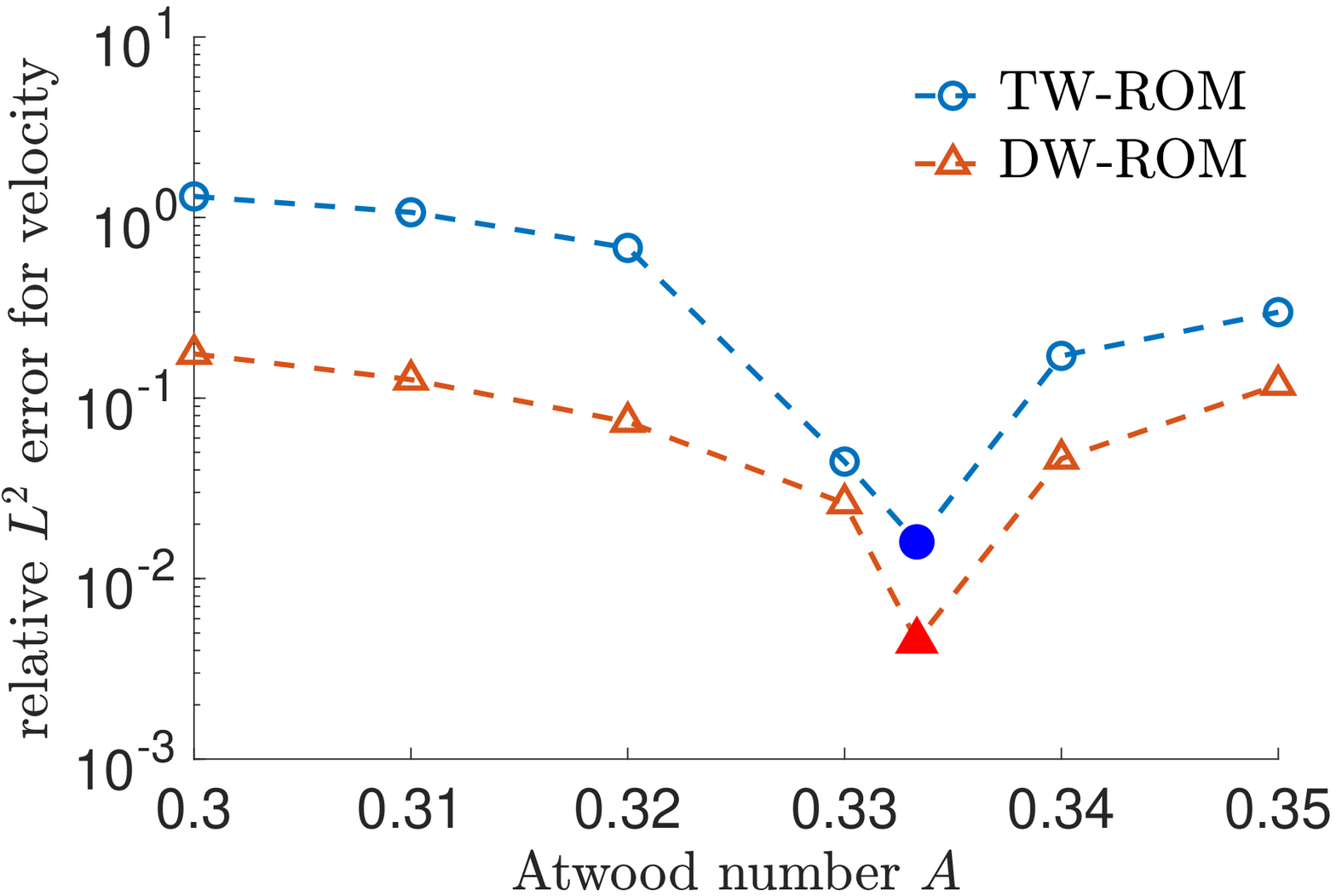}
\hspace{0.05\linewidth}
\includegraphics[width=0.45\linewidth]{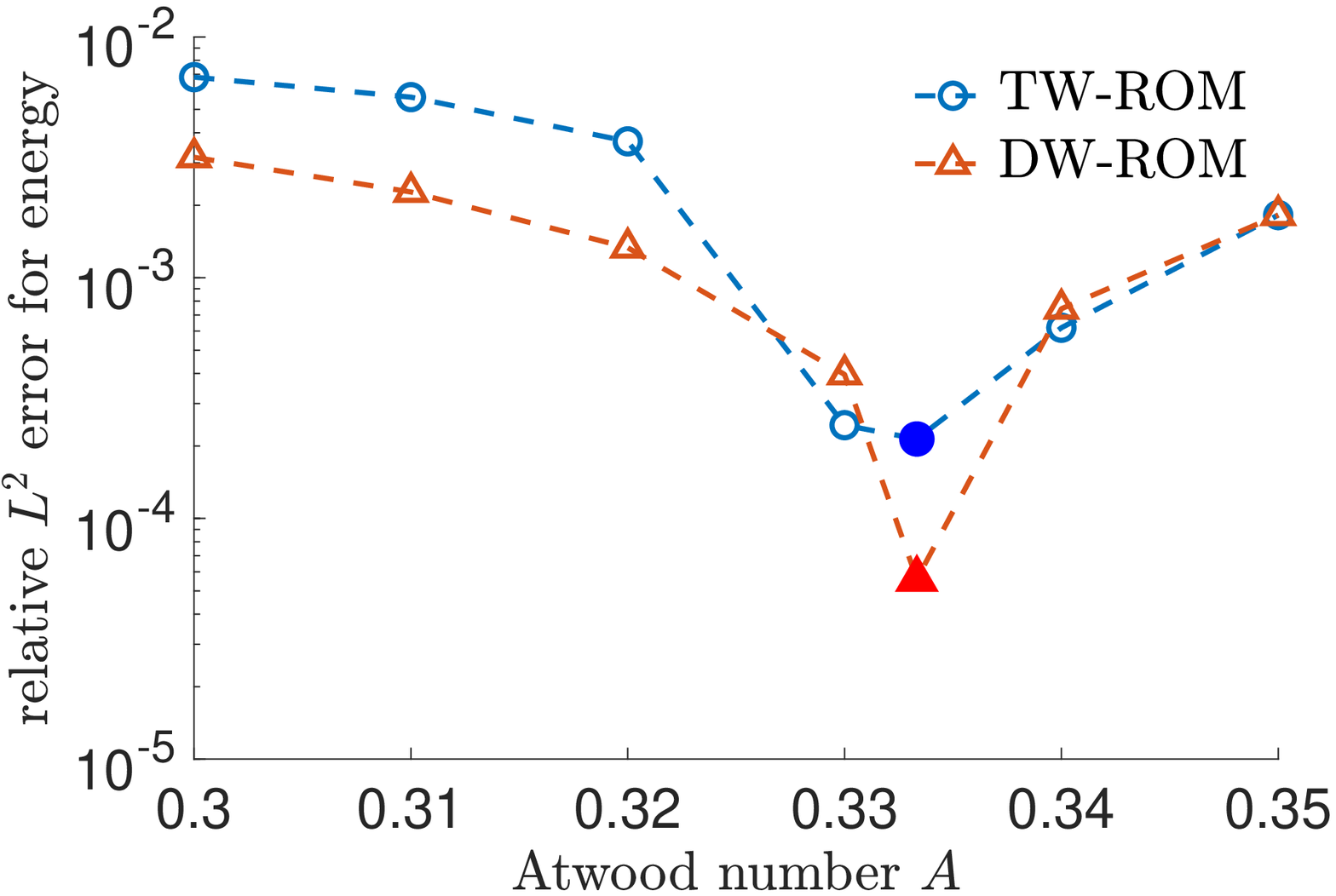}
\vspace{0.2in}
\includegraphics[width=0.45\linewidth]{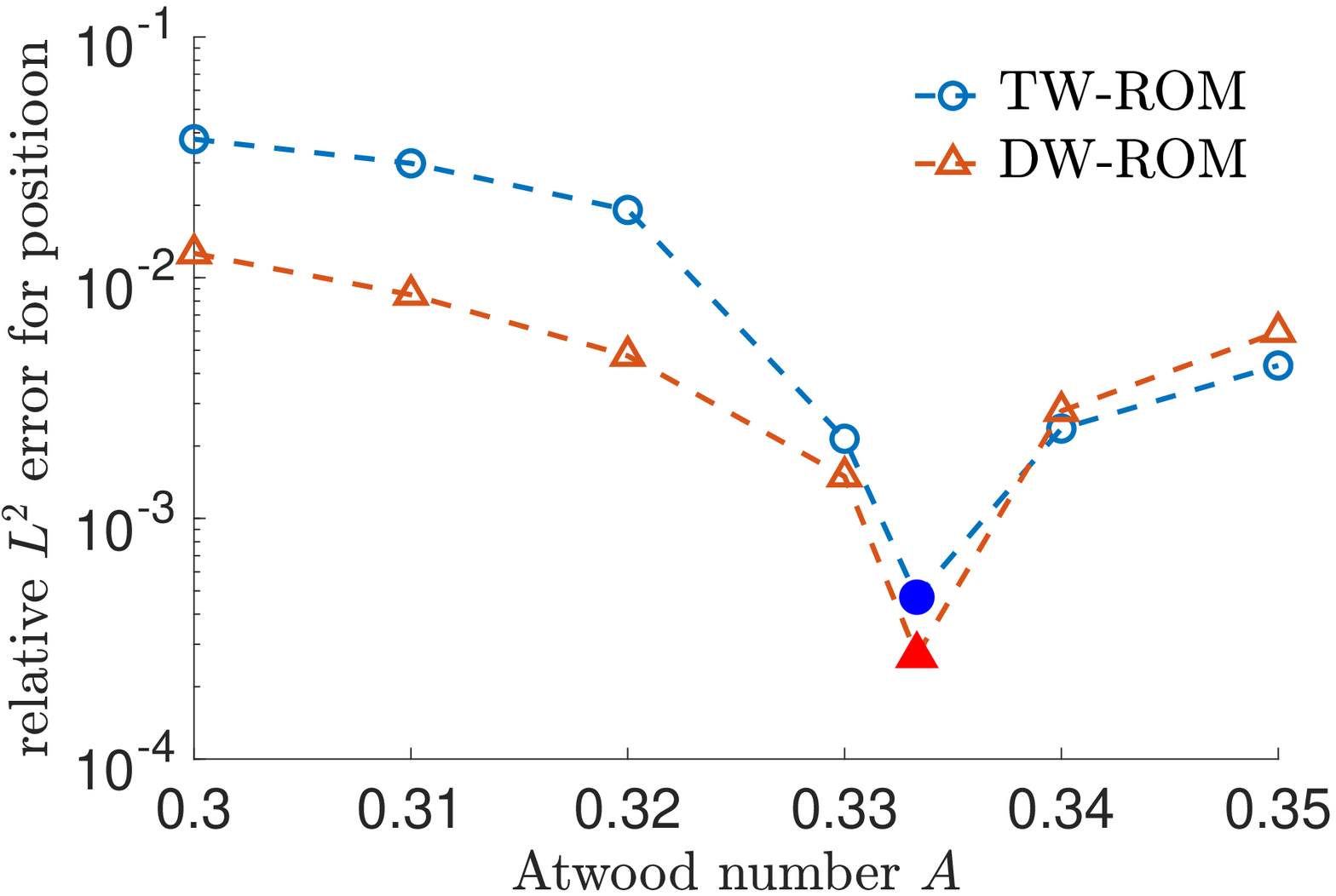}
\hspace{0.05\linewidth}
\includegraphics[width=0.45\linewidth]{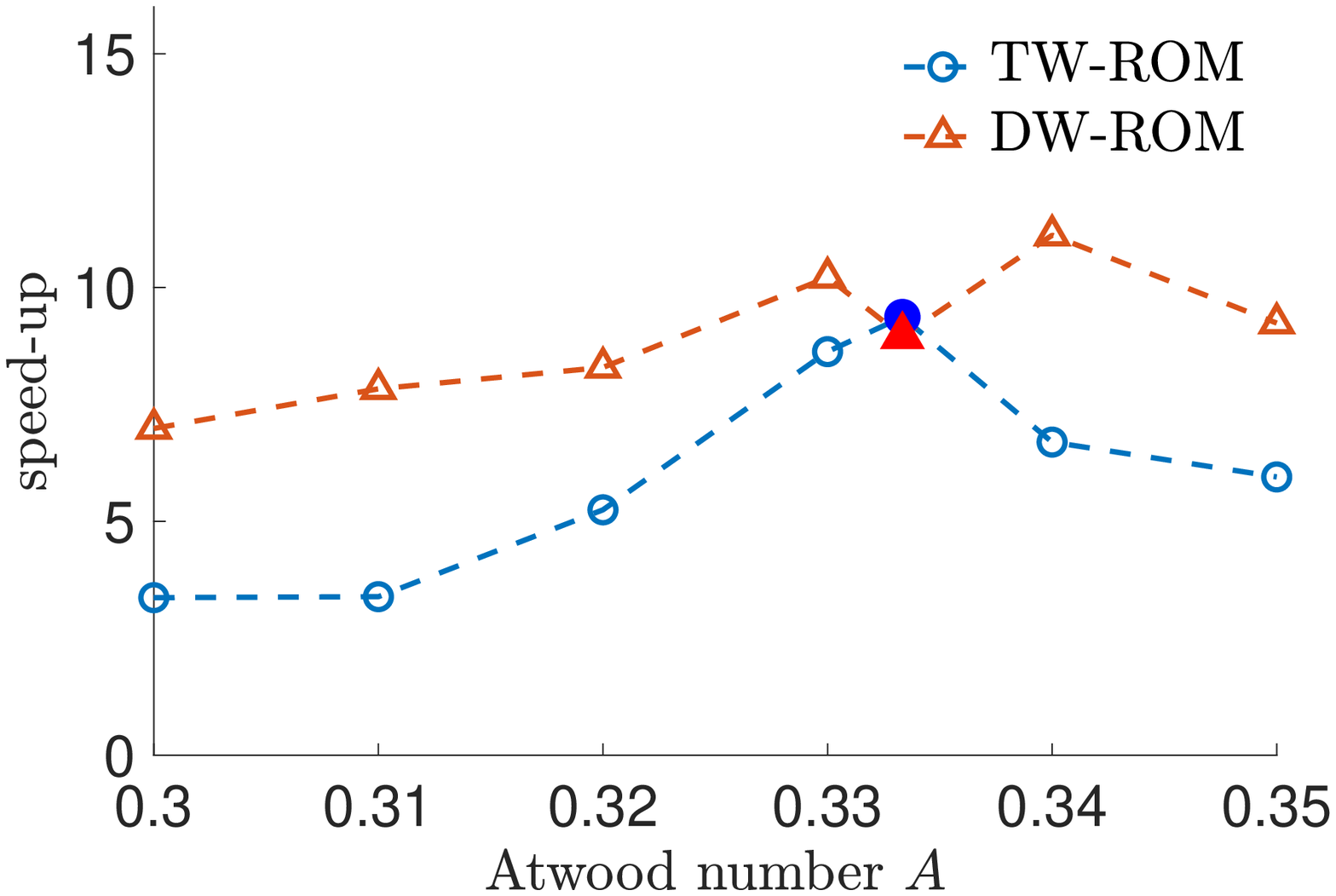}
  \caption{{\bf Extrapolation test}: ROM performance comparison for long-time simulation in 
Rayleigh--Taylor instability problem at varying Atwood number with $\nparam = 1$ and $\singularValueThreshold = 10^{-10}$:
relative $L^2$ error for velocity (top-left), 
relative $L^2$ error for energy (top-right), 
relative $L^2$ error for position (bottom-left), and 
speed-up (bottom-right). 
Using the penetration distance improves the solution accuracy at the extrapolatory cases.
}
\label{fig:mid_par1_ef10}
\end{figure}

In the first experiment, we use $\nparam = 1$ training parameter, namely $\param_1 = 1/3$,
with the final time $\finalTime(\param_1) = 3.2$ for snapshot sampling. 
We take $\ntimestepWindow = 20$ for decomposing the indicator range by Algorithm~\ref{alg:decompose}, 
and $\singularValueThreshold \in \{ 10^{-4}, 10^{-10} \}$ for performing proper orthogonal decomposition. 
The ROM simulation is performed at the various parameters 
$\param \in \paramDomain = [0.3, 0.35]$ with final time $\finalTimeROM(\param) \equiv 3$
and the oversampling ratio $\factorROMforceOneSample=\factorROMforceTvSample=5$ for hyper-reduction in the online phase. 
Figure~\ref{fig:mid_par1_ef4} and Figure~\ref{fig:mid_par1_ef10} show the results in solution accuracy and speed-up 
with  $\singularValueThreshold = 10^{-4}$ and $\singularValueThreshold = 10^{-10}$ respectively.
The filled markers correspond to the reproductive case, i.e. $\param = \param_1 = 1/3$,
while the other markers correspond to {\it extrapolatory} cases at selected testing points 
$\param \in \paramDomain \setminus \{\param_1\}$. 
The results of TW-ROM, i.e. using physical time as the indicator in Section~\ref{sec:tw}, 
and DW-ROM, i.e. using penetration distance as the indicator in Section~\ref{sec:dw}, 
are depicted in blue and red, respectively.
From the error plots for all the solution fields, 
it can be observed that the solution accuracy degenerates when the testing points 
are farther away from the training point $\param_1 = 1/3$, in both cases of using 
the physical time or penetration distance as the indicator $\windowIndicator$. 
Both approaches provide an acceptable solution accuracy at the reproductive case $\param = 1/3$, 
and using the penetration distance gives better results, where the error in velocity is less than 1\%. 
However, the approaches have significantly different performance at the extrapolatory cases 
which are of greater interest in the practical use of reduced order models. 
The numerical results show that the physical time indicator generally has worse accuracy, 
especially for the velocity field at the extrapolatory cases, 
where the error at $\param \in [0.3, 0.32]$ are higher than 40\% with $\singularValueThreshold = 10^{-4}$ 
as shown in Figure~\ref{fig:mid_par1_ef4}, 
and almost 100\% with $\singularValueThreshold = 10^{-10}$ 
as shown in Figure~\ref{fig:mid_par1_ef10}.
On the other hand, the penetration distance is a more reliable indicator to 
assign the reduced order model in the online phase for extrapolatory parameters in general, 
especially for lower Atwood numbers $\param \in [0.3, 0.32]$. 
Moreover, the penetration distance indicator is superior to the physical time indicator in terms of speed-up. 

\begin{figure}[ht!]
\centering
\includegraphics[width=0.45\linewidth]{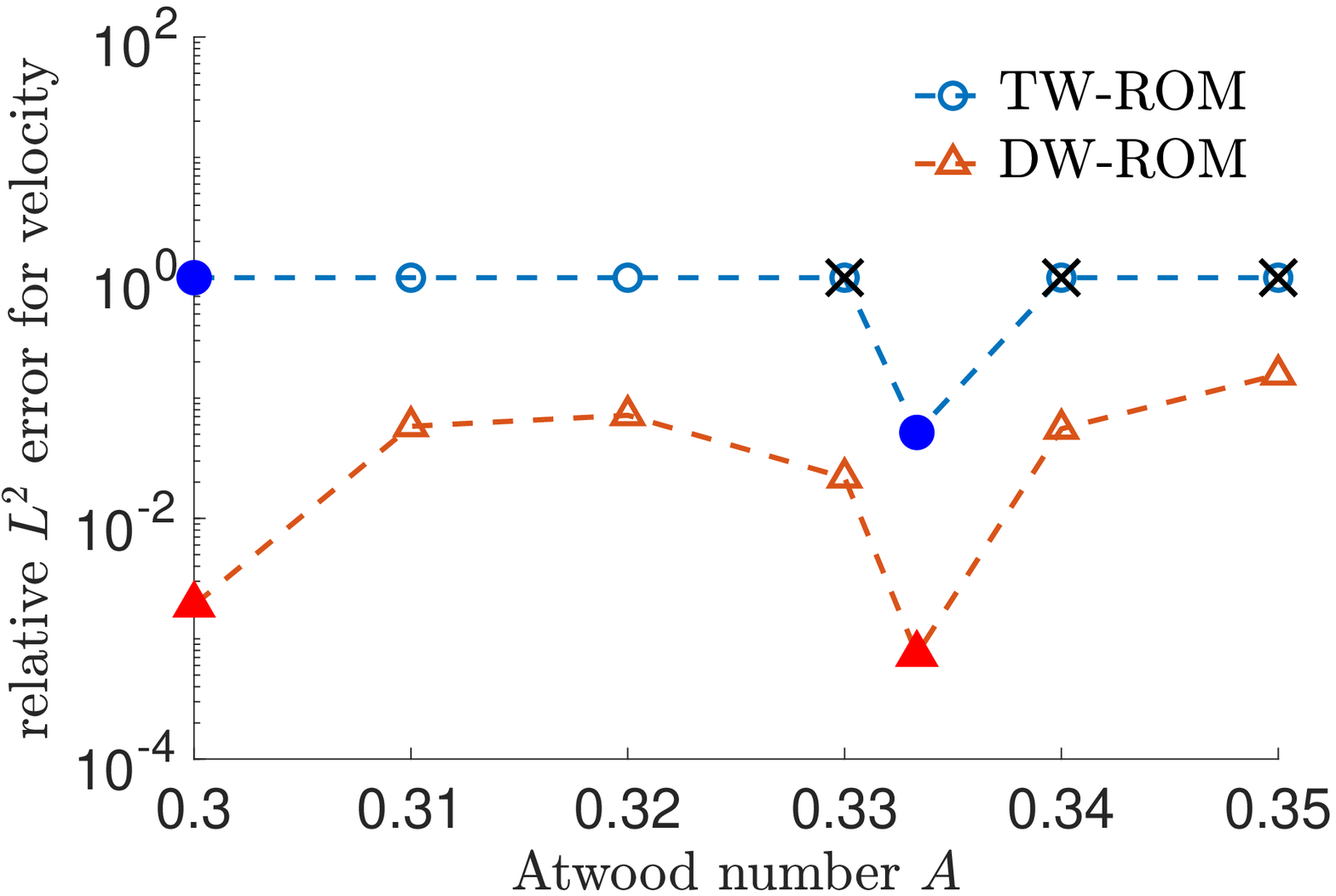}
\hspace{0.05\linewidth}
\includegraphics[width=0.45\linewidth]{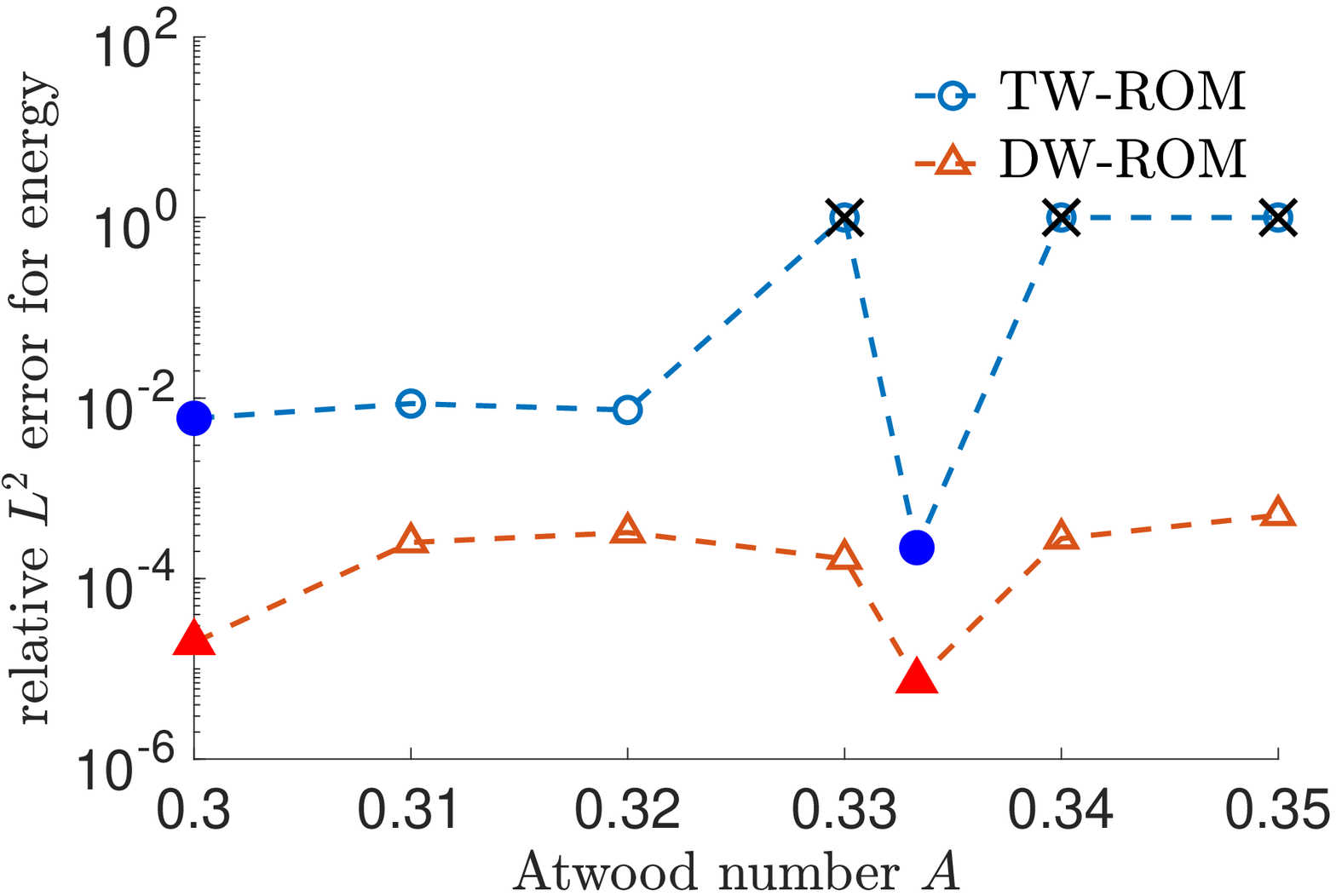}
\vspace{0.2in}
\includegraphics[width=0.45\linewidth]{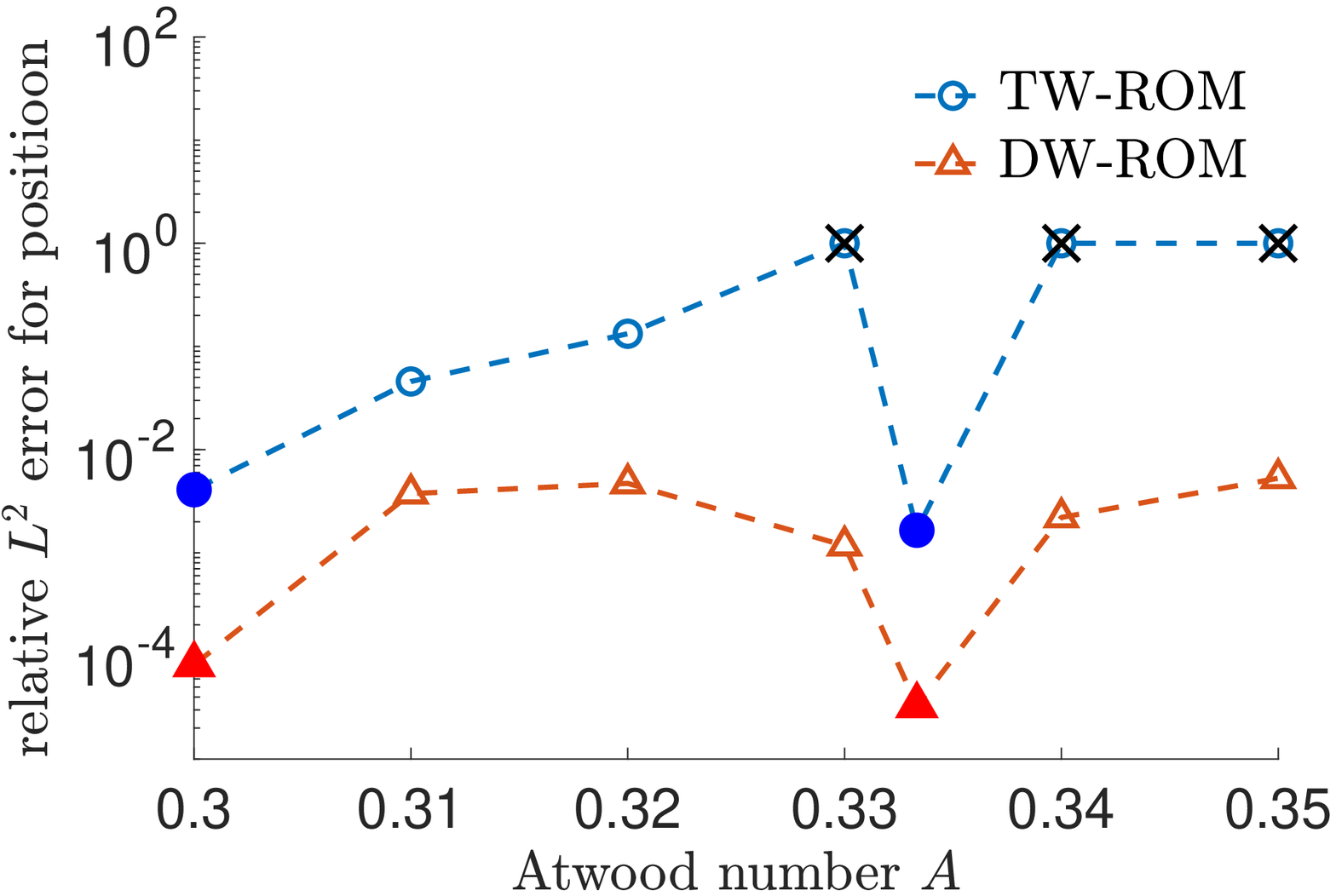}
\hspace{0.05\linewidth}
\includegraphics[width=0.45\linewidth]{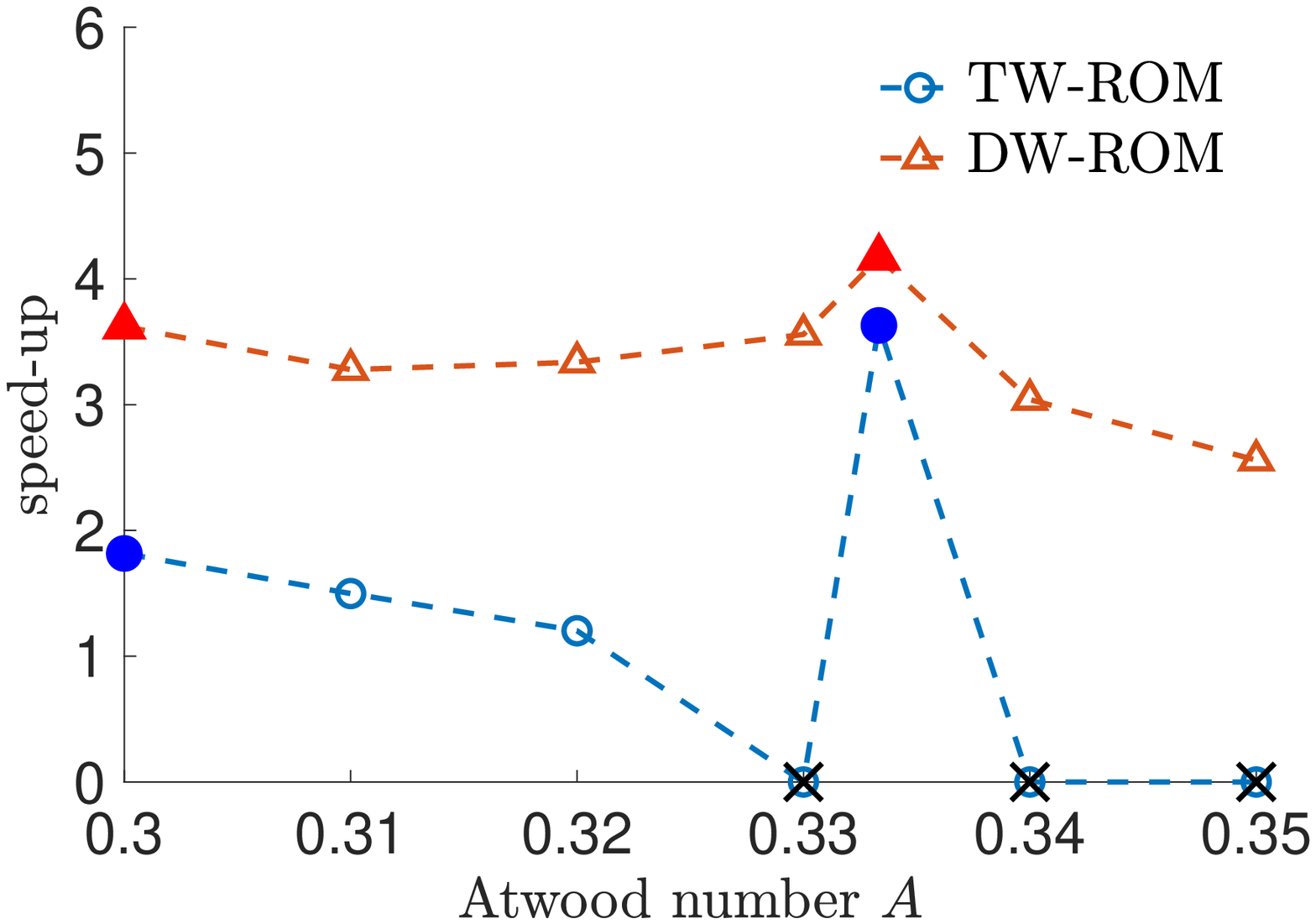}
\caption{ROM performance comparison for long-time simulation in 
Rayleigh--Taylor instability problem at varying Atwood number with $\nparam = 2$ and $\singularValueThreshold = 10^{-10}$:
relative $L^2$ error for velocity (top-left), 
relative $L^2$ error for energy (top-right), 
relative $L^2$ error for position (bottom-left), and 
speed-up (bottom-right). 
Using the penetration distance improves the solution accuracy and stability at the extrapolatory cases.
}
\label{fig:mid_par2}
\end{figure}

In the second experiment, we use $\nparam = 2$ training parameters, namely $\param_1 = 1/3$ and $\param_2 = 0.3$,
with the final time $\finalTime(\param_1)  = \finalTime(\param_2) = 3.2$ for snapshot sampling. 
We take $\ntimestepWindow = 20$ for decomposing the indicator range by Algorithm~\ref{alg:decompose}, 
and $\singularValueThreshold = 10^{-10}$ for performing proper orthogonal decomposition. 
The ROM simulation is performed at the various parameters 
$\param \in \paramDomain = [0.3, 0.35]$ with final time $\finalTimeROM(\param) \equiv 3$
and the oversampling ratio $\factorROMforceOneSample=\factorROMforceTvSample=10$ for hyper-reduction in the online phase. 
Figure~\ref{fig:mid_par2} shows the solution accuracy and speed-up. 
The filled markers correspond to the reproductive case, i.e. $\param \in \{\param_1, \param_2\}$,  
while the other markers correspond to extrapolatory cases at selected testing points 
$\param \in \paramDomain \setminus \{\param_1, \param_2\}$
that do not contribute snapshot samples to the basis construction in the offline phase. 
The results of TW-ROM, i.e. using physical time as the indicator in Section~\ref{sec:tw}, 
and DW-ROM, i.e. using penetration distance as the indicator in Section~\ref{sec:dw}, 
are depicted in blue and red, respectively. 
In this experiment, using the physical time as the indicator yields numerical instability, 
which is due to poor approximation properties of the reduced subspaces. 
It is therefore evidential that the decomposition of the solution manifold by the physical time 
does not produce submanifolds with small Kolmogorov $n$-width. 
The unstable cases at $\param \in \{0.33,0.34,0.35\}$ are indicated by the crossed markers 
with a 100\% error and zero speed-up in the blue line plot corresponding to TW-ROM. 
Even at the stable cases, the error in velocity is almost 100\%, which again suggests that 
using the physical time does not provide good approximation for the velocity.
Meanwhile, using the penetration distance as the indicator produces reasonable results. 
The error in the solution attains local minima at the training parameters $\{\param_1, \param_2\}$, 
and grows as the testing parameter $\param$ gets farther away from the training parameters.

\subsection{Comparison of partition indicators on final time}\label{sec:exp-final-time}

In this experiment, we investigate the performance of 
the ROM with the different partition indicators at different time instances in an extrapolatory problem setting. 
The final time is chosen to be sufficiently long 
for the instability to develop into a highly nonlinear regime with the formation of a vortex at the fluid interface. 
The FOM discretization is set on the same mesh refinement level of 4
and the same finite element polynomial degree of 2 and 1 for the kinematic and thermodynamic spaces, respectively.
The procedure of building and applying the reduced order model follows exactly as in Section~\ref{sec:exp-fom}. 
We use $\nparam = 1$ training parameter, namely $\param_1 = 1/3$,
with the final time $\finalTime(\param_1) = 3.6$ for snapshot sampling. 
We take $\ntimestepWindow = 20$ for decomposing the indicator range by Algorithm~\ref{alg:decompose}. 
The ROM simulation is performed at the unseen parameter $\param = 0.33$ 
with various final time $\finalTimeROM(\param) \in [3.0, 3.6]$ 
and the oversampling ratio $\factorROMforceOneSample=\factorROMforceTvSample=2$ for hyper-reduction in the online phase. 
To avoid complicating the results, we fix all ROM parameters. 

\begin{figure}[ht!]
\centering
\includegraphics[width=0.45\linewidth]{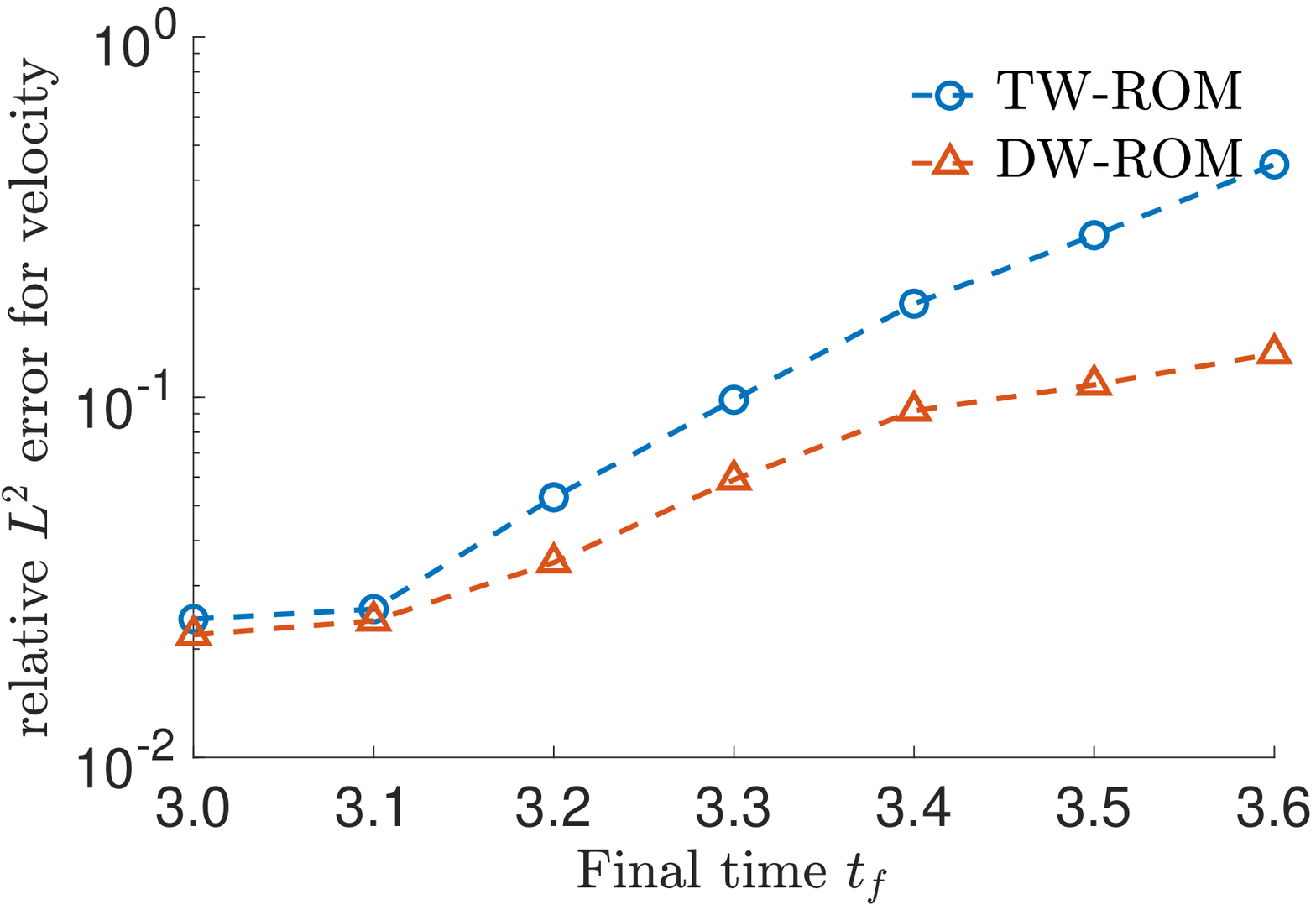}
\hspace{0.05\linewidth}
\includegraphics[width=0.45\linewidth]{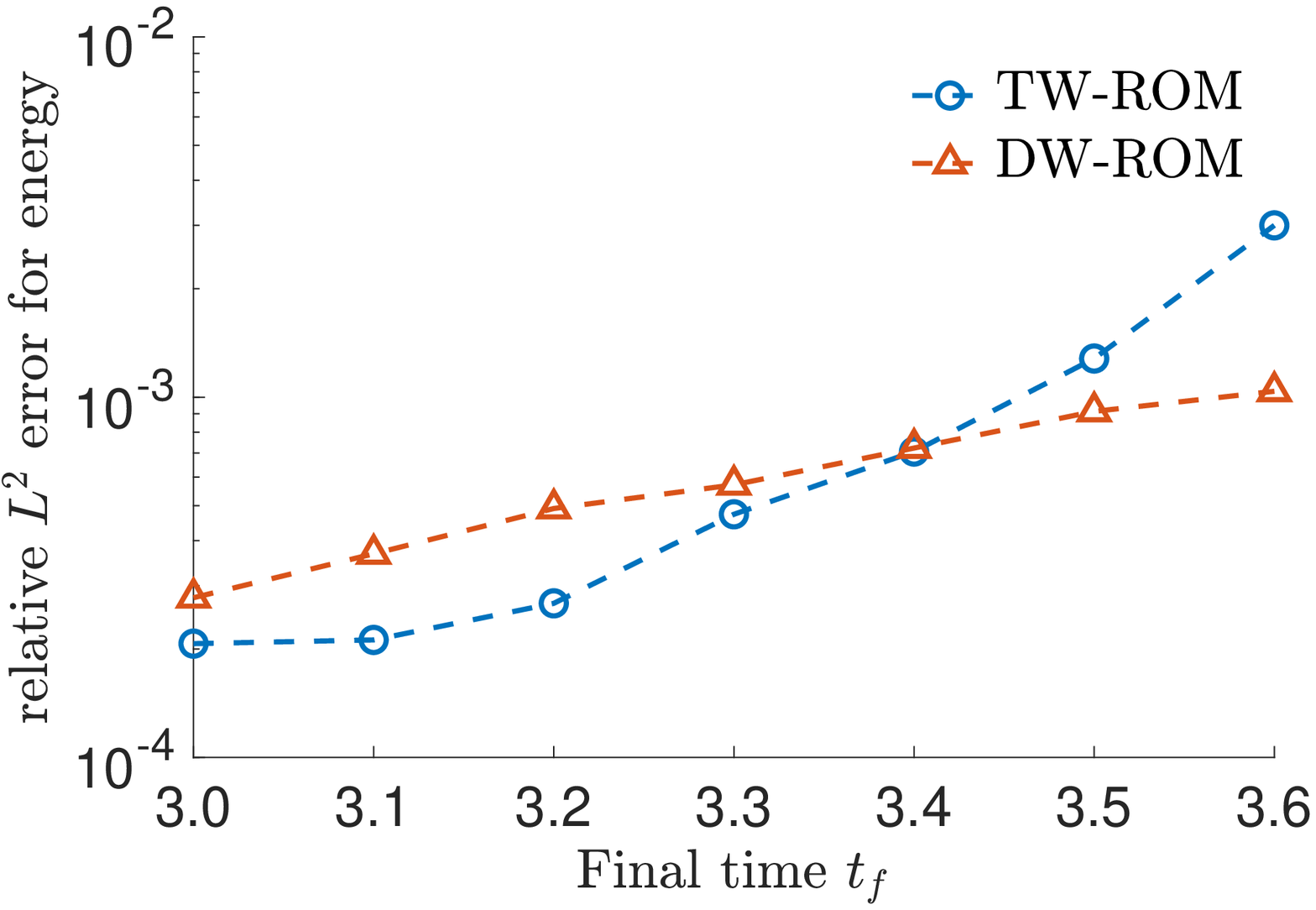}
\vspace{0.2in}
\includegraphics[width=0.45\linewidth]{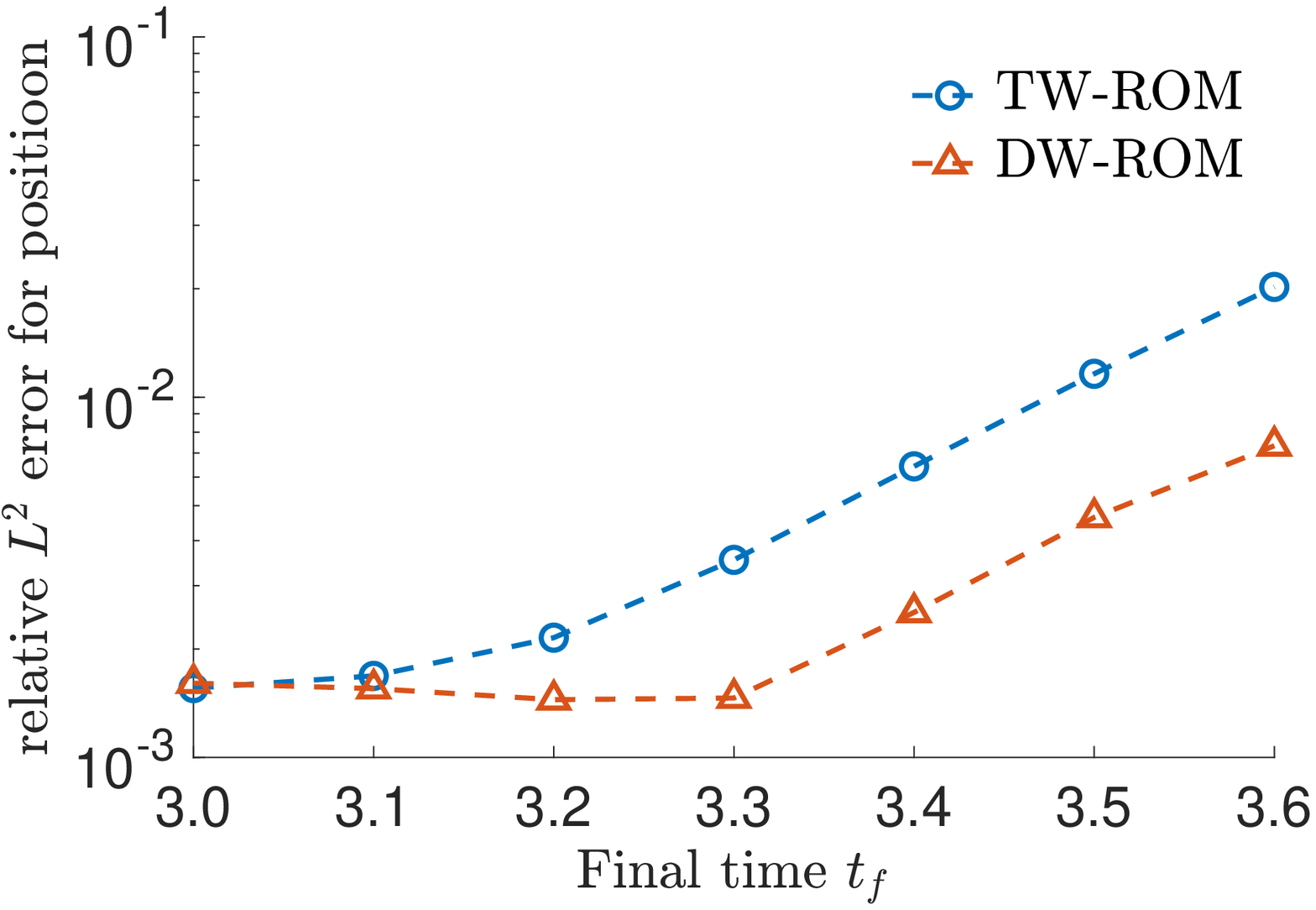}
\hspace{0.05\linewidth}
\includegraphics[width=0.45\linewidth]{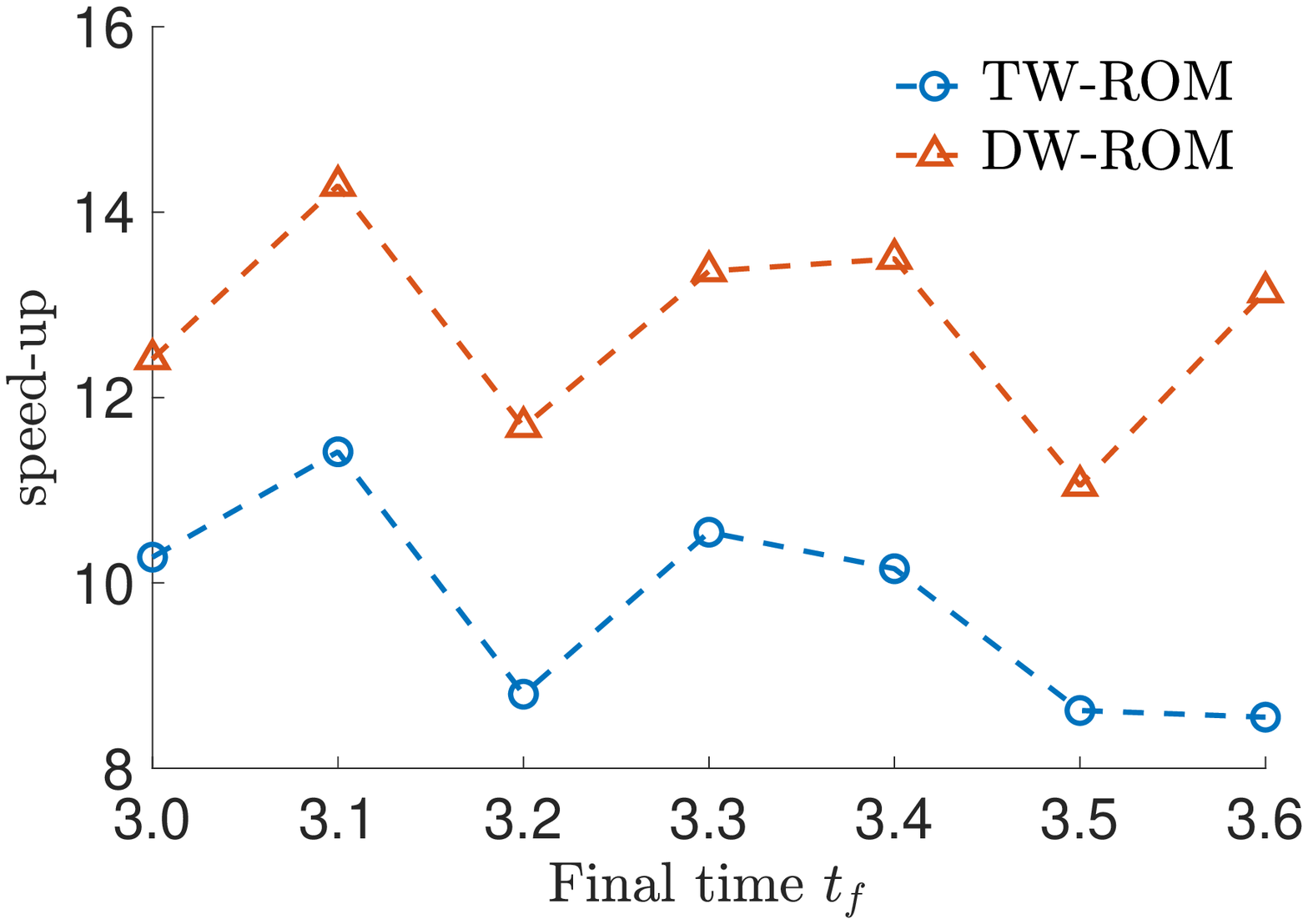}
\caption{ROM performance comparison for long-time simulation in 
  Rayleigh--Taylor instability problem with varying final time: 
relative $L^2$ error for velocity (top-left), 
relative $L^2$ error for energy (top-right), 
relative $L^2$ error for position (bottom-left), and 
speed-up (bottom-right). 
Using the penetration distance improves the solution accuracy 
and speed-up at all tested oversampling ratio.}
\label{fig:reuslts-final-time}
\end{figure}

Figure~\ref{fig:reuslts-final-time} shows the results in terms of solution accuracy and speed-up. 
The results of TW-ROM, i.e. using physical time as the indicator in Section~\ref{sec:tw}, 
and DW-ROM, i.e. using penetration distance as the indicator in Section~\ref{sec:dw}, 
are depicted in blue and red, respectively. 
It can be seen that the solution error accumulates with time in both decomposition mechanisms. 
At the final time $\finalTime = 3.6$, the error in velocity is around 44\% when the physical time is used as the indicator, 
while it is around 13\% when the penetration distance is used. 
These results highlight the importance of an appropriate indicator 
for accurate local ROM approximation, especially in highly complex nonlinear dynamics, 
such as the long-time simulation of Rayleigh--Taylor instability with the formation of vortices. 
As we will see in Section~\ref{sec:exp-oversampling}, the solution accuracy can be 
further improved by adding more samples in hyper-reduction when the penetration distance is used as the indicator.  
It is also important to point out that 
when the penetration distance is used as the indicator, a speed-up of 12-15 times is achieved at varying time instances, 
which are higher than 8-12 times when the physical time is used. 

\subsection{Comparison of partition indicators on oversampling ratio}\label{sec:exp-oversampling}

In this experiment, we further investigate the performance of 
the ROM with the different partition indicators and various oversampling ratios of hyper-reduction in an extrapolatory problem setting. 
The final time is set to be sufficiently long to allow the formation of a vortex 
at the fluid interface in the highly nonlinear dynamics. 
The FOM discretization is set on the same mesh refinement level of 4
and the same finite element polynomial degree of 2 and 1 for the kinematic and thermodynamic spaces, respectively. 
The procedure of building and applying the reduced order model follows exactly as in Section~\ref{sec:exp-fom}. 
We use $\nparam = 1$ training parameter, namely $\param_1 = 1/3$,
with the final time $\finalTime(\param_1) = 3.6$ for snapshot sampling. 

\begin{figure}[ht!]
\centering
\includegraphics[width=0.45\linewidth]{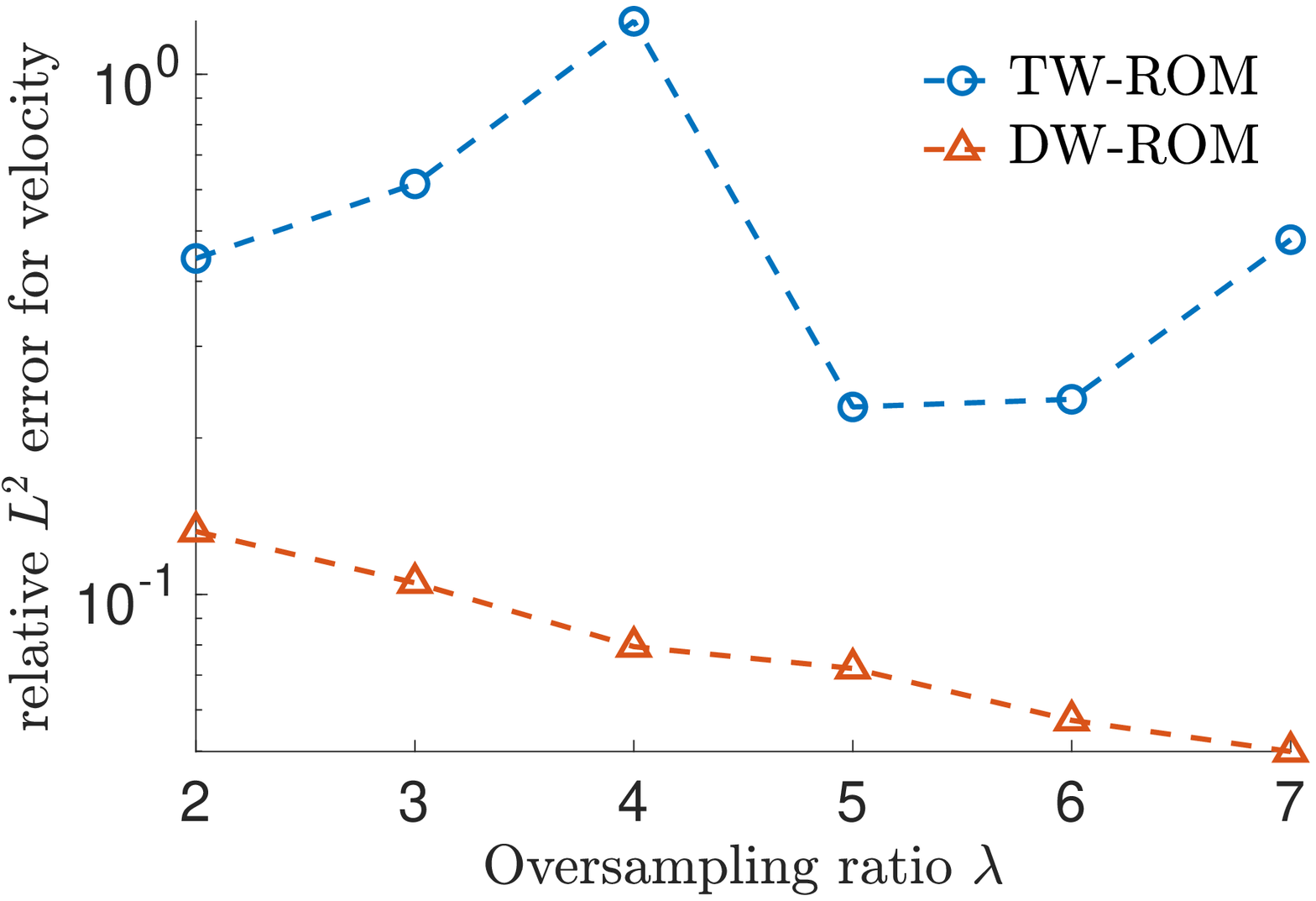}
\hspace{0.05\linewidth}
\includegraphics[width=0.45\linewidth]{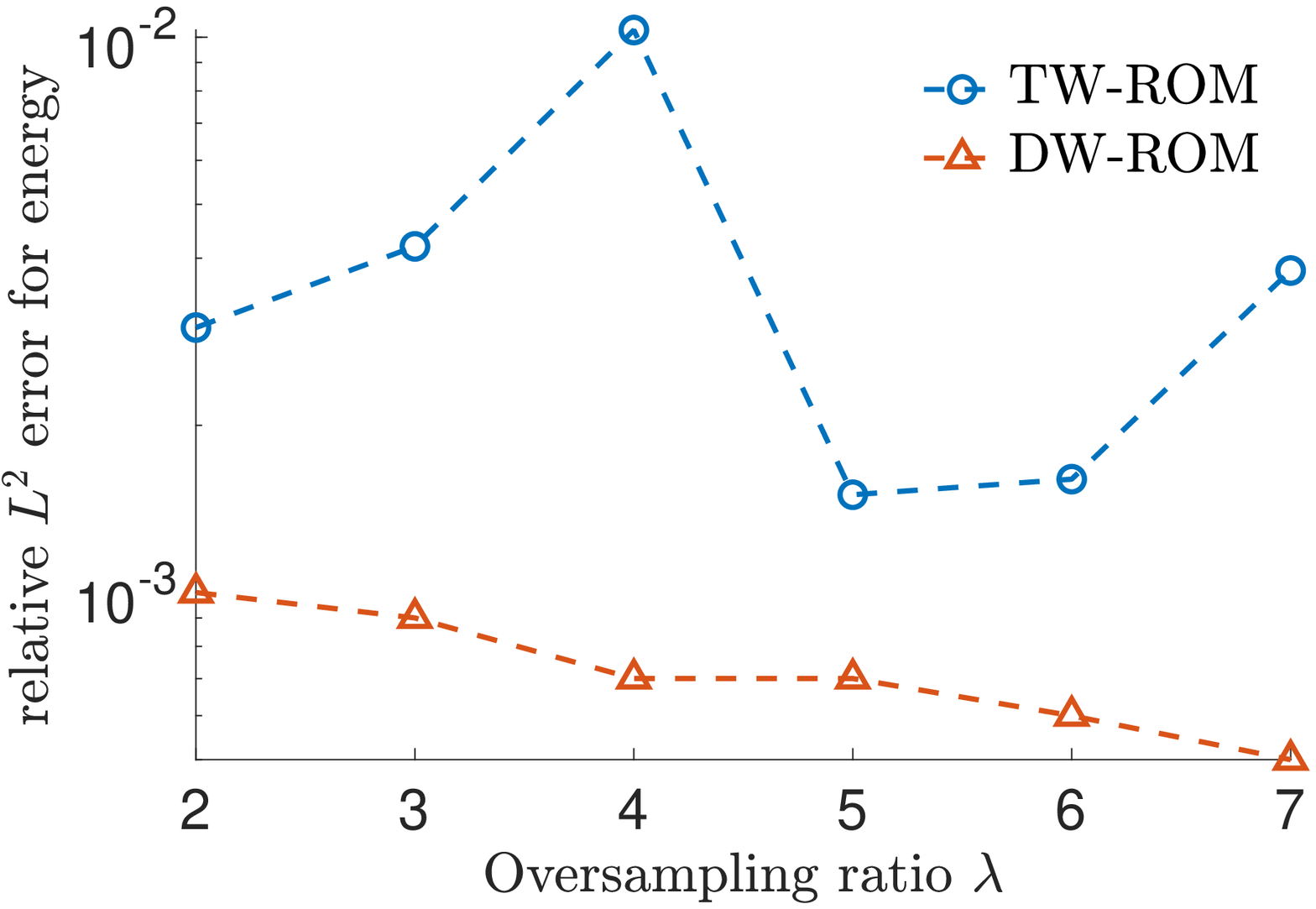}
\vspace{0.2in}
\includegraphics[width=0.45\linewidth]{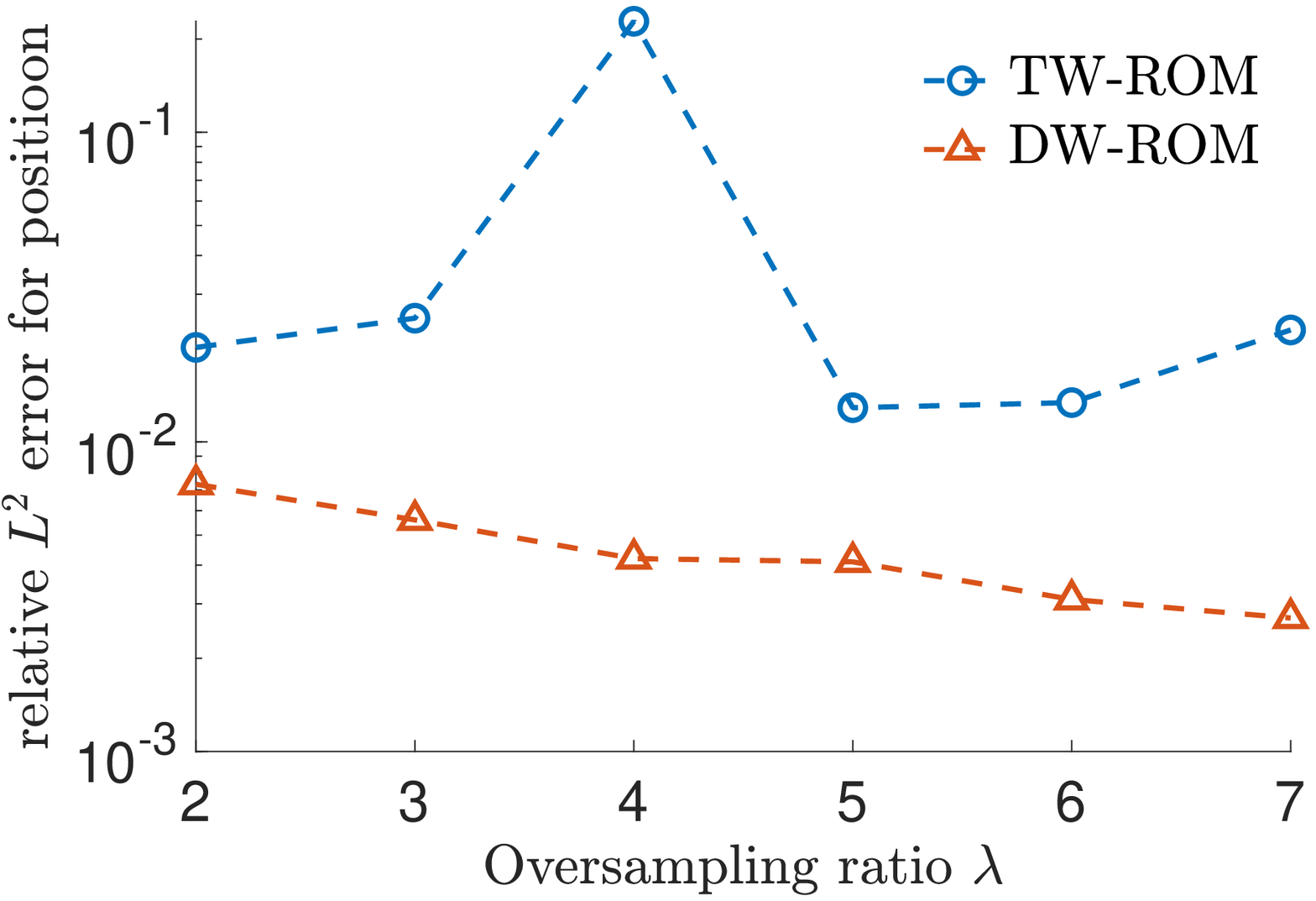}
\hspace{0.05\linewidth}
\includegraphics[width=0.45\linewidth]{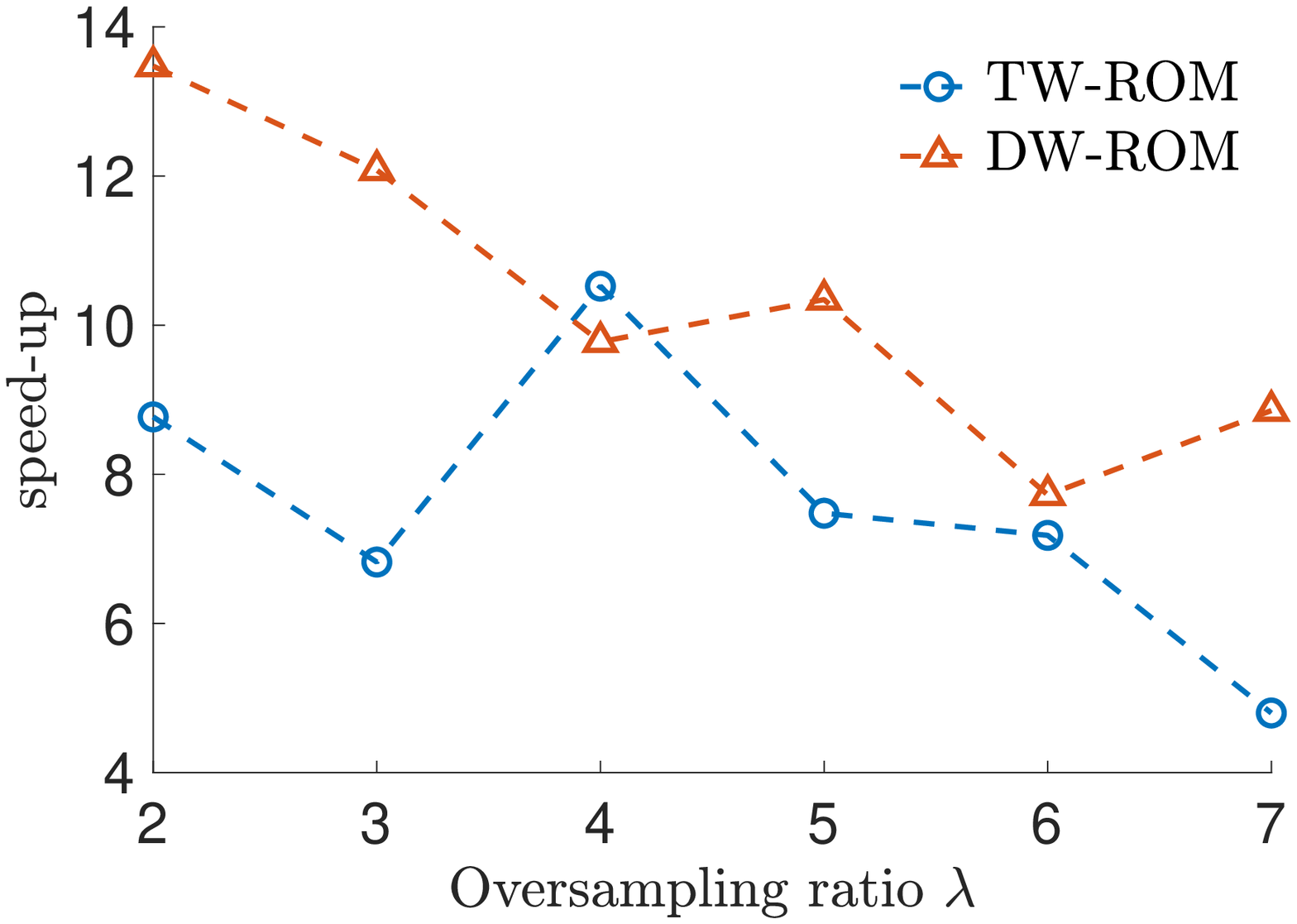}
\caption{ROM performance comparison for long-time simulation in 
Rayleigh--Taylor instability problem with varying oversampling ratio: 
relative $L^2$ error for velocity (top-left), 
relative $L^2$ error for energy (top-right), 
relative $L^2$ error for position (bottom-left), and 
speed-up (bottom-right). 
Using the penetration distance improves the solution accuracy 
and speed-up at all tested oversampling ratio.}
\label{fig:reuslts-oversampling}
\end{figure}

First, we take $\ntimestepWindow = 20$ for decomposing the indicator range by Algorithm~\ref{alg:decompose}.
The ROM simulation is performed at the unseen parameter $\param = 0.33$ 
with final time $\finalTimeROM(\param) = 3.6$ and varying oversampling ratio 
$\factorROMforceOneSample=\factorROMforceTvSample=\lambda \in [2,7]$ for hyper-reduction in the online phase. 
To avoid complicating the results, we fix all ROM parameters except the oversampling ratio. 
Figure~\ref{fig:reuslts-oversampling} shows the results in terms of solution accuracy and speed-up
with $\ntimestepWindow = 20$ and various oversampling ratio $\lambda$. 
The results of TW-ROM, i.e. using physical time as the indicator in Section~\ref{sec:tw}, 
and DW-ROM, i.e. using penetration distance as the indicator in Section~\ref{sec:dw}, 
are depicted in blue and red, respectively. 
It can be seen that using the penetration distance as the indicator does not only provide 
more accurate approximation in all the solution variables but also better speed-up in the ROM simulation 
compared to using the physical time. 
Moreover, when the penetration distance is used, we can observe a decaying trend in the solution error and speed-up 
with increasing oversampling ratio $\lambda$, which provides a tuneable solution accuracy 
against the computational expense. 
On the other hand, the solution accuracy does not necessarily improve 
with increasing oversampling ratio $\lambda$ when the physical time is used. 
At the oversampling ratio $\lambda = 7$, the error in velocity is around 48\% when the physical time is used as the indicator, 
while it is around 5\% when the penetration distance is used. 

\begin{figure}[ht!]
\centering
\includegraphics[width=0.96\linewidth]{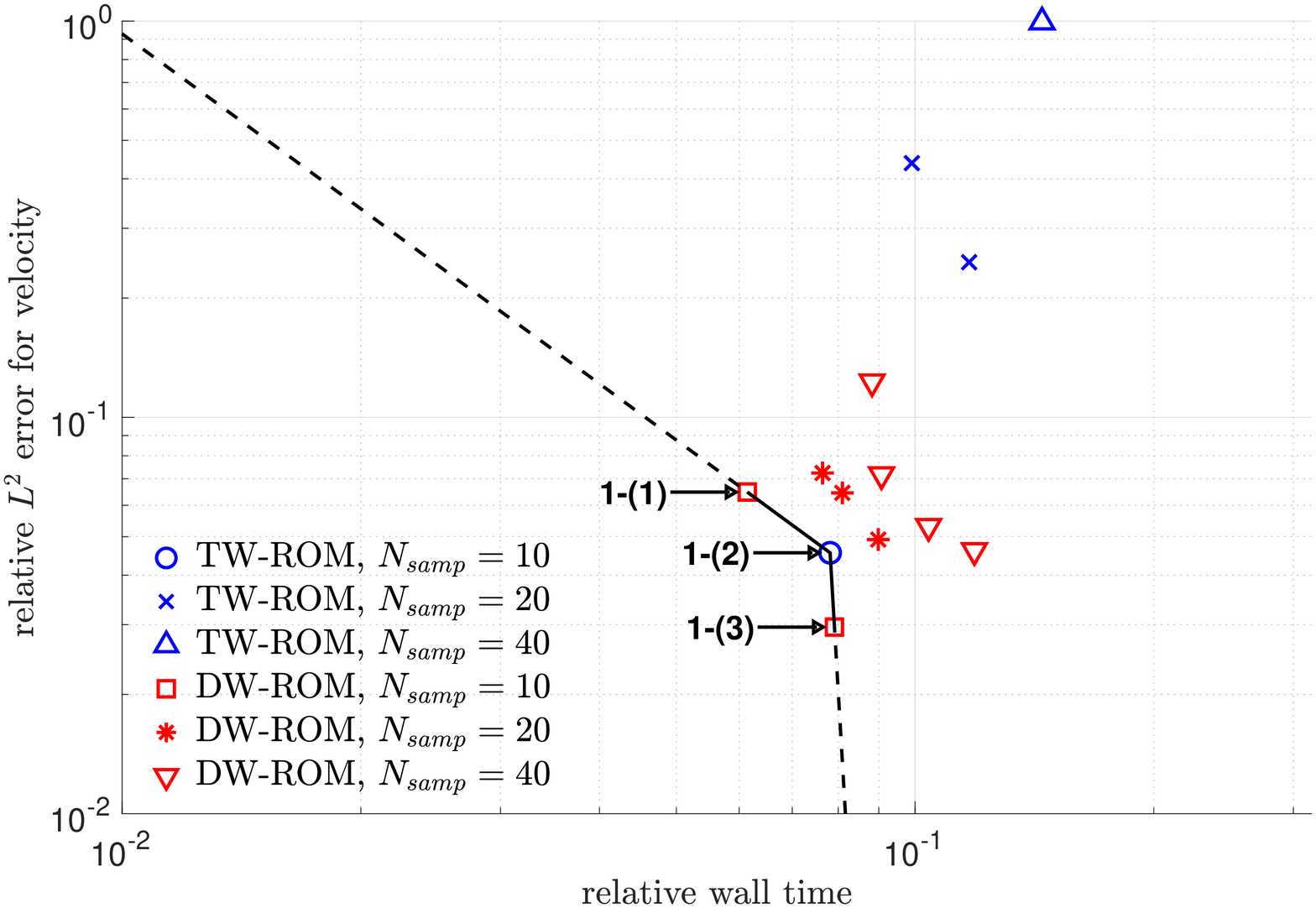}
\caption{ROM performance comparison for long-time simulation in 
Rayleigh--Taylor instability problem at $\param = 0.33$ with varying oversampling ratio. 
Relative $L^2$ error for velocity versus relative wall time for varying ROM parameters.
  }
\label{fig:pareto_af33}
\end{figure}

\begin{figure}[ht!]
\centering
\includegraphics[width=0.96\linewidth]{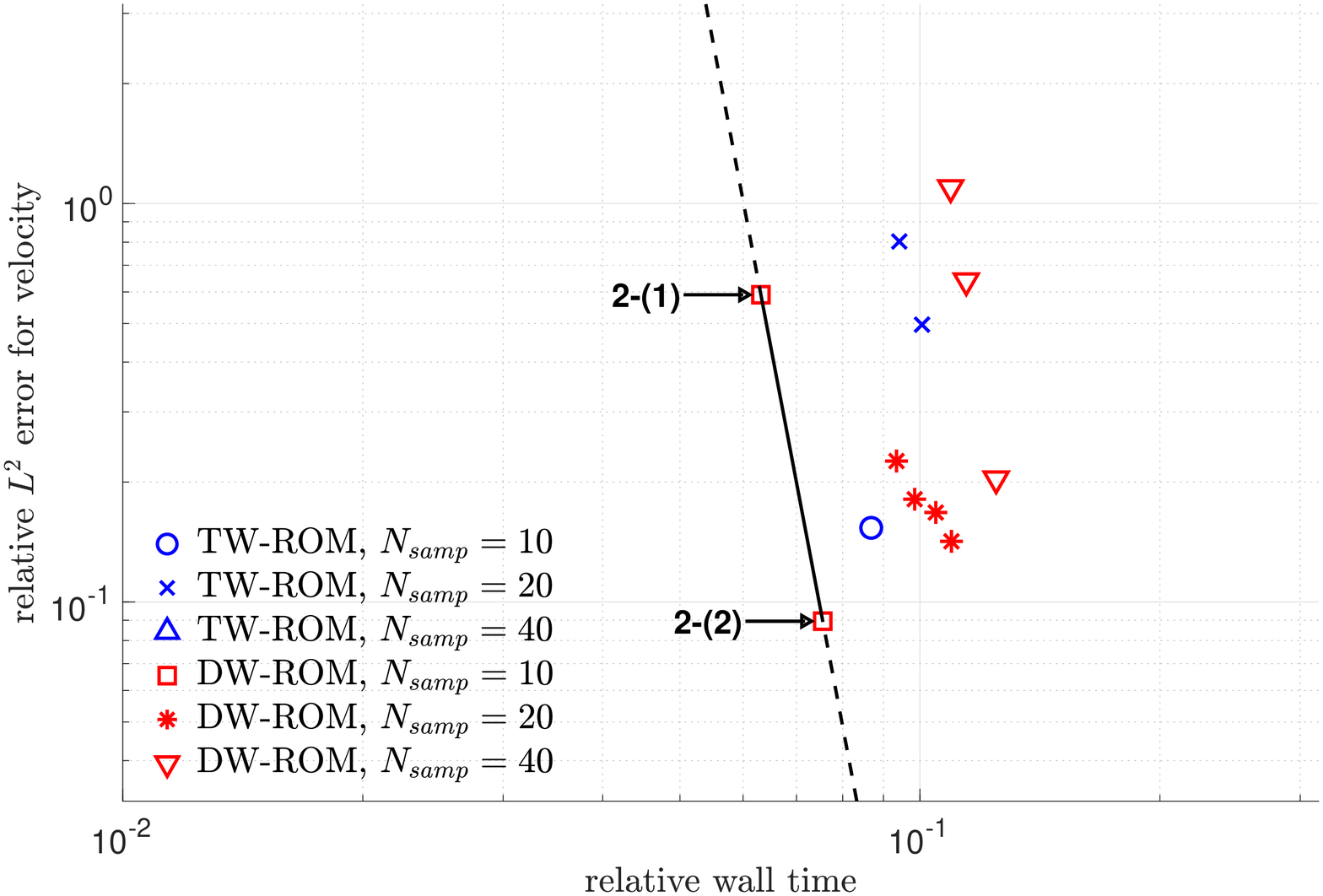}
\caption{ROM performance comparison for long-time simulation in 
Rayleigh--Taylor instability problem at $\param = 0.32$ with varying oversampling ratio. 
Relative $L^2$ error for velocity versus relative wall time for varying ROM parameters.
  }
\label{fig:pareto_af32}
\end{figure}

Finally, we provide a thorough comparison of different choices in ROM simulation. 
We use different combinations of indicator and $\ntimestepWindow \in \{10, 20, 40\}$ in Algorithm~\ref{alg:decompose}. 
Again, the ROM simulation is performed at the unseen parameter $\param = \{0.32,0.33\}$ 
with final time $\finalTimeROM(\param) = 3.6$ and varying oversampling ratio 
$\factorROMforceOneSample=\factorROMforceTvSample=\lambda \in [2,7]$ for hyper-reduction in the online phase. 
In Figure~\ref{fig:pareto_af33} and Figure~\ref{fig:pareto_af32}, we construct a Pareto front, which is characterized by
the ROM user-defined input values that minimize the competing objectives of
relative $L^2$ error for velocity and relative wall time, with all the combinations of ROM parameters 
at $\param = 0.33$ and $\param = 0.32$, respectively. 
It can be observed that using the penetration distance as the indicator performs better in general, 
and the error is less sensitive to the number of intermediate samples in a subinterval $\ntimestepWindow$. 
Meanwhile, the error increases dramatically with $\ntimestepWindow$ when the physical time is used as the indicator. 
Moreover, with $\param = 0.32$ and $\ntimestepWindow = 40$, 
TW-ROM is numerically unstable with any choice of oversampling ratio $\lambda$. 
An overall Pareto front that selects the ROM parameters that are Pareto-optimal across all groups
is also illustrated. Table~\ref{tab:pareto_labels}
reports the ROM user-defined input values that yielded the results on the
overall Pareto fronts. The results of this experiment also suggest that 
the penetration distance is a more reliable indicator than the physical time 
to produce a remarkable speed-up as well as accurate approximations 
with the appropriate use of hyper-reduction. 

\begin{table}[ht!]
\centering
\begin{tabular}{|c||c|c|c||c|c|}
\hline
Label & 1-(1) & 1-(2) & 1-(3) & 2-(1) & 2-(2)  \\
\hline
Decomposition & DW & TW & DW & DW & DW \\
$\ntimestepWindow$ & 10 & 10 & 10 & 10 & 10 \\ 
$\lambda$ & 2 & 2 & 4 & 2 & 4 \\ 
\hline
\end{tabular}
\caption{ROM user-defined input values yielding Pareto-optimal performance for
  long-time simulation in the Rayleigh--Taylor instability problem with varying oversampling ratio. 
  Figure~\ref{fig:pareto_af33}--\ref{fig:pareto_af32} provide labels.}
\label{tab:pareto_labels}
\end{table}

\section{Conclusion}\label{sec:conclusion}
In this paper, we have developed reduced order model techniques to accelerate 
the simulation of Rayleigh--Taylor instability in compressible gas dynamics. 
The Euler equation is used to model the compressible gas dynamics 
in a complex multimaterial setting. 
Moreover, we study the effects of Atwood number as a problem parameter on the physical quantities. 
Using the curvilinear finite element method for solving the governing equations in a moving Lagrangian frame, 
the numerical solutions of Rayleigh--Taylor instability are accurate but computationally expensive. 
We use POD for solution basis construction and SNS and DEIM for hyper-reduction. 
Furthermore, we introduce a general framework for the temporal domain partition and 
temporally-local reduced order model construction with varying Atwood number. 
We compare two approaches in this framework, namely decomposing the solution manifold 
by physical time and penetration distance respectively. 
Numerical results show that the penetration distance indicator provides better results than the physical time indicator, 
in terms of both speed-up and solution accuracy. 
In the future, we will carry out further investigation on reliable indicators for 
parameter-dependent temporal domain partition which lead to 
efficient reduced order models for advection-dominated problems. 
Although this work focused on a single parameter specific to the application, the general methodology
can in principle be extended to other parameters for other problems.

\section*{Acknowledgments}
This work was performed at Lawrence Livermore National Laboratory.
Lawrence Livermore National Laboratory is operated by Lawrence
Livermore National Security, LLC, for the U.S. Department of Energy,
National Nuclear Security Administration under Contract DE-AC52-07NA27344
and LLNL-JRNL-830909.

\section*{Disclaimer}
This document was prepared as an account of work sponsored by an agency of the
United States government.  Neither the United States government nor Lawrence
Livermore National Security, LLC, nor any of their employees makes any warranty,
expressed or implied, or assumes any legal liability or responsibility for the
accuracy, completeness, or usefulness of any information, apparatus, product, or
process disclosed, or represents that its use would not infringe privately owned
rights.  Reference herein to any specific commercial product, process, or
service by trade name, trademark, manufacturer, or otherwise does not
necessarily constitute or imply its endorsement, recommendation, or favoring by
the United States government or Lawrence Livermore National Security, LLC.  The
views and opinions of authors expressed herein do not necessarily state or
reflect those of the United States government or Lawrence Livermore National
Security, LLC, and shall not be used for advertising or product endorsement
purposes.

\appendix
\section{Command line options of LaghosROM}\label{sec:appendix} 
For the purpose of reproducible research, we present some examples of the command line options 
for simulating the Rayleigh--Taylor instability problem using LaghosROM. 
Due to rapid software development in the repository, the commands lines are subject to change. 
However, we try our best to present the command lines compatible with recent versions of different
dependent softwares and maintain a simple usage of the program. The following
command lines are compatible with the recent commits of 
the \texttt{master} branch of MFEM
\footnote{GitHub page, {\it https://github.com/mfem/mfem}, commit 39022bc}, 
the \texttt{master} branch of libROM
\footnote{GitHub page, {\it https://github.com/LLNL/libROM}, commit 6258824}, 
and the \texttt{rom} branch of Laghos
\footnote{GitHub page, {\it https://github.com/CEED/Laghos/tree/rom}, commit 26b3206}. 
In order to use the ROM capability of Laghos, a user has to 
navigate to the \texttt{rom} subdirectory. 

\subsection{Problem specification}\label{sec:cmd_prob}
First, we present the commands lines for the FOM user-defined input values 
reported in Section~\ref{sec:numerical} using the executable \texttt{laghos}. 
These command lines can also be used along with 
command line options for ROM user-defined input values. 
We use the long-time simulation for illustration. 
The following command line is used to specify the full-order Rayleigh-Taylor instability problem 
with the final time $\finalTime = 3.6$, refinement level of 4, and 
the same finite element polynomial degree of 2 and 1 for the kinematic and thermodynamic space respectively. 
\begin{verbatim}
./laghos -p 7 -m data/rt2D.mesh -s 7 -tf 3.6 -rs 4 -ok 2 -ot 1
\end{verbatim}
This basic command line is appended by different options in different stages, 
namely computing FOM reference solution, 
sampling snapshots, constructing reduced basis and operators, 
computing ROM solution, and evaluating error. 
For computing FOM reference solution at the testing Atwood number $\param = 0.33$
to evaluate the error of the ROM solution obtained in Section~\ref{sec:cmd_online}, one appends 
\begin{verbatim}
-writesol -af 0.33
\end{verbatim}

\subsection{Offline computations}\label{sec:cmd_offline}
Next, we present the command line options for the offline 
ROM construction in Rayleigh-Taylor instability problem, 
which involves sampling snapshots and constructing reduced basis. 
The options are appended to the command lines for
specifying the full-order problem in Section~\ref{sec:cmd_prob}.  
To build the ROM with one training Atwood number  
$\param_1 = 1/3$, and using penetration distance as the indicator
and $\ntimestepWindow = 20$ in Algorithm~\ref{alg:decompose}, 
one appends 
\begin{verbatim}
-offline -romsns -ef 0.9999 -sample-stages -nwinsamp 21 -loctype distance
\end{verbatim}
Here and thereafter, one drops \texttt{-loctype distance} in all the command lines 
for using physical time as the indicator instead. 

On the other hand, to build the ROM with two training Atwood number 
$\param_1 = 1/3$ and $\param_2 = 0.3$, one appends 
\begin{verbatim}
-offline -romsns -sample-stages -rpar 0 -sdim 200000 -loctype distance
-offline -romsns -sample-stages -rpar 1 -af 0.3 -sdim 200000 -loctype distance
\end{verbatim}
and then uses a separate executable \texttt{merge}
\begin{verbatim}
./merge -nset 2 -romsns -ef 0.9999 -nwinsamp 21 -loctype distance
\end{verbatim}

\subsection{Online computations}\label{sec:cmd_online}
Finally, we present the command line options for the online 
ROM simulation in Rayleigh-Taylor instability problem, 
which involves constructing reduced operators, 
computing ROM solution and evaluating error. 
The options are appended to the command lines for
specifying the full-order problem in Section~\ref{sec:cmd_prob}.  
Below we present the ROM commands for using SNS hyper-reduction 
with over-sampling ratio $\factorROMforceOneSample = \factorROMforceTvSample = 2$
on $\nwindow^{\text{off}} = 3863$ subintervals in the indicator range (from Algorithm~\ref{alg:decompose})
at the testing Atwood number $\param = 0.33$. 
For hyper-reduction preprocessing in the online phase, one appends
\begin{verbatim}
-online -romhrprep -romsns -sfacv 2 -sface 2 -nwin 3863 -af 0.33 -loctype distance
\end{verbatim}
For ROM simulation in the online phase, one appends 
\begin{verbatim}
-online -romhr -romsns -sfacv 2 -sface 2 -nwin 3863 -af 0.33 -loctype distance
\end{verbatim}
Finally, for solution postprocessing and calculating the relative error
between the FOM solution and the ROM solution, one appends 
\begin{verbatim}
-restore -soldiff -romsns -nwin 3863 -af 0.33 -loctype distance
\end{verbatim}

\bibliographystyle{unsrt}
\bibliography{references}

\end{document}

%% file: commands.tex
\newcommand{\bmat}[1]{\begin{bmatrix}#1\end{bmatrix}} 
\newcommand{\pmat}[1]{\begin{pmatrix}#1\end{pmatrix}} 

\newcommand{\dummyMatSymbol}{A}
\newcommand{\dummyElemSymbol}{a}
\newcommand{\dummyVecSymbol}{\dummyElemSymbol}

\newcommand{\dimensionSymbol}{d}
\newcommand{\solDomainSymbol}{\Omega}
\newcommand{\gradientSymbol}{\nabla}
\newcommand{\energySymbol}{e}
\newcommand{\velocitySymbol}{v}
\newcommand{\positionSymbol}{x}
\newcommand{\densitySymbol}{\rho}
\newcommand{\pressureSymbol}{p}
\newcommand{\gravitySymbol}{g}
\newcommand{\stressSymbol}{\sigma}
\newcommand{\artificialStressSymbol}{\stressSymbol_a}
\newcommand{\adiabaticIndexSymbol}{\gamma}
\newcommand{\normalSymbol}{n}
\newcommand{\identitySymbol}{I}

\newcommand{\initialSymbol}[1]{\widetilde{#1}}
\newcommand{\initialDomain}{\initialSymbol{\solDomainSymbol}}
\newcommand{\initialPosition}{\initialSymbol{\positionSymbol}}

\newcommand{\snapshotSymb}{E}

\newcommand{\subspaceSymbol}{\mathcal{S}}
\newcommand{\timeSymbol}{t}

\newcommand{\massMatSymbol}{M}

\newcommand{\oneSymbol}{1}
\newcommand{\forceSymbol}{F}

\newcommand{\kinematicSymbol}{V}
\newcommand{\thermodynamicSymbol}{E}
\newcommand{\feSpace}[1]{\mathcal{#1}}
\newcommand{\domain}[1]{\mathsf{#1}}
\newcommand{\paramSymbol}{A}
\newcommand{\mapto}{\rightarrow}
\newcommand{\sizeFOMsymbol}{N}
\newcommand{\sizeROMsymbol}{n}
\newcommand{\sizeSnapshot}{N_{\snapshotSymb}}

\newcommand{\nsmall}{\sizeROMsymbol}
\newcommand{\Span}[1]{\text{Span}\{#1\}}
\newcommand{\RR}[1]{\mathbb{R}^{#1}}

\newcommand{\timeIndex}{n}
\newcommand{\basisIndex}{i}
\newcommand{\paramIndex}{k}
\newcommand{\timeWindow}{\mathcal{T}}
\newcommand{\windowIndex}{j}
\newcommand{\paramDomainSymbol}{D}
\newcommand{\finalSymbol}{f}

\newcommand{\ROMBasisSymbol}{\mathbf{\Phi}}
\newcommand{\ROMBasisVecSymbol}{\mathbf{\phi}}
\newcommand{\ROMSamplingMatSymbol}{\mathbf{Z}}
\newcommand{\ROMSolutionSymbol}[1]{\widehat{#1}}
\newcommand{\ProjectionSymbol}{\mathcal{P}}

\newcommand{\relErrorSymbol}{\mathcal{\epsilon}}

\newcommand{\fullNorm}[1]{{\left\vert\kern-0.25ex\left\vert\kern-0.25ex\left\vert #1 
    \right\vert\kern-0.25ex\right\vert\kern-0.25ex\right\vert}}
\newcommand{\stateSymbol}{y}

\newcommand{\forceOne}{\forceSys^{\mathfrak{1}}}
\newcommand{\forceTv}{\forceSys^{t\velocitySymbol}}
\newcommand{\avgforceTv}{\bar{\forceSys}^{t\velocitySymbol}}
\newcommand{\forceOnek}[1]{\forceOne_{#1}}
\newcommand{\forceTvk}[1]{\forceTv_{#1}}
\newcommand{\avgforceTvk}[1]{\avgforceTv_{#1}}
\newcommand{\forceOneApprox}{\tilde{\forceSys}^{\mathfrak{1}}}
\newcommand{\forceTvApprox}{\tilde{\forceSys}^{t\velocitySymbol}}
\newcommand{\avgforceTvApprox}{\bar{\tilde{\forceSys}}^{t\velocitySymbol}}
\newcommand{\forceOneApproxk}[1]{\forceOneApprox_{#1}}
\newcommand{\forceTvApproxk}[1]{\forceTvApprox_{#1}}
\newcommand{\avgforceTvApproxk}[1]{\avgforceTvApprox_{#1}}

\newcommand{\energySubspace}{\subspaceSymbol_{\energySymbol}}
\newcommand{\velocitySubspace}{\subspaceSymbol_{\velocitySymbol}}
\newcommand{\positionSubspace}{\subspaceSymbol_{\positionSymbol}}
\newcommand{\forceOneSubspace}{\subspaceSymbol_{\forceOne}}
\newcommand{\forceTvSubspace}{\subspaceSymbol_{\forceTv}}

\newcommand{\paramDomain}{\domain{\paramDomainSymbol}}

\newcommand{\energy}{\boldsymbol{\energySymbol}}
\newcommand{\velocity}{\boldsymbol{\velocitySymbol}}
\newcommand{\position}{\boldsymbol{\positionSymbol}}
\newcommand{\state}{\boldsymbol{\stateSymbol}}

\newcommand{\velocityApprox}{\tilde{\velocity}}
\newcommand{\energyApprox}{\tilde{\energy}}
\newcommand{\positionApprox}{\tilde{\position}}
\newcommand{\stateApprox}{\tilde{\state}}

\newcommand{\relError}{\relErrorSymbol}

\newcommand{\massMat}{\boldsymbol{\massMatSymbol}}

\newcommand{\param}{\paramSymbol}
\newcommand{\oneVec}{\boldsymbol{\oneSymbol}}
\newcommand{\kinematicFE}{\feSpace{\kinematicSymbol}}
\newcommand{\thermodynamicFE}{\feSpace{\thermodynamicSymbol}}
\newcommand{\kinematicMassMat}{\massMat_{\kinematicFE}}
\newcommand{\thermodynamicMassMat}{\massMat_{\thermodynamicFE}}

\newcommand{\forceMat}{\boldsymbol{\forceSymbol}}

\newcommand{\forceSys}{\mathsf{\forceSymbol}}
\newcommand{\sizeKinematicFE}{\sizeFOMsymbol_{\kinematicFE}}
\newcommand{\sizeThermodynamicFE}{\sizeFOMsymbol_{\thermodynamicFE}}
\newcommand{\sizeROMforceOne}{\sizeROMsymbol_{\forceOne}}
\newcommand{\sizeROMforceTv}{\sizeROMsymbol_{\forceTv}}
\newcommand{\samplingOp}[1]{\breve{#1}}
\newcommand{\sizeROMforceOneSample}{\sizeROMsymbol_{\samplingOp{\forceOne}}}
\newcommand{\sizeROMforceTvSample}{\sizeROMsymbol_{\samplingOp{\forceTv}}}
\newcommand{\sizeROMvelocity}{\sizeROMsymbol_{\velocitySymbol}}
\newcommand{\sizeROMenergy}{\sizeROMsymbol_{\energySymbol}}
\newcommand{\sizeROMposition}{\sizeROMsymbol_{\positionSymbol}}

\newcommand{\factorROMforceOneSample}{\lambda_{\samplingOp{\forceOne}}}
\newcommand{\factorROMforceTvSample}{\lambda_{\samplingOp{\forceTv}}}
\newcommand{\sizeWholeFE}{\sizeFOMsymbol}
\newcommand{\sizeWholeROM}{\sizeROMsymbol}
\newcommand{\timek}[1]{\timeSymbol_{#1}}
\newcommand{\finalTime}{\timek{\finalSymbol}}
\newcommand{\ntimestep}{N_{\timeSymbol}}
\newcommand{\timestep}{\Delta\timeSymbol}
\newcommand{\timestepk}[1]{\timestep_{#1}}

\newcommand{\avgvelocityt}[1]{\bar{\velocity}_{#1}}
\newcommand{\statet}[1]{\state_{#1}}
\newcommand{\energyt}[1]{\energy_{#1}}
\newcommand{\velocityt}[1]{\velocity_{#1}}
\newcommand{\positiont}[1]{\position_{#1}}

\newcommand{\avgvelocityApproxt}[1]{\bar{\velocityApprox}_{#1}}
\newcommand{\stateApproxt}[1]{\stateApprox_{#1}}
\newcommand{\energyApproxt}[1]{\energyApprox_{#1}}
\newcommand{\velocityApproxt}[1]{\velocityApprox_{#1}}
\newcommand{\positionApproxt}[1]{\positionApprox_{#1}}
\newcommand{\ROMenergyt}[1]{\ROMenergy_{#1}}
\newcommand{\ROMvelocityt}[1]{\ROMvelocity_{#1}}
\newcommand{\ROMpositiont}[1]{\ROMposition_{#1}}

\newcommand{\nat}[1]{\mathbb{N}(#1)}
\newcommand{\identityMat}[1]{\boldsymbol{\identitySymbol}_{#1}}

\newcommand{\velocityBasis}{\ROMBasisSymbol_{\velocitySymbol}}
\newcommand{\energyBasis}{\ROMBasisSymbol_{\energySymbol}}
\newcommand{\positionBasis}{\ROMBasisSymbol_{\positionSymbol}}
\newcommand{\velocityBasisVec}{\ROMBasisVecSymbol_{\velocitySymbol}}
\newcommand{\energyBasisVec}{\ROMBasisVecSymbol_{\energySymbol}}
\newcommand{\positionBasisVec}{\ROMBasisVecSymbol_{\positionSymbol}}
\newcommand{\velocityBasisVeck}[1]{\velocityBasisVec^{#1}}
\newcommand{\energyBasisVeck}[1]{\energyBasisVec^{#1}}
\newcommand{\positionBasisVeck}[1]{\positionBasisVec^{#1}}
\newcommand{\forceSysEnergy}{\forceSys_{\energySymbol}}
\newcommand{\forceSysVelocity}{\forceSys_{\velocitySymbol}}
\newcommand{\forceSysPosition}{\forceSys_{\positionSymbol}}
\newcommand{\forceTvBasis}{\ROMBasisSymbol_{\forceTv}}
\newcommand{\forceOneBasis}{\ROMBasisSymbol_{\forceOne}}
\newcommand{\forceTvBasisVec}{\ROMBasisVecSymbol_{\forceTv}}
\newcommand{\forceOneBasisVec}{\ROMBasisVecSymbol_{\forceOne}}
\newcommand{\forceTvBasisVeck}[1]{\forceTvBasisVec^{#1}}
\newcommand{\forceOneBasisVeck}[1]{\forceOneBasisVec^{#1}}
\newcommand{\forceTvSamplingMat}{\ROMSamplingMatSymbol_{\forceTv}}
\newcommand{\forceOneSamplingMat}{\ROMSamplingMatSymbol_{\forceOne}}

\newcommand{\ROMstate}{\ROMSolutionSymbol{\state}}
\newcommand{\ROMenergy}{\ROMSolutionSymbol{\energy}}
\newcommand{\ROMvelocity}{\ROMSolutionSymbol{\velocity}}
\newcommand{\ROMposition}{\ROMSolutionSymbol{\position}}
\newcommand{\ROMforceOne}{\ROMSolutionSymbol{\forceOne}}
\newcommand{\ROMforceTv}{\ROMSolutionSymbol{\forceTv}}
\newcommand{\ROMforceSys}{\ROMSolutionSymbol{\forceSys}}
\newcommand{\ROMforceSysEnergy}{\ROMforceSys_{\energySymbol}}
\newcommand{\ROMforceSysVelocity}{\ROMforceSys_{\velocitySymbol}}
\newcommand{\ROMforceSysPosition}{\ROMforceSys_{\positionSymbol}}
\newcommand{\offsetSymbol}{\text{os}}
\newcommand{\energyOS}{\energy_{\offsetSymbol}}
\newcommand{\velocityOS}{\velocity_{\offsetSymbol}}
\newcommand{\positionOS}{\position_{\offsetSymbol}}

\newcommand{\timekROM}[1]{\widetilde{\timeSymbol}_{#1}}
\newcommand{\finalTimeROM}{\timekROM{\finalSymbol}}
\newcommand{\ntimestepROM}{\widetilde{N}_{\timeSymbol}}

\newcommand{\forceOneObliqProjMat}{\ProjectionSymbol_{\forceOne}}

\newcommand{\windowIndicator}{\Psi}
\newcommand{\snapshotGroup}{\mathcal{G}}
\newcommand{\windowk}[1]{T_{#1}}
\newcommand{\nwindow}{N_w}

\newcommand{\sizeROMforceOneWindow}[1]{\sizeROMsymbol_{\forceOne}^{#1}}
\newcommand{\sizeROMforceTvWindow}[1]{\sizeROMsymbol_{\forceTv}^{#1}}
\newcommand{\sizeROMforceOneSampleWindow}[1]{\sizeROMsymbol_{\samplingOp{\forceOne}}^{#1}}
\newcommand{\sizeROMforceTvSampleWindow}[1]{\sizeROMsymbol_{\samplingOp{\forceTv}}^{#1}}
\newcommand{\energyOSWindow}[1]{\energy_{\offsetSymbol}^{#1}}
\newcommand{\velocityOSWindow}[1]{\velocity_{\offsetSymbol}^{#1}}
\newcommand{\positionOSWindow}[1]{\position_{\offsetSymbol}^{#1}}
\newcommand{\velocityBasisWindow}[1]{\ROMBasisSymbol_{\velocitySymbol}^{#1}}
\newcommand{\energyBasisWindow}[1]{\ROMBasisSymbol_{\energySymbol}^{#1}}
\newcommand{\positionBasisWindow}[1]{\ROMBasisSymbol_{\positionSymbol}^{#1}}
\newcommand{\forceOneBasisWindow}[1]{\ROMBasisSymbol_{\forceOne}^{#1}}
\newcommand{\forceTvBasisWindow}[1]{\ROMBasisSymbol_{\forceTv}^{#1}}
\newcommand{\ROMenergyWindow}[1]{\ROMenergy^{#1}}
\newcommand{\ROMvelocityWindow}[1]{\ROMvelocity^{#1}}
\newcommand{\ROMpositionWindow}[1]{\ROMposition^{#1}}

\newcommand{\forceTvSamplingMatWindow}[1]{\forceTvSamplingMat^{#1}}
\newcommand{\forceOneSamplingMatWindow}[1]{\forceOneSamplingMat^{#1}}
\newcommand{\ntimestepWindow}{N_{\text{sample}}}
\newcommand{\timeIndexWindow}[1]{\timeIndex_{#1}}

\newcommand{\relerrorVelocityt}[1]{\relError_{\velocitySymbol,#1}}
\newcommand{\relerrorEnergyt}[1]{\relError_{\energySymbol,#1}}
\newcommand{\relerrorPositiont}[1]{\relError_{\positionSymbol,#1}}

\newcommand{\snapshots}{\boldsymbol \snapshotSymb}
\newcommand{\rkStage}{r}

\newcommand{\dummyMat}{\boldsymbol{\dummyMatSymbol}}
\newcommand{\dummyVec}{\boldsymbol{\dummyVecSymbol}}
\newcommand{\leftSingularMat}{\boldsymbol{U}}
\newcommand{\leftSingularVec}{\boldsymbol{u}}
\newcommand{\leftSingularVeck}[1]{\leftSingularVec_{#1}}
\newcommand{\rightSingularMat}{\boldsymbol{V}}
\newcommand{\singularValueMat}{\boldsymbol{\Sigma}}
\newcommand{\singularValue}{\sigma}
\newcommand{\singularValueThreshold}{\delta_{\singularValue}}
\newcommand{\nparam}{\nsmall_\paramSymbol}

\newcommand{\conditionnumber}{\kappa}
\newcommand{\unitvecSymbol}{\mathfrak{e}}
\newcommand{\unitveck}[1]{\unitvecSymbol_{#1}}
\newcommand{\tuningparam}{\eta}
\DeclareMathOperator*{\argmin}{arg\,min}